\documentclass[11pt]{article}

\usepackage[utf8]{inputenc}
\usepackage[T1]{fontenc}
\usepackage{geometry}
\usepackage{mathpazo}
\usepackage{amsmath,amsthm,mathtools,amssymb}
\usepackage{graphicx}
\usepackage{booktabs}
\usepackage{subcaption}
\usepackage{enumitem}
\usepackage{csquotes}
\usepackage[
    backend=biber,
    style=apa6,
    sorting=nyt
]{biblatex}
\usepackage{titlesec}
\usepackage[hidelinks]{hyperref}

\newtheorem{theorem}{Theorem}[section]
\newtheorem{proposition}[theorem]{Proposition}
\newtheorem{lemma}[theorem]{Lemma}

\newtheorem{assumption}[theorem]{Assumption}

\theoremstyle{definition}
\newtheorem{definition}[theorem]{Definition}
\theoremstyle{remark}
\newtheorem{remark}[theorem]{Remark}
\newtheorem{example}[theorem]{Example}

\DeclarePairedDelimiter{\abs}{\lvert}{\rvert}
\DeclarePairedDelimiter{\norm}{\lVert}{\rVert}
\DeclarePairedDelimiter{\sprod}{\langle}{\rangle}

\numberwithin{equation}{section}

\makeatletter
\renewcommand{\l@section}{\@dottedtocline{1}{0em}{1.8em}}
\renewcommand{\l@subsection}{\@dottedtocline{2}{1.8em}{2.8em}}
\renewcommand{\l@subsubsection}{\@dottedtocline{3}{4.6em}{3.7em}}
\makeatother

\titleformat{\section}
  {\normalfont\Large\bfseries}
  {\thesection}{1em}{}
\titleformat{\subsection}
  {\normalfont\large\bfseries}
  {\thesubsection}{1em}{}
\titleformat{\subsubsection}
  {\normalfont\normalsize\bfseries}
  {\thesubsubsection}{1em}{}

\begin{document}

\title{An Euler scheme for backward stochastic differential equations via the Wiener chaos decomposition}
\author{ Pere Diaz-Lozano \\\vspace{-0.2em}
    \small Department of Mathematics \\\vspace{-0.2em}
   \small University of Oslo \\\vspace{-0.2em}
  \small \texttt{peredl@math.uio.no} \and Giulia Di Nunno \\\vspace{-0.2em}
    \small Department of Mathematics \\\vspace{-0.2em}
   \small University of Oslo \\\vspace{-0.2em}
  \small \texttt{giulian@math.uio.no} 
  }

\date{}
\maketitle

\begin{abstract}
The Euler scheme is a standard time discretization for BSDEs, but its implementation hinges on approximating conditional expectations and the associated martingale terms at each time step. We propose an implementation based on the Wiener chaos decomposition to approximate these quantities. In contrast to many numerical schemes that rely on a finite-dimensional Markovian representation, our approach accommodates arbitrary $\mathcal{F}_T$-measurable square-integrable terminal conditions. We provide a comprehensive convergence analysis under additional Malliavin regularity assumptions and illustrate the method on several numerical examples, including genuinely non-Markovian problems arising, for instance, in the pricing and hedging of contingent claims under rough-volatility models.
\end{abstract}

\noindent\textbf{Keywords:} Backward stochastic differential equations, non-Markovian BSDEs, numerical methods, Wiener chaos decomposition, Malliavin calculus

\section{Introduction}

We are interested in the numerical approximation of the solution of a backward stochastic differential equation (BSDE for short) of the following form:
\begin{flalign}\label{BSDE}
    Y_{t} = \xi + \int_{t}^{T} f(s, Y_{s}, Z_{s}) ds - \int_{t}^{T} Z_{s} \cdot dB_{s}, \quad t \in [0,T],
\end{flalign}
where $B$ is a $d$-dimensional standard Brownian motion, $\mathbb{F} = (\mathcal{F}_{t})_{t \in [0,T]}$ is its natural and augmented filtration, $\xi$ is an $\mathbb{R}$-valued $\mathcal{F}_{T}$-measurable random variable, and the generator $f \colon [0,T] \times \Omega \times \mathbb{R} \times \mathbb{R}^{d} \to \mathbb{R}$ is a measurable $\mathbb{F}$-adapted random field. Informally, a solution to this equation is a pair of adapted processes $(Y,Z) = (Y_{t}, Z_{t})_{t \in [0,T]}$ satisfying equation (\ref{BSDE}). BSDEs were introduced by \textcite{BismutJean_Michel1973Ccfi} in the case of a linear generator, and were later generalized by \textcite{PARDOUX199055} to the non-linear case. Several applications of BSDEs in finance can be found, for example, in \textcite{ElKarouiN.1997BSDE}, where they are used for pricing and hedging financial products. We also mention their connection with dynamic risk measures; see \textcite{RosazzaGianinEmanuela2006Rmvg}.

There is extensive research on the numerical approximation of \eqref{BSDE}, with many methods relying on an Euler-type discretization over a grid
$\pi \coloneqq \{0=t_0<\cdots<t_m=T\}$. A standard backward Euler scheme is defined by $Y_{t_m}^\pi=\xi$ and, for $i=m,\dots,1$ and $t\in[t_{i-1},t_i)$,
\begin{flalign}\label{Euler scheme i}
Y_t^\pi
&=
\mathbb E\bigl(F_i^\pi\mid\mathcal F_t\bigr),
\qquad
\sprod{Z}_i^\pi
=
\mathbb E\left(
\frac{1}{\Delta_{i+1}}
\int_{t_i}^{t_{i+1}} Z_s^\pi\,ds
\,\middle|\,
\mathcal F_{t_i}
\right),
\quad i<m.
\end{flalign}
Here $\Delta_i\coloneqq t_i-t_{i-1}$, $\sprod{Z}_m^\pi=0$, and
\begin{flalign}\label{Euler scheme ii}
F_i^\pi
&=
Y_{t_i}^\pi
+
\Delta_i
f\bigl(t_i,Y_{t_i}^\pi,\sprod{Z}_i^\pi\bigr).
\end{flalign}
The process $(Z_s^\pi)_{s\in[t_{i-1},t_i)}$ is obtained from the martingale representation of $F_i^\pi$. We refer to \textcite{BouchardBruno2004DaaM} and \textcite{zhang_bsde_method} for the analysis of such schemes; related ideas already appear in \textcite{Chevance_1997} for generators independent of $Z$ and in \textcite{Bally1997} for random time meshes.

Much of the literature studies this scheme in a finite-dimensional Markovian
setting. More precisely, suppose that, for each
$(t,x)\in[0,T]\times\mathbb R^{d_X}$, there exists an
$\mathbb R^{d_X}$-valued process $(X_s^{t,x})_{s\in[t,T]}$ adapted to the
filtration generated by the Brownian increments after time $t$, satisfying
$X_t^{t,x}=x$ and the flow property
\begin{gather*}
X_r^{t,x}
=
X_r^{s,X_s^{t,x}},
\qquad
0\leq t\leq s\leq r\leq T,
\qquad
\mathbb P\text{-a.s.}
\end{gather*}
Setting $X=X^{0,x_0}$, one assumes that
$\xi=\varphi(X_T)$ and that the generator is of the form
$f(t,\omega,y,z)=\widetilde f(t,X_t(\omega),y,z)$, for deterministic
functions $\varphi$ and $\widetilde f$. Under the standard assumptions, the
discrete solution can then be represented through deterministic functions of
the forward state,
\begin{gather}\label{eq:deterministic-functions-numerics-BSDE}
Y_{t_i}^\pi=u_i(X_{t_i}),
\qquad
\sprod{Z}_i^\pi=v_i(X_{t_i}),
\end{gather}
so that the numerical problem reduces to approximating $u_i$ and $v_i$. Such a finite-dimensional representation is not always available. This occurs, for instance, in rough-volatility models driven by Volterra processes, as illustrated by the pricing and hedging problems considered in Section~\ref{section: numerical examples}; see also
\textcite{bayer_2022,FukasawaHorvathTankov2023Hedging}.

In particular, the classical finite-dimensional Markovian setting considered in much of the literature is obtained when the family is generated by a forward SDE,
\begin{flalign}\label{forward SDE}
X_s^{t,x}
=
x
+
\int_t^s b(r,X_r^{t,x})\,dr
+
\int_t^s \sigma(r,X_r^{t,x})\,dB_r,
\qquad
t\leq s\leq T.
\end{flalign}
Under the usual conditions ensuring strong existence and pathwise uniqueness, the resulting family satisfies the flow property. When the forward SDE cannot be simulated exactly, it may instead be approximated by a standard time discretization method, such as Euler--Maruyama. In \textcite{BouchardBruno2004DaaM}, the convergence of the backward Euler scheme was proved to be of order $1/2$ under standard Lipschitz assumptions in this setting. In \textcite{zhang_bsde_method}, the same order was obtained when the terminal condition is a Lipschitz functional of the whole path of $X$. A notable exception is \textcite{HuYaozhong2011MCFB}, where the order $1/2$ is established for general terminal conditions and random generators under appropriate Malliavin regularity assumptions.

Turning to implementation, the main difficulty in \eqref{Euler scheme i}--\eqref{Euler scheme ii} is the computation of the conditional expectations and martingale terms at each time step. Within the Markovian framework described above, these quantities can be approximated using least-squares regression 
\parencite{BouchardBruno2004DaaM,GobetEmmanuel2005ARMC,LemorJean-Philippe2006RoCo,One_step_malliavin}, Malliavin weights \parencite{BouchardBruno2004OtMa,GOBETEMMANUEL2016Aobs}, or cubature techniques \parencite{CrisanD.2012Sbsd}. More recently, neural networks have been used to approximate the corresponding deterministic functions in \eqref{eq:deterministic-functions-numerics-BSDE}, leading to approaches with strong performance when the forward state $X$ is high-dimensional; see \textcite{HanJiequn2018Shpd,HanJiequn2020Cotd,HureCome2020Dbsf,dinunno2024deepoperatorbsde}. These neural-network approaches can also be interpreted as non-linear least-squares Monte Carlo. However, a difficulty with regression-based implementations is that an approximation error in $Y$ is introduced directly at every step of the backward recursion. As the time grid is refined, the number of such errors therefore increases, and it is not immediate that their cumulative contribution vanishes as $|\pi|\to0$; see \textcite{HureCome2020Dbsf,GermainPhamWarin2022Approximation}.

For a comprehensive survey of numerical methods for BSDEs, see \textcite{ChessariJared2023Nmfb}.

Instead of introducing a forward process and approximating deterministic functions of its state, we directly expand each $F_i^\pi$ in Wiener chaos. More precisely, for $i=1,\dots,m$, one can write
\begin{gather*}
F_i^\pi
=
\sum_{k\geq0}
\sum_{|a|=k}
d_a^i
\prod_{j\geq1}
H_{a_j}
\left(
\int_0^{t_i}h_j^i(s)\cdot dB_s
\right),
\end{gather*}
where $H_k$ denotes the Hermite polynomial of degree $k$, $(h_j^i)_{j\geq1}$ is an orthonormal basis of $L^2([0,t_i];\mathbb R^d)$, $|a|=\sum_{j\geq1}a_j$ for any finitely supported sequence of non-negative integers $a=(a_j)_{j\geq1}$, and $(d_a^i)_a$ denotes the corresponding family of chaos coefficients.

Considering a finite-dimensional truncation of the form
\begin{gather*}
F_i^\pi
\approx
\sum_{k=0}^{P}
\sum_{|a|=k}
d_a^i
\prod_{j=1}^{M(i)}
H_{a_j}
\left(
\int_0^{t_i}h_j^i(s)\cdot dB_s
\right),
\qquad
a=(a_1,\dots,a_{M(i)}),
\end{gather*}
and approximating the coefficients $(d_a^i)_a$ by Monte Carlo yields an implementable numerical scheme. For a suitable choice of the functions $(h_j^i)_{1\leq j\leq M(i)}$, both the conditional expectations and the martingale terms appearing in \eqref{Euler scheme i}--\eqref{Euler scheme ii} can be computed explicitly from this representation. This construction does not rely on the finite-dimensional Markovian representation described above and can therefore handle general square-integrable terminal conditions.

Among the existing numerical schemes, the approach of \textcite{BriandPhilippe2014SOBB} is particularly relevant to the present work. It also accommodates general square-integrable terminal conditions and uses truncated Wiener chaos decompositions to compute the conditional expectations and martingale terms arising in successive Picard iterations. Our approach differs in several respects. In particular, our convergence analysis explicitly accounts for all the approximation parameters entering the implementation, including the discretization of the time integrals. Moreover, the generic constant in our estimate is independent of the approximation parameters, and our regularity assumptions are of fixed order, whereas both the constants and the required regularity in \textcite{BriandPhilippe2014SOBB} depend on the Picard and chaos truncation levels. From a computational perspective, our scheme avoids rerunning the full Picard iteration when additional samples are required and is less demanding in terms of computation time and memory usage, as illustrated in Example~2.

The remainder of the paper is organized as follows. After a brief review of background concepts in Section \ref{section: Preliminaries}, we present the full implementation of the numerical scheme in Section \ref{section: Description of the algorithm}. Its convergence analysis is given in Section \ref{section: Convergence analysis}, where we quantify the errors introduced by each approximation step. We conclude with numerical examples in Section \ref{section: numerical examples}, including an empirical verification of the convergence rates established above and applications to the pricing and hedging of contingent claims under rough-volatility models. Some technical proofs are deferred to the Appendix.

\section{Preliminaries}\label{section: Preliminaries}

\subsection{Definitions and notation}

Let $B = (B_{t})_{t \in [0,T]}$ be a standard $d$-dimensional Brownian motion defined on a complete probability space $(\Omega, \mathcal{F}, \mathbb{P})$. Let $p \geq 1$, $n \in \mathbb{N}$. We recall the following spaces and notational conventions:

\begin{itemize}
    \item $\mathbb{F} = (\mathcal{F}_{t})_{ t \in [0,T]}$ the filtration generated by the Brownian motion $B$ and augmented with the $\mathbb{P}$-null sets. We assume $\mathcal{F} = \mathcal{F}_{T}$.

    \item $\mathcal{P}(\mathbb{F})$ the progressive $\sigma$-algebra on $[0,T]\times\Omega$.

    \item $L^{p}(\mathcal{F}_t)$, $t \in [0,T]$, the space of all $\mathcal{F}_t$-measurable random variables $X\colon \Omega \to \mathbb{R}$ satisfying $\norm{X}_{L^p}^p \coloneqq \mathbb{E}\!\left[\abs{X}^{p}\right] < \infty$.

    \item $\mathbb{E}_{t}(X)$ denotes $\mathbb{E}(X \vert \mathcal{F}_{t})$ for any $X \in L^{1}(\mathcal{F}_T)$ and $t \in [0,T]$.

    \item $S^{p}(\mathbb{R}^{n})$ the space of all càdlàg adapted processes $\varphi \colon \Omega \times [0,T] \to \mathbb{R}^{n}$ such that $\norm{\varphi}_{S^{p}}^{p} \coloneqq \mathbb{E}\!\left[\sup_{t \in [0,T]}\abs{\varphi_t}^{p}\right] < \infty$.

    \item $H^{p}(\mathbb{R}^{n})$, $p \geq 2$, the space of all $\mathbb{F}$-progressively measurable processes $\varphi \colon \Omega \times [0,T] \to \mathbb{R}^{n}$ such that $\norm{\varphi}_{H^{p}}^{p} \coloneqq \mathbb{E}\!\left[\left(\int_{0}^{T}\abs{\varphi_t}^{2}\,dt\right)^{p/2}\right] < \infty$.

    \item $C_{\text{pol}}^{\infty}$ the space of smooth functions $f: \mathbb{R}^{n} \to \mathbb{R}$ such that $f$ and all its partial derivatives have polynomial growth.
\end{itemize}

We also recall some basic definitions related to Malliavin calculus:

\begin{itemize}

    \item $\mathbb{H} = L^{2}([0,T]; \mathbb{R}^{d})$ the space of functions $h\colon[0,T] \to \mathbb{R}^{d}$ equipped with the inner product
    $\langle h,g\rangle_{\mathbb{H}} \coloneqq \int_0^T h(s)\cdot g(s)\,ds$
    and the associated norm $\norm{h}_{\mathbb{H}}$.
    
    \item For $h = (h^{1}, \dots, h^{d}) \in \mathbb{H}$, we define the real-valued random variable
    \begin{gather*}
        B(h) \coloneqq \int_{0}^{T} h(s)\cdot dB_s
        =
        \sum_{j=1}^{d}
        \int_{0}^{T} h^{j}(s)\,dB_s^{j}.
    \end{gather*}

   \item $\mathcal{S}$ denotes the class of smooth random variables, which have the form
   \begin{gather*}
       F = \phi(B(h_{1}), \dots, B(h_{n})),
   \end{gather*}
   where $\phi\in C_{\text{pol}}^{\infty}(\mathbb{R}^{n})$ and $h_{1},\dots,h_n\in\mathbb H$.

   \item The Malliavin derivative of a random variable $F \in \mathcal{S}$ can be defined as the $d$-dimensional stochastic process
   \begin{gather*}
       D_{s}^{j} F = \sum_{i=1}^{n} \partial_i\phi\big(B(h_{1}), \dots, B(h_{n})\big)h_{i}^{j}(s), \quad s \in [0,T], \quad j=1,\dots,d.
   \end{gather*}
   Notice that $DF$ can be regarded as a random variable taking values in $\mathbb{H}$.

   \item More generally, for $k \geq 2$, we define the $k$-th Malliavin derivative of $F \in \mathcal{S}$ along the direction $\alpha = (\alpha_{1}, \dots, \alpha_{k}) \in \{1, \dots, d\}^{k}$ evaluated at $\mathbf{s} = (s_{1}, \dots, s_{k}) \in [0,T]^{k}$ as the real-valued random variable
   \begin{gather*}
        D_{\mathbf{s}}^{\alpha} F
        \coloneqq
        \sum_{i_{1}, \dots, i_{k} = 1}^{n}
        \partial_{i_1\cdots i_k}^{k}\phi\big(B(h_{1}), \dots, B(h_{n})\big)
        h_{i_{1}}^{\alpha_{1}}(s_{1}) \cdots
        h_{i_{k}}^{\alpha_{k}}(s_{k}).
    \end{gather*}
    The $k$-th Malliavin derivative of $F$ evaluated at $\mathbf{s}$ is then defined as $D_{\mathbf{s}}^{k} F = (D_{\mathbf{s}}^{\alpha} F)_{\alpha \in \{1, \dots, d\}^{k}}$. Notice that if we treat $\mathbf{s}$ as a parameter, $D^{k} F$ can be viewed as an $\mathbb{H}^{\otimes k}$-valued r.v., which by construction is symmetric in the $(s_{1}, \alpha_{1}) \dots, (s_{k}, \alpha_{k})$ variables.
    
    \item $\mathbb{D}^{k,p}$, $p \geq 2$, denotes the closure of $\mathcal{S}$ with respect to the norm given by
    \begin{flalign*}
        \norm{F}_{\mathbb{D}^{k,p}}^{p}
        &\coloneqq
        \mathbb{E}\!\left[\abs{F}^{p}\right]
        +
        \sum_{j=1}^{k}
        \mathbb{E}\!\left[\norm{D^{j}F}_{\mathbb{H}^{\otimes j}}^{p}\right],
    \end{flalign*}
    where 
    \begin{flalign*}
        \norm{D^jF}_{\mathbb{H}^{\otimes j}}^{p}
        =
        \left(
        \sum_{\alpha\in\{1,\dots,d\}^{j}}
        \int_{[0,T]^{j}}\abs{D_{\mathbf{s}}^{\alpha}F}^{2}\,d\mathbf{s}
        \right)^{p/2}.
    \end{flalign*}
    \item $\mathbb{L}_{a}^{1,p}$, $p\geq 2$, denotes the space of progressively measurable processes
    $u:[0,T]\times\Omega\to\mathbb{R}$ such that $u_t\in\mathbb{D}^{1,p}$ for $dt$-a.e. $t\in[0,T]$ and
    \[
        \mathbb{E}\!\left[
        \left(\int_0^T\abs{u_t}^{2}\,dt\right)^{p/2}
        +
        \left(\int_0^T\int_0^T\abs{D_su_t}^{2}\,ds\,dt\right)^{p/2}
    \right]<\infty.
\]
\end{itemize}

Finally, when working with random fields $f: [0,T]\times\Omega\times\mathbb{R}\times\mathbb{R}^{d}\to\mathbb{R}$, we use the notation $(D_s f)(t,Y,Z)$ to indicate that the Malliavin derivative is first taken with respect to $\omega$ for fixed $(t,y,z)$, and that the resulting random field is then evaluated at $(y,z)=(Y,Z)$. This should not be confused with the Malliavin derivative of the composition $f(t,Y,Z)$, which we denote by $D_s f(t,Y,Z)$. The same convention is used for higher-order Malliavin derivatives and for Malliavin derivatives of the partial derivatives of $f$.

\subsection{BSDEs and their numerical approximation}
\label{Section Numerical Scheme}

We recall the following standard well-posedness result for BSDEs; see
\textcite[Theorem 4.3.1]{ZhangJianfeng2017BSDE}.

\begin{theorem}\label{theorem 2}
Let
$f \colon [0,T] \times \Omega \times \mathbb{R} \times \mathbb{R}^{d}
\to \mathbb{R}$
be $\mathcal{P}(\mathbb{F}) \otimes
\mathcal{B}(\mathbb{R} \times \mathbb{R}^{d})$-measurable and uniformly
Lipschitz continuous in $(y,z)$, that is, assume that there exists a
constant $[f]_{L}>0$ such that, for $dt \otimes d\mathbb{P}$-a.e.
$(t,\omega)$,
\[
\abs{
f(t,\omega,y_{1},z_{1})
-
f(t,\omega,y_{2},z_{2})
}
\leq
[f]_{L}
\big(
\abs{y_{1}-y_{2}}
+
\abs{z_{1}-z_{2}}
\big)
\]
for all $(y_{1},z_{1}),(y_{2},z_{2})
\in \mathbb{R} \times \mathbb{R}^{d}$. Assume moreover that $\int_{0}^{T}\abs{f(t,0,0)}dt \in L^{2}(\mathcal{F}_{T})$ and that $\xi\in L^{2}(\mathcal{F}_{T})$. Then there exists a unique pair
\[
(Y,Z)
\in
S^{2}(\mathbb{R})
\times
H^{2}(\mathbb{R}^{d})
\]
satisfying \eqref{BSDE} for every $t\in[0,T]$, $\mathbb{P}$-a.s.
\end{theorem}

Throughout the paper, we use the backward Euler scheme introduced in \eqref{Euler scheme i}--\eqref{Euler scheme ii}. For later use, we make explicit the martingale representation defining $(Z_s^\pi)_{s\in[t_{i-1},t_i)}$:
\begin{gather}\label{conditional expectation and martingale term}
F_i^\pi
=
\mathbb{E}_{t_{i-1}}\big(F_i^\pi\big)
+
\int_{t_{i-1}}^{t_i}
Z_s^\pi\cdot dB_s,
\qquad
i=1,\dots,m.
\end{gather}
It follows that
\begin{gather*}
\sprod{Z}_{i}^{\pi}
=
\frac{1}{\Delta_{i+1}}
\mathbb{E}_{t_i}
\Big[
F_{i+1}^{\pi}
\big(
B_{t_{i+1}}-B_{t_i}
\big)
\Big],
\qquad
i=1,\dots,m-1.
\end{gather*}
Moreover, if $F_i^\pi\in\mathbb{D}^{1,2}$, the Clark--Ocone formula yields
\begin{gather*}
Z_t^\pi
=
\mathbb{E}_t\big[D_tF_i^\pi\big]
=
D_t\mathbb{E}_t\big(F_i^\pi\big)
=
D_tY_t^\pi,
\qquad
t\in[t_{i-1},t_i),
\end{gather*}
where the identities hold $dt\otimes d\mathbb{P}$-a.e. Thus, implementing the scheme amounts to computing conditional expectations and their Malliavin derivatives. In the following sections, we perform these computations through finite Wiener chaos expansions of the random variables $F_i^\pi$.

\subsection{The Wiener chaos decomposition}\label{subsection: chaos decomposition}

We now present an elementary introduction to the Wiener chaos decomposition of $L^2(\mathcal{F}_{T})$ using Hermite polynomials. We refer to \textcite[Section 1.1]{NualartDavidTMCa} or \textcite[Chapter 1]{DiNunnoGiulia2009MCfL} for more details on this topic. 

Let $H_{n}(x)$ denote the $n$th Hermite polynomial, defined by 
\begin{gather*}
    H_{n}(x) = \frac{(-1)^{n}}{n!} e^{\frac{x^{2}}{2}} \frac{d^{n}}{dx^{n}}\big( e^{\frac{-x^{2}}{2}} \big), \quad n \geq 1,
\end{gather*}
and $H_{0}(x) = 1$. These are the coefficients of the expansion of the power-series in $t$ of the function $F(x,t) = \exp\big(tx - \frac{t^{2}}{2} \big)$. In fact, we have
\begin{flalign*}
    F(x,t) = \exp \Big( \frac{x^{2}}{2} - \frac{1}{2}(x-t)^{2} \Big)
    &= \sum_{n=0}^{\infty} t^{n} H_{n}(x).
\end{flalign*}
Furthermore, since $\partial_{x} F(x,t) = t \exp \Big( \frac{x^{2}}{2} - \frac{1}{2}(x-t)^{2} \Big)$, one has that $H_{n}'(x) = H_{n-1}(x)$ for all $n \geq 0$, with the convention $H_{-1}(x) \equiv 0$. 

Let $\mathcal{A}$ denote the set of all sequences $a=(a_1,a_2,\dots)$, $a_j\in\mathbb{N}\cup\{0\}$, such that all but finitely many terms vanish. For $a\in\mathcal{A}$, set $a!=\prod_{j\geq1}a_j!$ and $|a|=\sum_{j\geq1}a_j$. For $M\in\mathbb{N}$, set $\mathcal{A}_M\coloneqq\{a\in\mathcal{A}:a_j=0\ \forall j>M\}$. For $a\in\mathcal{A}$, define
\begin{gather*}
    \Phi_a
    \coloneqq
    \sqrt{a!}
    \prod_{j\geq1}
    H_{a_j}\Big(
    \int_0^T h_j(s)\cdot dB_s
    \Big),
\end{gather*}
where $(h_j)_{j\geq1}$ is a fixed orthonormal basis of $\mathbb H$. Notice that the above product is well defined because $H_{0}(x) = 1$ and $a_{j} \neq 0$ only for a finite number of indices.

For each $n\geq0$, let $\mathcal{H}_n$ denote the closed linear subspace of $L^2(\mathcal{F}_T)$ generated by
$\{\Phi_a:a\in\mathcal{A},\,|a|=n\}$, called the Wiener chaos of order $n$. We then have the following Wiener chaos decomposition; see \textcite[Theorem 1.1.1 and Proposition 1.1.1]{NualartDavidTMCa}.

\begin{theorem}
The space $L^2(\mathcal{F}_T)$ is the orthogonal direct sum of the spaces $\mathcal{H}_n$, $n\geq0$, and
$\{\Phi_a:a\in\mathcal{A},\,|a|=n\}$ is a complete orthonormal system in $\mathcal{H}_n$. Consequently, every $F\in L^2(\mathcal{F}_T)$ admits the decomposition
\begin{gather}
    F
    =
    \sum_{n\geq0}
    \sum_{\substack{a\in\mathcal{A}\\ |a|=n}}
    \sprod{\Phi_a,F}_{L^2(\mathcal{F}_T)}
    \Phi_a,
    \label{chaos expansion hermite}
\end{gather}
where the series converges in $L^2(\mathcal{F}_T)$.
\end{theorem}

Equivalently, we may write
\begin{gather*}
    F
    =
    \sum_{n\geq0}
    \sum_{\substack{a\in\mathcal{A}\\ |a|=n}}
    d_a
    \prod_{j\geq1}
    H_{a_j}\Big(
    \int_0^T h_j(s)\cdot dB_s
    \Big),
    \quad
    d_a
    \coloneqq
    a!\,
    \mathbb{E}\Bigg[
    F
    \prod_{j\geq1}
    H_{a_j}\Big(
    \int_0^T h_j(s)\cdot dB_s
    \Big)
    \Bigg].
\end{gather*}

\begin{definition}
Let $F\in L^{2}(\mathcal{F}_{T})$ have representation \eqref{chaos expansion hermite}. For $n,P\in\mathbb{N}\cup\{0\}$ and $M\in\mathbb{N}$, define the projection onto the $n$th chaos by
\begin{flalign*}
    \mathcal{P}_{n}(F)
    &\coloneqq
    \sum_{\substack{a\in\mathcal{A}\\ |a|=n}}
    d_a
    \prod_{j\geq1}
    H_{a_j}\Big(
    \int_{0}^{T} h_j(s)\cdot dB_s
    \Big).
\end{flalign*}
We also define the truncation at chaos order $P$ and its finite-dimensional version involving only the first $M$ Gaussian variables, respectively, by
\begin{flalign*}
    \mathcal{C}_{P}(F)
    &\coloneqq
    \sum_{n=0}^{P}\mathcal{P}_{n}(F),\quad
    \mathcal{C}_{P}^{M}(F)
    \coloneqq
    \sum_{n=0}^{P}
    \sum_{\substack{a\in\mathcal{A}_M\\ |a|=n}}
    d_a
    \prod_{j=1}^{M}
    H_{a_j}\Big(
    \int_{0}^{T}h_j(s)\cdot dB_s
    \Big).
\end{flalign*}
\end{definition}

For $M\in\mathbb N$, let $\Pi_M$ denote the orthogonal projection of $\mathbb H$ onto the space spanned by $h_1,\dots,h_M$, that is,
\begin{gather*}
    \Pi_M h
    \coloneqq
    \sum_{j=1}^M
    \sprod{h,h_j}_{\mathbb H}h_j.
\end{gather*}
We use the same notation for its pointwise extension to $L^2(\Omega;\mathbb H)$.

Notice that $\mathcal{P}_{n}$ and $\mathcal{C}_{P}$ do not depend on the choice of the orthonormal basis, whereas $\mathcal{C}_{P}^{M}$ and the projection $\Pi_M$ do. We suppress this dependence in the notation.

We finally record some properties of the truncated chaos decomposition.

\begin{lemma}\label{lemma inequality chaos truncation}
Let $F \in L^{2}(\mathcal{F}_{T})$. Then, for all $P, M \in \mathbb{N}$, we have that:

\begin{enumerate}[label=(\roman*)]
\item $\mathbb{E}(|\mathcal{C}_{P}^{M}(F)|^{2}) \leq \mathbb{E}(|\mathcal{C}_{P}(F)|^{2}) \leq \mathbb{E}(|F|^{2})$.
\item If $F \in \mathbb{D}^{k,2}$, then for any $\alpha \in \{1, \dots, d\}^{k}$ we have $D_{\mathbf{s}}^{\alpha} \mathcal{C}_{P}(F) = \mathcal{C}_{P-k}(D_{\mathbf{s}}^{\alpha} F)$ for a.e. $\mathbf{s} \in [0,T]^{k}$, where we use the convention $\mathcal{C}_{Q}(\cdot) = 0$ for $Q < 0$.
\item If $F \in \mathbb{D}^{1,2}$, then $\mathbb{E}\big(\norm{D\mathcal{C}_{P}^{M}(F)}_{\mathbb{H}}^{2}\big) \leq \mathbb{E}\big(\norm{DF}_{\mathbb{H}}^{2}\big)$.
\end{enumerate}
\end{lemma}

\subsection{Error due to the truncation of the chaos decomposition}
\label{subsection: Error due to the truncation of the chaos decomposition}

We study the errors introduced by truncating the order of the chaos decomposition and the number of basis functions.

\begin{lemma}\label{lemma truncation order chaos}
Let $k \geq 1$ and $F \in \mathbb{D}^{k,2}$. Then, for any $P \geq k-1$,
\begin{gather*}
    \mathbb{E}\big(
    \abs{F-\mathcal{C}_{P}(F)}^{2}
    \big)
    \leq
    \mathbb{E}\big(
    \norm{D^{k}F}_{\mathbb{H}^{\otimes k}}^{2}
    \big)
    \frac{(P-k+1)!}{(P+1)!}.
\end{gather*}
\end{lemma}

\begin{proof}
By orthogonality of the Wiener chaoses,
\begin{flalign*}
    \mathbb{E}\big(
    \abs{F-\mathcal{C}_{P}(F)}^{2}
    \big)
    &=
    \sum_{n\geq P+1}
    \mathbb{E}\big(
    \abs{\mathcal{P}_{n}(F)}^{2}
    \big)
    \\
    &\leq
    \frac{(P-k+1)!}{(P+1)!}
    \sum_{n\geq0}
    n(n-1)\cdots(n-k+1)
    \mathbb{E}\big(
    \abs{\mathcal{P}_{n}(F)}^{2}
    \big)
    \\
    &=
    \frac{(P-k+1)!}{(P+1)!}
    \mathbb{E}\big(
    \norm{D^{k}F}_{\mathbb{H}^{\otimes k}}^{2}
    \big),
\end{flalign*}
where the last identity follows from \textcite[p.~28]{NualartDavidTMCa}.
\end{proof}

\begin{lemma}\label{lemma truncation chaos M}
Let $F\in\mathbb{D}^{1,2}$. Then, for all $P,M\in\mathbb{N}$,
\begin{flalign}
    \mathbb{E}\big(
    \abs{
    (\mathcal{C}_P^{M}-\mathcal{C}_{P})(F)
    }^{2}
    \big)
    &\leq
    \mathbb{E}\big(
    \norm{
    (I-\Pi_M)DF
    }_{\mathbb{H}}^{2}
    \big).
    \label{general basis truncation bound}
\end{flalign}
\end{lemma}

\begin{proof}
The left hand side can be expressed as
\begin{gather*}
    \mathbb{E}\big(
    \abs{
    (\mathcal{C}_P^{M}-\mathcal{C}_{P})(F)
    }^{2}
    \big)
    =
    \sum_{\substack{a\in\mathcal{A}\\ |a|\leq P\\ \sum_{j>M}a_j\geq1}}
    \abs{
    \sprod{\Phi_a,F}_{L^2(\mathcal{F}_T)}
    }^{2}.
\end{gather*}
Using $H_n'=H_{n-1}$, we have
\begin{gather*}
    D\Phi_a
    =
    \sqrt{a!}
    \sum_{\substack{j\geq1\\a_j\geq1}}
    H_{a_j-1}\Big(
    \int_0^T h_j(s)\cdot dB_s
    \Big)
    \prod_{\ell\neq j}
    H_{a_\ell}\Big(
    \int_0^T h_\ell(s)\cdot dB_s
    \Big)
    h_j.
\end{gather*}
Since $F\in\mathbb{D}^{1,2}$, its differentiated chaos expansion converges in $L^2(\Omega;\mathbb{H})$. Consequently, by the orthonormality of $(h_j)_{j\geq1}$ and $(\Phi_a)_{a\in\mathcal{A}}$,
\begin{gather*}
    \mathbb{E}\big(
    \norm{
    (I-\Pi_M)DF
    }_{\mathbb{H}}^{2}
    \big)
    =
    \sum_{a\in\mathcal{A}}
    \abs{
    \sprod{\Phi_a,F}_{L^2(\mathcal{F}_T)}
    }^{2}
    \sum_{j>M}a_j.
\end{gather*}
Since $\sum_{j>M}a_j\geq1$ for every multi-index appearing in the basis-truncation error,
\begin{flalign*}
    \mathbb{E}\big(
    \abs{
    (\mathcal{C}_P^{M}-\mathcal{C}_{P})(F)
    }^{2}
    \big)
    &\leq
    \sum_{\substack{a\in\mathcal{A}\\ |a|\leq P\\ \sum_{j>M}a_j\geq1}}
    \abs{
    \sprod{\Phi_a,F}_{L^2(\mathcal{F}_T)}
    }^{2}
    \sum_{j>M}a_j
    \\
    &\leq
    \mathbb{E}\big(
    \norm{
    (I-\Pi_M)DF
    }_{\mathbb{H}}^{2}
    \big).
\end{flalign*}
\end{proof}

We next derive an explicit rate for the basis-truncation error associated with the parameter $M$, for the piecewise constant basis used in the implementation. We state the result for $d=1$; the extension to $d>1$ follows component-wise with the same rate.

\begin{lemma}\label{lemma approximation holder}
Let $F\in\mathbb{D}^{1,2}$, let
$\overline{\pi}=\{0=s_0<\cdots<s_M=T\}$ be a partition of $[0,T]$,
set $\delta_i=s_i-s_{i-1}$, and define
$h_i(t)=\mathbf{1}_{(s_{i-1},s_i]}(t)/\sqrt{\delta_i}$ for $i=1,\dots,M$.
Complete $(h_i)_{1\leq i\leq M}$ to an orthonormal basis of $\mathbb H$, with respect to which $\mathcal{C}_P^M$ is defined.

Assume that there exist $\kappa>0$ and $\beta\in(0,1]$ such that
\begin{gather}\label{Holder regularity Malliavin derivative}
    \mathbb{E}\big(
    \abs{D_sF-D_rF}^{2}
    \big)
    \leq
    \kappa\abs{s-r}^{2\beta},
    \qquad
    s,r\in[0,T].
\end{gather}
Then, for every $P\in\mathbb{N}$,
\begin{gather}
    \mathbb{E}\big(
    \abs{
    (\mathcal{C}_P^{M}-\mathcal{C}_{P})(F)
    }^{2}
    \big)
    \leq
    \frac{2\kappa T}
    {(2\beta+1)(2\beta+2)}
    |\overline{\pi}|^{2\beta}.
    \label{piecewise constant basis truncation bound}
\end{gather}
\end{lemma}

\begin{proof}
The orthogonal projection onto
$\operatorname{span}\{h_1,\dots,h_M\}$ is given by
\begin{gather*}
    \Pi_M DF
    =
    \sum_{i=1}^{M}
    \left(
    \frac{1}{\delta_i}
    \int_{s_{i-1}}^{s_i}D_rF\,dr
    \right)
    \mathbf{1}_{(s_{i-1},s_i]}.
\end{gather*}
Hence, by Jensen's inequality and \eqref{Holder regularity Malliavin derivative},
\begin{flalign*}
    \mathbb{E}\big(
    \norm{(I-\Pi_M)DF}_{\mathbb{H}}^{2}
    \big)
    &\leq
    \sum_{i=1}^{M}
    \frac{1}{\delta_i}
    \int_{s_{i-1}}^{s_i}
    \int_{s_{i-1}}^{s_i}
    \mathbb{E}\big(
    \abs{D_sF-D_rF}^{2}
    \big)\,dr\,ds \\
    &\leq
    \kappa
    \sum_{i=1}^{M}
    \frac{1}{\delta_i}
    \int_{s_{i-1}}^{s_i}
    \int_{s_{i-1}}^{s_i}
    \abs{s-r}^{2\beta}\,dr\,ds \\
    &=
    \frac{2\kappa}
    {(2\beta+1)(2\beta+2)}
    \sum_{i=1}^{M}\delta_i^{2\beta+1} \\
    &\leq
    \frac{2\kappa T}
    {(2\beta+1)(2\beta+2)}
    |\overline{\pi}|^{2\beta}.
\end{flalign*}
The result follows from Lemma \ref{lemma truncation chaos M}.
\end{proof}

\begin{example}
Let $b,\sigma\colon[0,T]\times\mathbb{R}^{k}\to\mathbb{R}^{k}$ be globally Lipschitz in $x$, $1/2$-Hölder continuous in $t$, and twice continuously differentiable in $x$ with bounded derivatives. Let $(X_t)_{t\in[0,T]}$ solve
\begin{gather*}
    X_t
    =
    x
    +
    \int_0^t b(u,X_u)\,du
    +
    \int_0^t\sigma(u,X_u)\,dB_u,
    \qquad
    x\in\mathbb{R}^{k},
\end{gather*}
and let $g\colon\mathbb{R}^{k}\to\mathbb{R}$ be continuously differentiable with bounded derivative. Standard estimates for the Malliavin derivative of an SDE, whose variational equation is described in \textcite[Section 2.2]{NualartDavidTMCa}, yield
\begin{gather*}
    \mathbb{E}\big(
    \abs{D_sX_t-D_rX_t}^{2}
    \big)
    \leq
    C\abs{s-r}, \qquad 0\leq s,r\leq t\leq T.
\end{gather*}
Using the Malliavin chain rule, it follows that both $g(X_T)$ and
$\int_0^T g(X_t)\,dt$ satisfy \eqref{Holder regularity Malliavin derivative} with $\beta=1/2$.
\end{example}

\section{Description of the algorithm}\label{section: Description of the algorithm}

We now proceed to describe an implementation of the Euler scheme for BSDEs introduced in \eqref{Euler scheme i}--\eqref{Euler scheme ii}. The random variables $F_i^\pi$ are approximated by finite Wiener chaos projections. Using at each time step the piecewise constant functions considered in Lemma~\ref{lemma approximation holder}, we obtain closed-form formulas for the conditional expectations and martingale terms in \eqref{conditional expectation and martingale term}.

\subsection{Description of the truncated basis}\label{sec: description basis}

For notational simplicity, we present the one-dimensional case ($d=1$); the corresponding definitions and formulae for the multi-dimensional case are provided in Section \ref{section: multidimensional formulas}.

Fix $P,M\in\mathbb N$. We let $\overline{\pi}_m \coloneqq \{s_j^m \coloneqq jT/M\}_{0\leq j\leq M}$ be the uniform partition of $[0,T]$ with $M$ subintervals. For any $i\in\{1,\dots,m\}$, we consider the partition of $[0,t_i]$ given by
\begin{gather*}
    \overline{\pi}_i
    \coloneqq
    \{s\in\overline{\pi}_m:s<t_i\}\cup\{t_i\}.
\end{gather*}
We set $M(i)\coloneqq\#\overline{\pi}_i-1$ and denote the elements of $\overline{\pi}_i$ by $\{s_j^i\}_{0\leq j\leq M(i)}$. Notice that $M(i)\leq M(i+1)$ for $i=1,\dots,m-1$.

For each $i\in\{1,\dots,m\}$, we define
\begin{gather}\label{truncated basis}
    h_j^i(t)
    \coloneqq
    \frac{\mathbf{1}_{(s_{j-1}^i,s_j^i]}(t)}{\sqrt{\delta_j^i}},
    \qquad
    \delta_j^i
    \coloneqq
    s_j^i-s_{j-1}^i,
    \qquad
    j=1,\dots,M(i).
\end{gather}
The family $(h_j^i)_{1\leq j\leq M(i)}$ is an orthonormal system in $\mathbb H$ and can be completed to an orthonormal basis of $\mathbb H$. Since each $h_j^i$ vanishes outside $[0,t_i]$, we have
\begin{gather*}
    B(h_j^i)
    =
    \int_0^T h_j^i(s)\,dB_s
    =
    \frac{B_{s_j^i}-B_{s_{j-1}^i}}{\sqrt{\delta_j^i}},
    \qquad
    j=1,\dots,M(i).
\end{gather*}
We denote by $\mathcal{C}_{P}^{\overline{\pi}_i}$ the finite-dimensional chaos projection of order $P$ associated with the partition $\overline{\pi}_i$. That is, for $F\in L^{2}(\mathcal{F}_{t_i})$,
\begin{flalign*}
    \mathcal{C}_{P}^{\overline{\pi}_i}(F)
    &\coloneqq
    \sum_{k=0}^{P}
    \sum_{\substack{a\in\mathcal{A}_{M(i)}\\ |a|=k}}
    d_a
    \prod_{j=1}^{M(i)}
    H_{a_j}\big(B(h_j^i)\big),\quad
    d_a
    \coloneqq
    a!\,
    \mathbb{E}\Bigg[
    F
    \prod_{j=1}^{M(i)}
    H_{a_j}\big(B(h_j^i)\big)
    \Bigg].
\end{flalign*}

\subsection{Description of the implementation}

Set $\theta\coloneqq(P,M)$. We define the recursion $Y_{t_m}^{\pi,\theta}\coloneqq\mathcal{C}_{P}^{\overline{\pi}_m}(\xi)$ and, for $i=m,\dots,1$,
\begin{gather*}
    F_i^{\pi,\theta}\coloneqq Y_{t_i}^{\pi,\theta}+\Delta_i f_i^{\pi,\theta},
    \quad
    f_i^{\pi,\theta}\coloneqq f(t_i,Y_{t_i}^{\pi,\theta},\sprod{Z}_i^{\pi,\theta}),
\end{gather*}
where
\begin{flalign*}
    Y_t^{\pi,\theta}
    =
    \mathcal{C}_{P}^{\overline{\pi}_i}(F_i^{\pi,\theta})
    -
    \int_t^{t_i}Z_s^{\pi,\theta}\cdot dB_s,
    \quad
    t\in[t_{i-1},t_i),
\end{flalign*}
and
\begin{gather*}
    \sprod{Z}_i^{\pi,\theta}
    \coloneqq
    \mathbb{E}_{t_i}\Bigg(
    \frac{1}{\Delta_{i+1}}
    \int_{t_i}^{t_{i+1}}Z_s^{\pi,\theta}\,ds
    \Bigg),
    \quad
    i=1,\dots,m-1,
\end{gather*}
with the convention $\sprod{Z}_m^{\pi,\theta}=0$. For $t\in[t_{i-1},t_i)$,
\begin{gather*}
    Y_t^{\pi,\theta}
    =
    \mathbb{E}_t\big(\mathcal{C}_{P}^{\overline{\pi}_i}(F_i^{\pi,\theta})\big),
\end{gather*}
\begin{gather*}
    \sprod{Z}_i^{\pi,\theta}
    =
    \frac{1}{\Delta_{i+1}}
    \mathbb{E}_{t_i}\Big[
    \mathcal{C}_{P}^{\overline{\pi}_{i+1}}(F_{i+1}^{\pi,\theta})
    \big(B_{t_{i+1}}-B_{t_i}\big)
    \Big],
    \quad
    i=1,\dots,m-1,
\end{gather*}
and $(Z_s^{\pi,\theta})_{s\in[t_{i-1},t_i)}$ is determined by the martingale representation of
$\mathcal{C}_{P}^{\overline{\pi}_i}(F_i^{\pi,\theta})$. Since a random variable with a finite chaos decomposition belongs to $\mathbb{D}^{1,2}$,
\begin{flalign*}
    Z_t^{\pi,\theta}
    \coloneqq
    D_t\mathbb{E}_t\big(\mathcal{C}_{P}^{\overline{\pi}_i}(F_i^{\pi,\theta})\big),
    \quad
    t\in[t_{i-1},t_i),
\end{flalign*}
where the identity holds $dt\otimes d\mathbb{P}$-a.e.

\subsection{Wiener chaos decomposition formulae}\label{subsection: Wiener chaos decomposition formulae}

We now present closed-form formulae for the computation of $Y_t^{\pi,\theta}$, $Z_t^{\pi,\theta}$, and $\sprod{Z}_i^{\pi,\theta}$ on $[t_{i-1},t_i)$. The formulae for the multi-dimensional case ($d>1$) are presented in Section~\ref{section: multidimensional formulas}.

For each $i\in\{1,\dots,m\}$, we denote by $d_a^i$ the chaos coefficient of $\mathcal{C}_{P}^{\overline{\pi}_i}(F_i^{\pi,\theta})$ associated with $a\in\mathcal{A}_{M(i)}$. Related formulae appear in \textcite[Proposition 2.7]{BriandPhilippe2014SOBB}.

\begin{proposition}\label{proposition Y and Z}
Let $i\in\{1,\dots,m\}$ and $r\in\{1,\dots,M(i)\}$ be such that $t\in[t_{i-1},t_i)\cap(s_{r-1}^i,s_r^i]$. Then
\begin{flalign*}
Y_t^{\pi,\theta}
=
d_0^i
+
\sum_{k=1}^P
\sum_{\substack{a\in\mathcal{A}_{M(i)}^r\\|a|=k}}
d_a^i
\prod_{j<r}H_{a_j}\big(B(h_j^i)\big)
\times
\Bigg(\frac{t-s_{r-1}^i}{\delta_r^i}\Bigg)^{a_r/2}
H_{a_r}\Bigg(\frac{B_t-B_{s_{r-1}^i}}{\sqrt{t-s_{r-1}^i}}\Bigg),\\
Z_t^{\pi,\theta}
=
(\delta_r^i)^{-1/2}
\sum_{k=1}^P
\sum_{\substack{a\in\mathcal{A}_{M(i)}^r\\|a|=k\\a_r>0}}
d_a^i
\prod_{j<r}H_{a_j}\big(B(h_j^i)\big)
\times
\Bigg(\frac{t-s_{r-1}^i}{\delta_r^i}\Bigg)^{(a_r-1)/2}
H_{a_r-1}\Bigg(\frac{B_t-B_{s_{r-1}^i}}{\sqrt{t-s_{r-1}^i}}\Bigg),
\end{flalign*}
where $\mathcal{A}_{M(i)}^r\coloneqq\{a\in\mathcal{A}_{M(i)}:a_j=0\text{ for all }j>r\}$. 

At $t=0$, we use the continuous extensions $Y_0^{\pi,\theta}=d_0^1$ and $Z_0^{\pi,\theta}=(\delta_1^1)^{-1/2}d_{e_1}^1$, where $e_1\coloneqq(1,0,\dots)\in\mathcal{A}_{M(1)}$.
\end{proposition}

We finally provide a closed-form formula for $\sprod{Z}_{i-1}^{\pi,\theta}$. The proof for the multi-dimensional case can be found in Proposition~\ref{prop:formula-for-scalarZ} in the Appendix.

\begin{proposition}\label{Proposition Z bar}
Let $i\in\{2,\dots,m\}$ and let $u$ be such that $t_{i-1}\in(s_{u-1}^i,s_u^i]$. Set
\begin{gather*}
C_1(i)\coloneqq\frac{s_u^i-t_{i-1}}{\sqrt{\delta_u^i}},
\qquad
C_2(i)\coloneqq\frac{t_{i-1}-s_{u-1}^i}{\delta_u^i}.
\end{gather*}
Then
\begin{flalign*}
\sprod{Z}_{i-1}^{\pi,\theta}
&=
\frac{1}{\Delta_i}
\sum_{k=1}^P
\Bigg\{
C_1(i)
\sum_{\substack{a\in\mathcal{A}_{M(i)}^u\\|a|=k\\a_u>0}}
d_a^i
\prod_{j<u}H_{a_j}\big(B(h_j^i)\big)
\times
\big[C_2(i)\big]^{(a_u-1)/2}
H_{a_u-1}\Bigg(
\frac{B_{t_{i-1}}-B_{s_{u-1}^i}}
{\sqrt{t_{i-1}-s_{u-1}^i}}
\Bigg)\\
&\quad+
\sum_{r=u+1}^{M(i)}
\sqrt{\delta_r^i}
\sum_{\substack{a\in\mathcal{A}_{M(i)}^{u,r}\\|a|=k}}
d_a^i
\prod_{j<u}H_{a_j}\big(B(h_j^i)\big)
\times
\big[C_2(i)\big]^{a_u/2}
H_{a_u}\Bigg(
\frac{B_{t_{i-1}}-B_{s_{u-1}^i}}
{\sqrt{t_{i-1}-s_{u-1}^i}}
\Bigg)
\Bigg\},
\end{flalign*}
where $\mathcal{A}_{M(i)}^{u,r}
\coloneqq
\{a\in\mathcal{A}_{M(i)}:a_r=1,\ a_j=0\text{ for all }j>u,\ j\neq r\}$. 
\end{proposition}

\subsection{Monte Carlo approximation of the chaos coefficients}\label{subsection: MC approximation}

To implement the method, one needs to compute the coefficients $d_a^i$ of the truncated chaos projection of $F_i^{\pi,\theta}\in L^2(\mathcal{F}_{t_i})$, which, for $a\in\mathcal{A}_{M(i)}$ with $|a|\leq P$, are given by
\begin{flalign}
\frac{d_a^i}{a!}
&=
\mathbb{E}\Bigg[
F_i^{\pi,\theta}
\prod_{j=1}^{M(i)}H_{a_j}\big(B(h_j^i)\big)
\Bigg]\notag\\
&=
\mathbb{E}\Bigg[
Y_{t_i}^{\pi,\theta}
\prod_{j=1}^{M(i)}H_{a_j}\big(B(h_j^i)\big)
\Bigg]
+
\Delta_i\mathbb{E}\Bigg[
f(t_i,Y_{t_i}^{\pi,\theta},\sprod{Z}_i^{\pi,\theta})
\prod_{j=1}^{M(i)}H_{a_j}\big(B(h_j^i)\big)
\Bigg].
\label{eq:expect}
\end{flalign}
For $i=1,\dots,m-1$, Proposition~\ref{proposition Y and Z} gives the coefficients of $\mathcal{C}_{P}^{\overline{\pi}_i}(Y_{t_i}^{\pi,\theta})$ in closed form. Indeed, $t_i\in(s_{M(i)-1}^{i+1},s_{M(i)}^{i+1}]$ and, for $a\in\mathcal{A}_{M(i)}$, viewed as an element of $\mathcal{A}_{M(i+1)}$,
\begin{gather*}
a!\,
\mathbb{E}\Bigg[
Y_{t_i}^{\pi,\theta}
\prod_{j=1}^{M(i)}H_{a_j}\big(B(h_j^i)\big)
\Bigg]
=
d_a^{i+1}
\Bigg(
\frac{t_i-s_{M(i)-1}^{i+1}}{\delta_{M(i)}^{i+1}}
\Bigg)^{a_{M(i)}/2}.
\end{gather*}
Hence, apart from the coefficients of the terminal projection $\mathcal{C}_{P}^{\overline{\pi}_m}(\xi)$, we only need to approximate the second expectation in \eqref{eq:expect}. We do so by Monte Carlo simulation, as in \textcite{BriandPhilippe2014SOBB}.

Let $N\in\mathbb{N}$ be the number of samples used to approximate these expectations. We use different Brownian motions at each time step and therefore consider:
\begin{itemize}
    \item $\overline{B}=(B^{i,n})_{0\leq i\leq m,\,1\leq n\leq N}$ and $B^\xi=(B^{\xi,n})_{1\leq n\leq N}$ are an $(m+1)N$-dimensional and an $N$-dimensional Brownian motion, respectively, independent of $B$ and of each other. We set $B^{i,0}=B$ for $0\leq i\leq m$ and $B^{\xi,0}=B$. For $h\in\mathbb H$, we write $B^{i,n}(h)\coloneqq\int_0^T h(s)\,dB_s^{i,n}$ and $B^{\xi,n}(h)\coloneqq\int_0^T h(s)\,dB_s^{\xi,n}$.

    \item For $0\leq i\leq m$ and $1\leq n\leq N$, $\mathbb{F}^{i,n}=(\mathcal{F}_t^{i,n})_{t\in[0,T]}$ denotes the natural augmented filtration generated by $B^{i,n}$. Similarly, $\mathbb{F}^{\xi,n}=(\mathcal{F}_t^{\xi,n})_{t\in[0,T]}$ denotes the natural augmented filtration generated by $B^{\xi,n}$. We use the conventions $\mathbb{F}^{i,0}=\mathbb{F}^{\xi,0}=\mathbb{F}$.

    \item For $1\leq i\leq m$ and $1\leq n\leq N$, let $f^{i,n}$ be copies of $f$ such that the pairs $(B^{i,n},f^{i,n})$ are mutually independent copies of $(B,f)$ and $f^{i,n}(t,y,z)$ is $\mathcal{F}_t^{i,n}$-measurable for every $(t,y,z)$. We set $f^{i,0}=f$.

    \item Let $(B^{\xi,n},\xi^n)_{1\leq n\leq N}$ be mutually independent copies of $(B,\xi)$, independent of the previous families, and set $\xi^0=\xi$.

    \item $(\overline{\Omega},\overline{\mathcal{F}},\overline{\mathbb{P}},\overline{\mathbb{F}})$ is the enlarged filtered probability space carrying all the objects above.

    \item $\overline{\mathbb{E}}$ denotes expectation on this probability space, and $\overline{D}^{(i,n)}$ denotes the Malliavin derivative with respect to $B^{i,n}$, with $\overline{D}^{(i,0)}=D$.

    \item For $0\leq i\leq m$, $\mathbb{F}^i=(\mathcal{F}_t^i)_{t\in[0,T]}$ denotes the natural augmented filtration generated by $(B^{i,1},\dots,B^{i,N})$. Similarly, $\mathbb{F}^\xi=(\mathcal{F}_t^\xi)_{t\in[0,T]}$ denotes the natural augmented filtration generated by $(B^{\xi,1},\dots,B^{\xi,N})$.

    \item For $i=1,\dots,m$, set $\mathcal{F}^{i\leq}\coloneqq\mathcal{F}_{t_i}^i\vee\mathcal{F}_{t_{i+1}}^{i+1}\vee\cdots\vee\mathcal{F}_{t_m}^m\vee\mathcal{F}_{t_m}^\xi$, with the convention $\mathcal{F}^{m+1\leq}=\mathcal{F}_{t_m}^\xi$.

    \item For $i=1,\dots,m$, $n=0,\dots,N$, and $a\in\mathcal{A}_{M(i)}$, set
    \[
    H_a^{(i,n)}\coloneqq\prod_{j=1}^{M(i)}H_{a_j}\big(B^{i,n}(h_j^i)\big),
    \qquad
    \widehat{H}_a^{(i,n)}\coloneqq\prod_{j=1}^{M(i)}H_{a_j}\big(B^{i-1,n}(h_j^i)\big).
    \]
    Similarly, for $n=0,\dots,N$ and $a\in\mathcal{A}_{M(m)}$, set $H_a^{(\xi,n)}\coloneqq\prod_{j=1}^{M(m)}H_{a_j}\big(B^{\xi,n}(h_j^m)\big)$.
\end{itemize}
The random variables $H_a^{(i,n)}$ are used to approximate the coefficient $d_a^i$, whereas $\widehat{H}_a^{(i,n)}$ are used to evaluate the approximate chaos decomposition at step $i$.

We then recursively define the following approximation of the truncated chaos decomposition of $F_i^{\pi,\theta}$.

\begin{definition}
For $n=0,\dots,N$, let
\begin{gather*}
    \widehat{\mathcal{C}}_{P}^{\overline{\pi}_m,n}(\xi) \coloneqq \hat{d}_0^\xi + \sum_{k=1}^{P}\sum_{\substack{a\in\mathcal{A}_{M(m)}\\|a|=k}}\hat{d}_a^\xi H_a^{(m,n)}, \quad \frac{\hat{d}_a^\xi}{a!}\coloneqq\frac{1}{N}\sum_{l=1}^{N}\xi^l\times H_a^{(\xi,l)}.
\end{gather*}
We then recursively define $\big(\widehat{Y}_{t_m}^{\pi,\theta,n},\sprod{\widehat{Z}}_m^{\pi,\theta,n}\big)\coloneqq\big(\widehat{\mathcal{C}}_{P}^{\overline{\pi}_m,n}(\xi),0\big)$. For $i\in\{m,\dots,1\}$ and $n=0,\dots,N$, set $\widehat{f}_i^{\pi,\theta,n}\coloneqq f^{i,n}(t_i,\widehat{Y}_{t_i}^{\pi,\theta,n},\sprod{\widehat{Z}}_i^{\pi,\theta,n})$ and $\widehat{F}_i^{\pi,\theta,n}\coloneqq\widehat{Y}_{t_i}^{\pi,\theta,n}+\Delta_i\widehat{f}_i^{\pi,\theta,n}$. Whenever $n$ is omitted from the superscript, it is understood that $n=0$. For $t\in[t_{i-1},t_i)$, define
\begin{gather*}
    \widehat{Y}_t^{\pi,\theta,n}=\overline{\mathbb{E}}\big(\widehat{\mathcal{C}}_{P}^{\overline{\pi}_i,n}(\widehat{F}_i^{\pi,\theta})\mid\mathcal{F}_t^{i-1,n}\vee\mathcal{F}^{i\leq}\big), \quad \widehat{Z}_t^{\pi,\theta,n}=\overline{D}_t^{(i-1,n)}\widehat{Y}_t^{\pi,\theta,n},
\end{gather*}
where
\begin{gather*}
    \widehat{\mathcal{C}}_{P}^{\overline{\pi}_i,n}(\widehat{F}_i^{\pi,\theta})\coloneqq\hat{d}_0^i+\sum_{k=1}^{P}\sum_{\substack{a\in\mathcal{A}_{M(i)}\\|a|=k}}\hat{d}_a^i\widehat{H}_a^{(i,n)},
\end{gather*}
with
\begin{flalign*}
    \frac{\hat{d}_a^i}{a!}\coloneqq\overline{\mathbb{E}}\Big(\widehat{Y}_{t_i}^{\pi,\theta,0}\times H_a^{(i,0)}\mid\mathcal{F}^{i+1\leq}\Big)+\frac{\Delta_i}{N}\sum_{l=1}^{N}f^{i,l}(t_i,\widehat{Y}_{t_i}^{\pi,\theta,l},\sprod{\widehat{Z}}_i^{\pi,\theta,l})\times H_a^{(i,l)}.
\end{flalign*}
For $i=2,\dots,m$ and $n=0,\dots,N$, $\sprod{\widehat{Z}}_{i-1}^{\pi,\theta,n}$ is defined by the formula in Proposition~\ref{Proposition Z bar}, replacing $d_a^i$ and $B$ with $\hat{d}_a^i$ and $B^{i-1,n}$, respectively.
\end{definition}

\begin{remark}\label{remark MC formulae}
The random variables $\hat{d}_a^i$ are $\mathcal{F}^{i\leq}$-measurable. Moreover, for each $n=1,\dots,N$, $\hat{d}_a^i$ is independent of $\widehat{H}_a^{(i,n)}$, since they are built from different Brownian families.

Consequently, the random variables $\{\widehat{Y}_t^{\pi,\theta,n},\widehat{Z}_t^{\pi,\theta,n},\sprod{\widehat{Z}}_{i-1}^{\pi,\theta,n}\}$ can be obtained via the formulae in Propositions~\ref{proposition Y and Z} and~\ref{Proposition Z bar} by replacing $d_a^i$ with $\hat{d}_a^i$ and using the Brownian motion $B^{i-1,n}$ in place of $B$.
\end{remark}

\begin{remark}
Under Assumption~\ref{main assumption convergence rate},
$\widehat{\mathcal C}_{P}^{\overline{\pi}_i,0}(\widehat F_i^{\pi,\theta})$
is square-integrable and $\mathcal F_{t_i}\vee\mathcal F^{i\leq}$-measurable.
Since $\mathcal F^{i\leq}$ is independent of $B$, the conditional martingale representation theorem yields
\[
    \widehat Y_{t_{i-1}}^{\pi,\theta,0}
    =
    \overline{\mathbb E}\left(
        \widehat{\mathcal C}_{P}^{\overline{\pi}_i,0}
        (\widehat F_i^{\pi,\theta})
        \,\middle|\,
        \mathcal F_{t_{i-1}}\vee\mathcal F^{i\leq}
    \right)
    =
    \widehat{\mathcal C}_{P}^{\overline{\pi}_i,0}
    (\widehat F_i^{\pi,\theta})
    -
    \int_{t_{i-1}}^{t_i}
    \widehat Z_t^{\pi,\theta,0}\,dB_t,
\]
where $(\mathcal F_t)_{t\in[0,T]}$ denotes the natural augmented filtration of $B$.
\end{remark}

\section{Convergence analysis}\label{section: Convergence analysis}

In this section, we present a detailed convergence analysis of the method introduced above. Although the numerical scheme is formulated for general square-integrable terminal conditions, the convergence analysis is carried out under the additional regularity assumption stated below. For notational simplicity, we consider the one-dimensional case ($d=1$); the arguments carry over directly to the multi-dimensional setting. Some proofs are deferred to Section~\ref{Proofs section convergence analysis} in the Appendix.

The assumption under which we carry out the convergence analysis is the following.

\begin{assumption}\label{main assumption convergence rate}
    Fix $q>4$ and $\beta_\xi,\beta_f\in[1/2,1]$.
    \begin{enumerate}[label=(\roman*)]
        \item $\xi \in \mathbb{D}^{2,q}$, and there exists $\mu_\xi >0$ such that, for all $s,r \in [0,T]$,
        \begin{flalign*}
           \mathbb{E}\big( \abs{D_s\xi-D_r\xi}^{2} \big) \leq \mu_\xi\abs{s-r}^{2\beta_\xi}&, \qquad
           \sup_{0 \leq s \leq T} \mathbb{E}\big( \abs{D_s\xi}^{q}\big) < \infty, \\
           & \sup_{0 \leq s,r \leq T} \mathbb{E}\big( \abs{D_sD_r\xi}^{q}\big) < \infty.
        \end{flalign*}

        \item For $dt\otimes d\mathbb{P}$-a.e. $(t,\omega)$, the map
        $(y,z)\mapsto f(t,\omega,y,z)$ is twice continuously differentiable, with all first- and second-order partial derivatives uniformly bounded in $(t,\omega,y,z)$. Moreover,
        $f(\cdot,0,0)\in H^q(\mathbb R)$.

        \item For every $(y,z)$, the processes
        $f(\cdot,y,z)$, $\partial_y f(\cdot,y,z)$, and
        $\partial_z f(\cdot,y,z)$ belong to $\mathbb{L}_{a}^{1,q}$.
        Moreover, for every $(t,y,z)$, $f(t,y,z)\in\mathbb D^{2,q}$ and,
        for every $s\in[0,T]$ and $(y,z)$,
        $(D_s f)(\cdot,y,z)\in\mathbb{L}_{a}^{1,q}$.
        Furthermore, for every $0\leq s\leq t\leq T$, the map
        $(y,z)\mapsto(D_s f)(t,y,z)$ has continuous partial derivatives.

        There exist non-negative random fields
        $\Gamma_s(t)$ and $K_{s,r}(t)$, and a constant $\mu_f>0$, such that, a.s., for all $s,r\leq t$ and $(y,z)$,
        \begin{gather*}
            \abs{(D_s f)(t,y,z)}+\abs{(D_s\partial_y f)(t,y,z)}+\abs{(D_s\partial_z f)(t,y,z)}
            \leq \Gamma_s(t),\\
            \abs{(D_s f)(t,y,z)-(D_r f)(t,y,z)}
            \leq K_{s,r}(t),
        \end{gather*}
        with
        \[
            \sup_{0\leq s\leq t\leq T}\mathbb E\abs{\Gamma_s(t)}^q<\infty,
            \qquad
            \sup_{s\vee r\leq t\leq T}\mathbb E|K_{s,r}(t)|^2
            \leq \mu_f\abs{s-r}^{2\beta_f}.
        \]
        In addition, denoting by $(Y,Z)$ the solution to \eqref{BSDE},
        \[
            \sup_{0\leq s,r\leq T}
            \mathbb E\left[
                \left(
                    \int_{s\vee r}^{T}
                    \abs{(D_sD_r f)(t,Y_t,Z_t)}^2\,dt
                \right)^{q/2}
            \right]
            <\infty.
        \]

        \item There exists $[f]_H>0$ such that, a.s.,
        \[
            \abs{f(s,y,z)-f(r,y,z)}
            \leq [f]_H\abs{s-r}^{1/2},
            \qquad s,r\in[0,T],\quad (y,z)\in\mathbb R^2.
        \]
    \end{enumerate}
\end{assumption}

Throughout the remainder of this section, we consider time partitions $\pi=\{0=t_0<\cdots<t_m=T\}$ underlying the Euler scheme \eqref{Euler scheme i}--\eqref{Euler scheme ii} and satisfying
\[
    \max_{1\leq i\leq m-1}\frac{\Delta_i}{\Delta_{i+1}}\leq L
\]
for some fixed constant $L>0$, and write $|\pi|\coloneqq\max_{1\leq i\leq m}\Delta_i$.

\subsection{Error given by the time discretization}

For the time-discretization error, we use the convergence result established in \textcite{HuYaozhong2011MCFB}. Indeed, Assumption~\ref{main assumption convergence rate} is stronger than what is needed for this result and implies Assumption~2.2 therein with $p=2$, together with the time-regularity condition~(3.4). Their result therefore yields the following theorem.

\begin{theorem}\label{Theorem main convergence rate}
    Let Assumption \ref{main assumption convergence rate} be satisfied. Then the BSDE \eqref{BSDE} is well-posed. Moreover, there exist constants $C,\eta>0$, independent of the partition $\pi$, such that
    \begin{gather*}
        \max_{0 \leq i \leq m-1}
        \mathbb{E}\Big(
        \sup_{t \in [t_i,t_{i+1}]}
        \abs{Y_t-Y_t^{\pi}}^{2}
        \Big)
        +
        \mathbb{E}\Big(
        \int_{0}^{T}
        \abs{Z_t-Z_t^{\pi}}^{2}\,dt
        \Big)
        \leq
        C\abs{\pi}
    \end{gather*}
    whenever $\abs{\pi}<\eta$. The constants $C$ and $\eta$ may depend on $T$, $L$, and the constants and norms appearing in Assumption \ref{main assumption convergence rate}.
\end{theorem}

\subsection{Error given by the truncation of the chaos}
\label{subsection: Error given by the truncation of the chaos}

We now focus on the error introduced by truncating the chaos decomposition at each step of the backward recursion. For each $i\in\{1,\dots,m\}$, let
\[
    \mathcal{G}_{i}^{M}
    :=
    \sigma\big(B(h_1^{i}),\dots,B(h_{M(i)}^{i})\big),
    \qquad
    f_i^{M}(y,z)
    :=
    \mathbb{E}\big(f(t_i,y,z)\mid\mathcal{G}_{i}^{M}\big).
\]
Choosing the Lipschitz version of the conditional expectation, $f_i^{M}$ is globally Lipschitz in $(y,z)$ with the same Lipschitz constant $[f]_L$ as $f$.

By the formulas in Subsection~\ref{subsection: Wiener chaos decomposition formulae},
the random variables $Y_{t_i}^{\pi,\theta}$ and
$\sprod{Z}_{i}^{\pi,\theta}$ are $\mathcal{G}_{i}^{M}$-measurable. We therefore set
\[
    f_i^{M,\pi,\theta}
    :=
    f_i^{M}\big(
        Y_{t_i}^{\pi,\theta},
        \sprod{Z}_{i}^{\pi,\theta}
    \big).
\]
Since $Y_{t_i}^{\pi,\theta}$ has a finite-dimensional chaos decomposition of order at most $P$,
\[
    \mathcal{C}_{P}^{\overline{\pi}_i}(F_i^{\pi,\theta})
    =
    Y_{t_i}^{\pi,\theta}
    +
    \Delta_i
    \mathcal{C}_{P}
    \big(f_i^{M,\pi,\theta}\big).
\]
Hence, the truncation error at each step consists of the projection error of the generator onto $\mathcal{G}_{i}^{M}$ and the truncation of its chaos order at level $P$.

We first compare the numerical scheme presented in
Subsection~\ref{Section Numerical Scheme} with the implementation outlined in
Section~\ref{section: Description of the algorithm}, assuming that the chaos coefficients can be computed exactly. The next lemma reduces the resulting error to the truncation of the terminal condition and the two errors identified above.

\begin{lemma}\label{lemma error implementation}
    There exist constants $C,\eta>0$, independent of $\pi$, $P$ and $M$, such that
  \begin{flalign*}
    &\max_{0 \leq i \leq m-1}
    \mathbb{E}\Big(
        \sup_{t\in[t_i,t_{i+1}]}
        \abs{Y_{t}^{\pi} - Y_{t}^{\pi,\theta}}^{2}
    \Big)
    +
    \mathbb{E}\Big(
        \int_{0}^{T}
        \abs{Z_{u}^{\pi} - Z_{u}^{\pi,\theta}}^{2}\,du
    \Big)
    \\
    &\leq
    C\left\{
        \mathbb{E}\Big(
            \abs{\xi-\mathcal{C}_{P}^{\overline{\pi}_m}(\xi)}^{2}
        \Big)
        +
        \sum_{i=1}^{m}
        \Delta_i\,
        \mathbb{E}\left(
            \abs{f_i^{\pi,\theta}-f_i^{M,\pi,\theta}}^{2}
            +
            \abs{
                f_i^{M,\pi,\theta}
                -
                \mathcal{C}_{P}(f_i^{M,\pi,\theta})
            }^{2}
        \right)
    \right\},
\end{flalign*}
    whenever $|\pi|<\eta$. The constants $C$ and $\eta$ may depend on $[f]_{L}$, $T$ and $L$.
\end{lemma}

The first term in the sum is specific to the randomness of the generator and measures the error introduced by its projection onto the retained Gaussian variables. The Malliavin regularity in Assumption~\ref{main assumption convergence rate}(iii) gives the same type of basis-truncation rate as Lemma~\ref{lemma approximation holder}.

\begin{lemma}\label{lemma generator basis truncation}
Let Assumption~\ref{main assumption convergence rate} hold. Then
\[
    \sum_{i=1}^{m}\Delta_i
    \mathbb{E}\big(
        \abs{f_i^{\pi,\theta}-f_i^{M,\pi,\theta}}^2
    \big)
    \leq
    C M^{-2\beta_f},
\]
where $C>0$ is independent of $\pi$, $P$ and $M$.
\end{lemma}

It remains to control the truncation of the chaos order. Lemma~\ref{lemma truncation order chaos} with $k=1$ yields
\begin{gather*}
    \mathbb{E}\big(
        \abs{
            f_i^{M,\pi,\theta}
            -
            \mathcal{C}_{P}(f_i^{M,\pi,\theta})
        }^{2}
    \big)
    \leq
    \frac{
        \mathbb{E}\big(
            \norm{D f_i^{M,\pi,\theta}}_{\mathbb{H}}^{2}
        \big)
    }{P+1}.
\end{gather*}
However, the numerator still depends on $P$ and $M$, since
$Y_{t_i}^{\pi,\theta}$ and $\sprod{Z}_{i}^{\pi,\theta}$ are themselves defined recursively through truncated chaos decompositions. Thus, to obtain a uniform estimate, we first establish energy-type bounds for their Malliavin derivatives.

\begin{lemma}\label{lem:malliavin-energy}
Let Assumption~\ref{main assumption convergence rate} hold. Then there exist constants $C,\eta>0$, independent of $\pi$, $P$ and $M$, such that
    \begin{gather}\label{eq:energy-estimates-Malliavin-Y-Z}
        \max_{1 \leq i \leq m}
        \mathbb{E}\big(
            \norm{D Y_{t_i}^{\pi,\theta}}_{\mathbb{H}}^{2}
        \big)
        +
        \sum_{i=1}^{m}\Delta_i
        \mathbb{E}\big(
            \norm{D\sprod{Z}_{i}^{\pi,\theta}}_{\mathbb{H}}^{2}
        \big)
        \leq
        C\left(
            1+
            \mathbb{E}\big(
                \norm{D\xi}_{\mathbb{H}}^{2}
            \big)
        \right),
    \end{gather}
    and
    \begin{gather}\label{eq:energy-estimates-malliavin-f}
        \sum_{i=1}^{m}\Delta_i
        \mathbb{E}\big(
            \norm{D f_i^{M,\pi,\theta}}_{\mathbb{H}}^{2}
        \big)
        \leq
        C\left(
            1+
            \mathbb{E}\big(
                \norm{D\xi}_{\mathbb{H}}^{2}
            \big)
        \right),
    \end{gather}
    whenever $|\pi|<\eta$. The constants $C$ and $\eta$ depend only on $L$, $T$ and the constants appearing in Assumption~\ref{main assumption convergence rate}.
\end{lemma}

\subsection{Error of the Monte Carlo approximation of the coefficients}\label{subsection: error montecarlo}

The goal of this subsection is to analyse the error due to the Monte Carlo approximation of the chaos coefficients described in Subsection~\ref{subsection: MC approximation}. To simplify the notation, let us write
\begin{gather*}
    \widehat{f}_{i}^{\pi,\theta,n}
    =
    f^{i,n}(t_i,\widehat{Y}_{t_i}^{\pi,\theta,n},\sprod{\widehat{Z}}_i^{\pi,\theta,n}),
    \quad
    \widehat{F}_{i}^{\pi,\theta,n}
    =
    \widehat{Y}_{t_i}^{\pi,\theta,n}
    +
    \Delta_i\widehat{f}_{i}^{\pi,\theta,n}.
\end{gather*}
Additionally, whenever we omit $n$ from the superscript, we mean that $n=0$, e.g. $\widehat{\mathcal{C}}_{P}^{\overline{\pi}_i,0}(\cdot)=\widehat{\mathcal{C}}_{P}^{\overline{\pi}_i}(\cdot)$, $\widehat{f}_{i}^{\pi,\theta}=\widehat{f}_{i}^{\pi,\theta,0}$, $\widehat{Y}_{t_i}^{\pi,\theta}=\widehat{Y}_{t_i}^{\pi,\theta,0}$, and so on.

For a square-integrable random variable $F$ on the enlarged probability space, define
\begin{gather*}
    \mathbf{V}_i(F)
    \coloneqq
    \sum_{k=0}^{P}
    \sum_{\substack{a\in\mathcal{A}_{M(i)}\\|a|=k}}
    a!\,\overline{\mathbb{V}}\big(F\times H_a^{(i,0)}\big),
\end{gather*}
where $\overline{\mathbb{V}}(\cdot)$ denotes the variance under $\overline{\mathbb{P}}$. Notice that $\mathbf{V}_i(F)$ may be infinite for a general square-integrable random variable $F$.

\begin{lemma}\label{lemma error implementation MC}
Let Assumption~\ref{main assumption convergence rate} hold and set $\Lambda_{P,M}\coloneqq 9^P\binom{P+M}{M}$. Then there exist constants $C,\eta>0$, independent of $\pi$, $P$, $M$ and $N$, such that, whenever $\abs{\pi}<\eta$,
\begin{flalign}
\max_{0\leq i\leq m-1}
\overline{\mathbb E}\Big(
    \sup_{t\in[t_i,t_{i+1}]}
    & \abs{Y_t^{\pi,\theta}-\widehat Y_t^{\pi,\theta}}^2
\Big)
+ 
\overline{\mathbb E}\Big(
    \int_0^T
    \abs{Z_u^{\pi,\theta}-\widehat Z_u^{\pi,\theta}}^2\,du
\Big)\notag\\
&\leq
\frac{C}{N}
\exp\left(
    \frac{C\Lambda_{P,M}}{N}
\right)
\left[
    \mathbf V_m(\xi)
    +
    \Lambda_{P,M}\abs{\pi}
\right].
\label{error given by MC}
\end{flalign}
Moreover,
\[
    \mathbf V_m(\xi)
    \leq
    3^P\binom{P+M}{M}
    \norm{\xi}_{L^q}^2.
\]
In particular, if $N\geq\Lambda_{P,M}$, then the right-hand side of \eqref{error given by MC} is bounded by $C\Lambda_{P,M}/N$. The constant $C$ may depend on $T$, $L$, $q$, and the constants and norms appearing in Assumption~\ref{main assumption convergence rate}.
\end{lemma}

\begin{remark}
The factor $\Lambda_{P,M}$ grows rapidly with $P$ and $M$, as the number of chaos coefficients to be estimated increases. Hence, the number of Monte Carlo samples $N$ should be increased accordingly to keep the Monte Carlo error small.
\end{remark}

\subsection{Main error estimate}

We conclude this section by combining the previous results to obtain the convergence rate of the fully implementable scheme.

\begin{theorem}\label{thm:main-conv-rate}
Let Assumption~\ref{main assumption convergence rate} hold, let $\theta=(P,M)\in\mathbb N^2$ and $N\in\mathbb N$. Then there exist constants $C,\eta>0$ and $C_{P,M}>0$ such that
\begin{flalign*}
&\max_{0\leq i\leq m-1}
\overline{\mathbb E}\Big(
    \sup_{t\in[t_i,t_{i+1}]}
    |Y_t-\widehat Y_t^{\pi,\theta}|^2
\Big)
+
\overline{\mathbb E}\Big(
    \int_0^T
    |Z_u-\widehat Z_u^{\pi,\theta}|^2\,du
\Big)\\
&\qquad\leq
C\left(
    |\pi|
    +
    \frac{1+\mathbb E\big(\|D\xi\|_{\mathbb H}^2\big)}{P+1}
    +
    M^{-2\min\{\beta_\xi,\beta_f\}}
    +
    \frac{C_{P,M}}{N}
\right)
\end{flalign*}
whenever $|\pi|<\eta$. The constants $C$ and $\eta$ are independent of $\pi$, $P$, $M$ and $N$, whereas $C_{P,M}$ is independent of $\pi$ and $N$ and may depend on $P$ and $M$. 
\end{theorem}

A few remarks are in order.

\begin{remark}
Compared with the Picard chaos scheme of \textcite{BriandPhilippe2014SOBB}, our convergence analysis differs in three main respects:
\begin{itemize}
    \item the regularity assumptions are of fixed order, whereas in \textcite{BriandPhilippe2014SOBB} the required Malliavin differentiability and smoothness orders increase with the Picard and chaos truncation levels.
    \item the discretization error of the time integrals is included in our estimates, while these integrals are assumed to be computed exactly in the analysis of \textcite{BriandPhilippe2014SOBB}.
    \item the constant $C$ in Theorem~\ref{thm:main-conv-rate} is independent of the approximation parameters, with the dependence of the Monte Carlo error on $P$ and $M$ isolated in the factor $C_{P,M}$.
\end{itemize}
In particular, our convergence result is stronger and holds under substantially more general regularity assumptions.
\end{remark}

\begin{remark}
We now turn to regression-based implementations of the Euler scheme in the finite-dimensional Markovian setting. At each backward step, the conditional expectation defining $Y_{t_i}^\pi$ is approximated by regression. The resulting error enters
\[
    F_i^\pi
    =
    Y_{t_i}^\pi
    +
    \Delta_i f(t_i,Y_{t_i}^\pi,\sprod{Z}_i^\pi)
\]
directly, rather than being multiplied by $\Delta_i$. Hence, as the partition is refined, these errors may accumulate over an increasing number of time steps and cannot in general be controlled by the usual discrete Gronwall argument unless the regression accuracy is refined accordingly; see, for example, \textcite{HureCome2020Dbsf,GermainPhamWarin2022Approximation}. By contrast, our implementation avoids this issue, since the corresponding conditional expectations are computed exactly once the finite chaos representations have been constructed.
\end{remark}

\section{Numerical Experiments}\label{section: numerical examples}

In this section, we present four numerical examples illustrating different aspects of the proposed method. Example~1 admits a closed-form solution and is used to verify numerically the convergence rates established in Section~\ref{section: Convergence analysis}. In Example~2, we compare our scheme and the Picard-based chaos method of \textcite{BriandPhilippe2014SOBB} against an independent benchmark provided by \textcite{BriandDelyonMemin2001}, and also compare their computational costs, observing a clear advantage of the Euler-based approach. Examples~3 and~4 illustrate applications to pricing and hedging under rough-volatility models. Except for Example~1, all the examples fall outside the classical finite-dimensional Markovian BSDE setting discussed in the Introduction; see, in particular, the discussion preceding \eqref{forward SDE}. We also remark that some of the examples fall outside the regularity assumptions of our convergence results and illustrate the performance of the method beyond its theoretical setting.

All experiments were conducted on a machine running Red Hat Enterprise Linux
8.10, equipped with an Intel(R) Xeon(R) Gold 5120 CPU and an NVIDIA A10 GPU
with 24 GB of memory. The algorithm was implemented in Python 3.9.19 using
PyTorch 2.3.0 for tensor computations and GPU acceleration. The code is
available in the associated GitHub repository.\footnote{\url{https://github.com/pere98diaz/An-Euler-scheme-for-BSDEs-using-the-Wiener-chaos-decomposition}}

\subsection{Example 1 -- Linear pricing and hedging in the Black--Scholes model}

We consider a one-dimensional Black--Scholes model under the physical probability measure $\mathbb{P}$, $dS_t=S_t(\mu\,dt+\sigma\,dB_t)$, $S_0>0$, together with the BSDE with linear driver
\begin{gather*}
    Y_t
    =
    \xi
    +
    \int_t^T
    \big(-rY_s-\frac{\mu-r}{\sigma} Z_s\big)\,ds
    -
    \int_t^T Z_s\,dB_s, \quad \xi
    =
    \exp\left(
        \frac{2}{T}\int_0^T \log S_t\,dt
    \right).
\end{gather*}
This BSDE describes the pricing by replication of the continuously monitored squared geometric-average claim $\xi$, where $Y$ represents its price process and $Z/(\sigma S)$ the corresponding replicating strategy; see, for example, \textcite{ElKarouiN.1997BSDE}.

One can check that this example satisfies Assumption~\ref{main assumption convergence rate}, with $\beta_\xi=1$ and, since the driver is deterministic, $\beta_f=1$. Hence, Theorem~\ref{thm:main-conv-rate} applies.

A closed-form expression for the solution can be obtained from the log-normal distribution of the continuously monitored geometric average; see, for example, \textcite{asian_zhang}. We use this expression as the reference solution in the numerical experiments.

We take $T=1$, $S_0=1$, $\mu=0.05$, $r=0.02$, and $\sigma=0.2$. For the convergence study in $m$, $M$, and $P$, we use the exact chaos coefficients of the terminal condition, which are available explicitly in this example, and report the empirical root mean-square error (RMSE) corresponding to the quantity on the left-hand side of Theorem~\ref{thm:main-conv-rate}.

\begin{figure}[t]
    \centering
    \begin{subfigure}{0.47\textwidth}
        \centering
        \includegraphics[width=\textwidth]{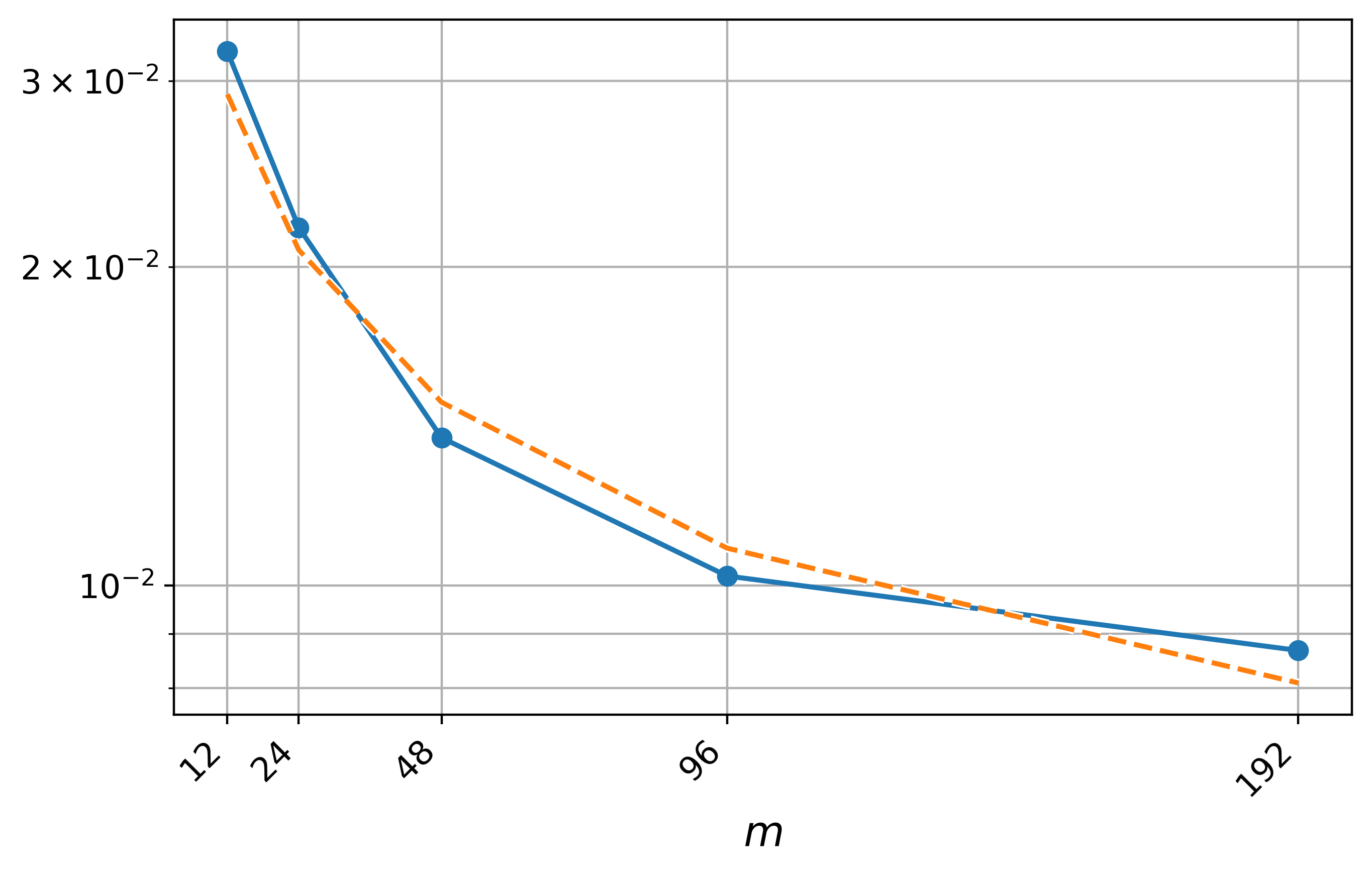}
        \caption{Varying $m$.}
        \label{fig:example1-m-rmse}
    \end{subfigure}
    \hfill
    \begin{subfigure}{0.47\textwidth}
        \centering
        \includegraphics[width=\textwidth]{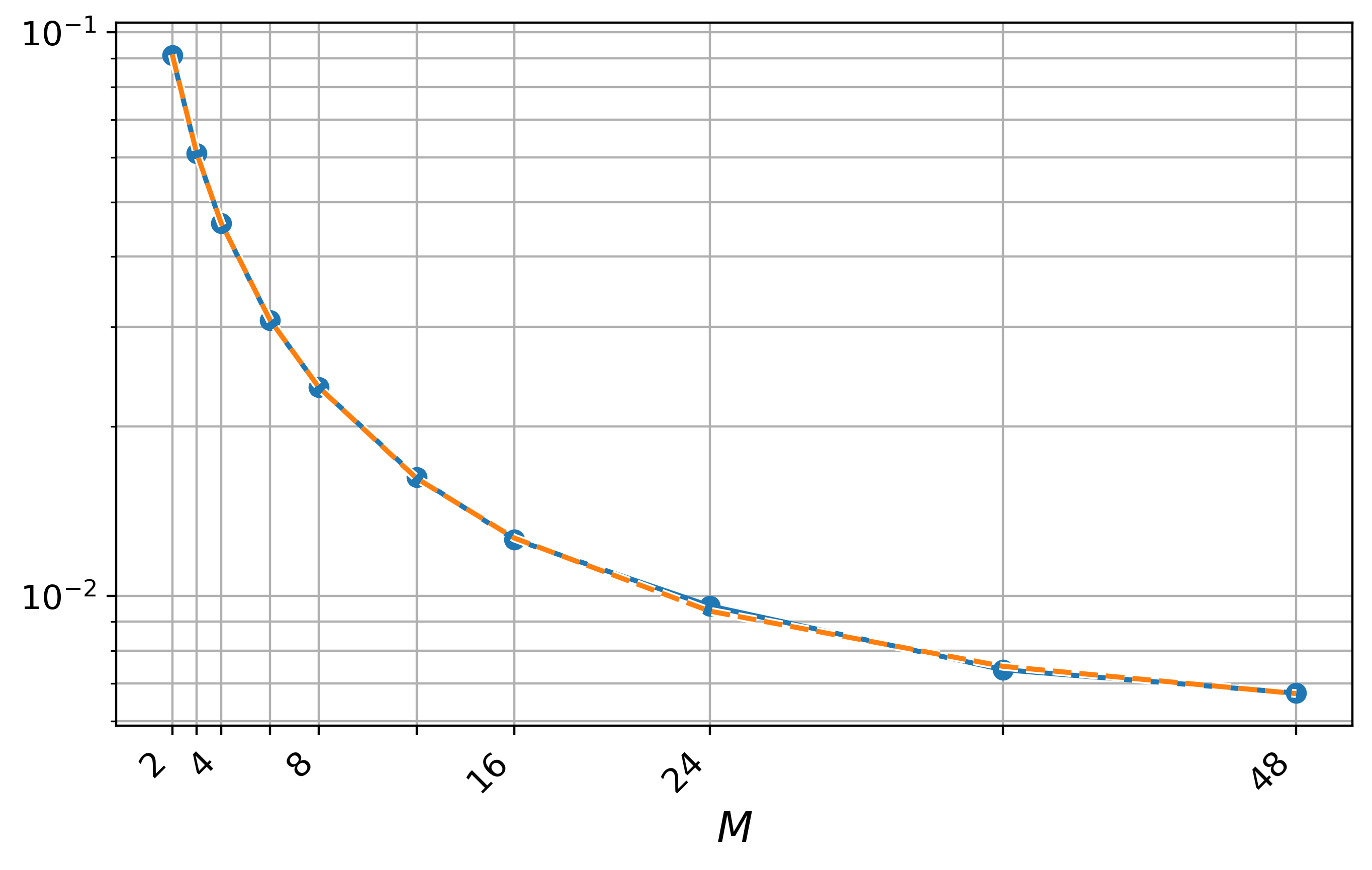}
        \caption{Varying $M$.}
        \label{fig:example1-M-rmse}
    \end{subfigure}

    \vspace{0.1em}

    \begin{subfigure}{0.47\textwidth}
        \centering
        \includegraphics[width=\textwidth]{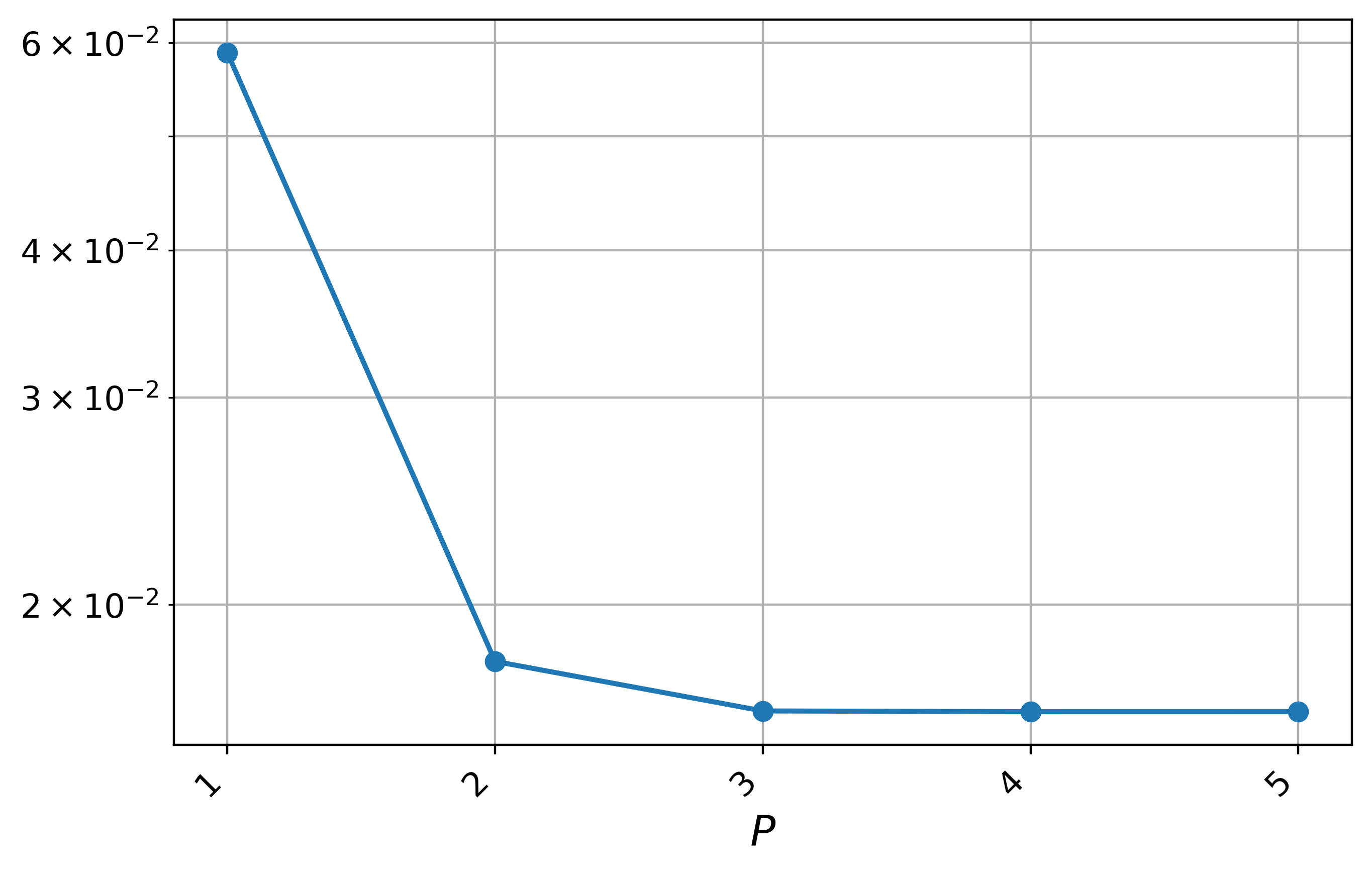}
        \caption{Varying $P$.}
        \label{fig:example1-P-rmse}
    \end{subfigure}
    \hfill
    \begin{subfigure}{0.47\textwidth}
        \centering
        \includegraphics[width=\textwidth]{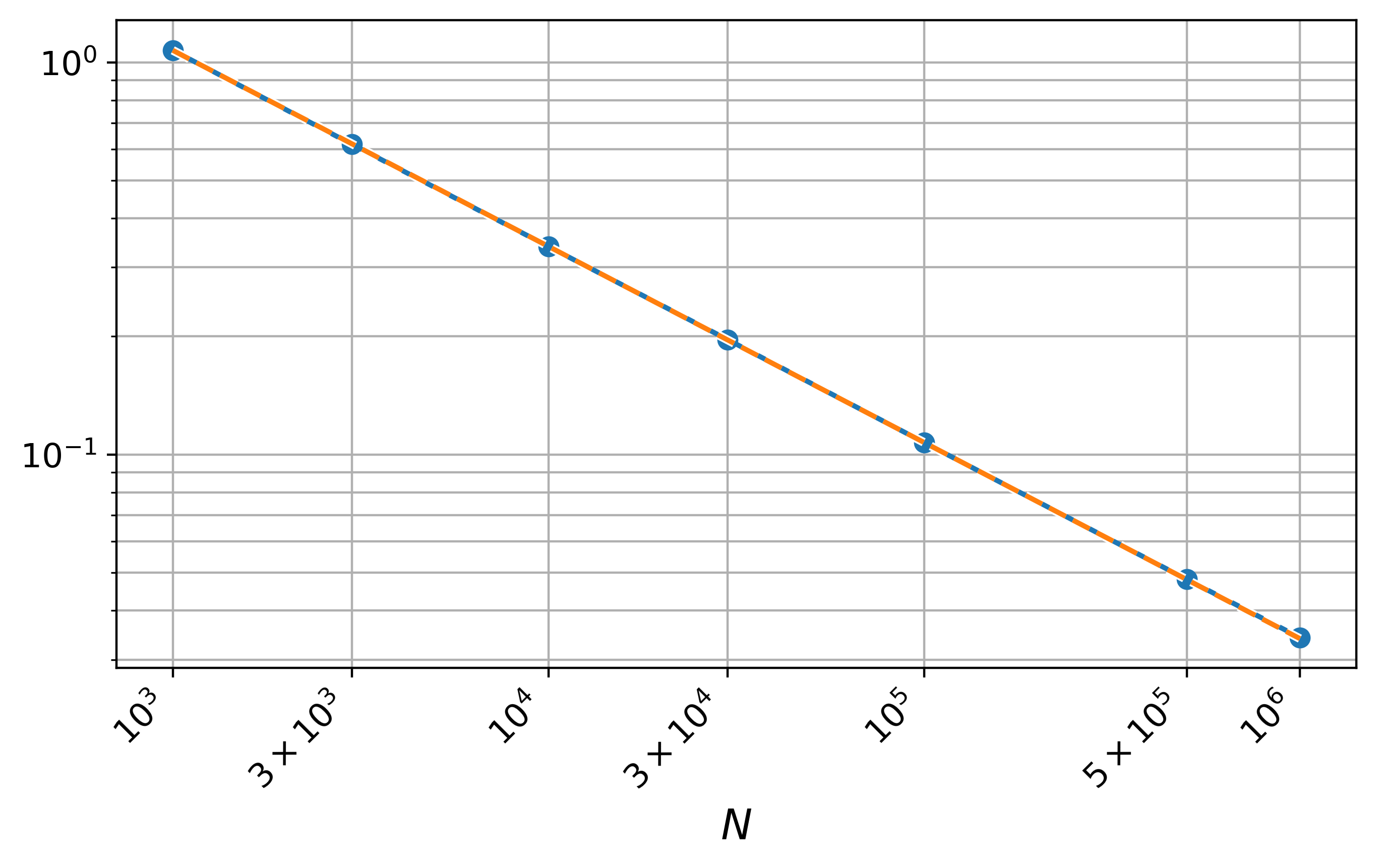}
        \caption{Varying $N$.}
        \label{fig:example1-N-rmse}
    \end{subfigure}

    \caption{Combined $(Y,Z)$ RMSE for Example~1 when varying $m$, $M$, $P$, and $N$. Dashed curves show the fits motivated by the theoretical rates in $m$ and $M$, and the $N^{-1/2}$ reference rate.}
    \label{fig:example1-convergence}
\end{figure}

Figure~\ref{fig:example1-convergence} shows the errors obtained by varying one approximation parameter at a time. For $m$, with $P=3$ and $M=24$ fixed, the fit $\sqrt{A/m+B}$ is consistent with the time-discretization term in Theorem~\ref{thm:main-conv-rate}; for $M$, with $P=3$ and $m=120$, the fit $\sqrt{AM^{-2}+B}$ is consistent with the basis-truncation term. Here, $A>0$ scales the contribution of the parameter being varied, while $B\geq0$ accounts for the error floor due to the remaining fixed parameters. When varying $P$, with $M=12$ and $m=120$, most of the improvement occurs from $P=1$ to $P=2$, after which the error reaches a plateau. Since the $(P+1)^{-1}$ term is only a generic upper bound and is not sharp for this terminal condition, we do not fit a rate in $P$.

To study the Monte Carlo approximation, we fix $P=3$, $M=12$, and $m=120$, and compare the solutions obtained with Monte Carlo and exact chaos coefficients over $100$ independent runs. The observed decay closely follows the expected $N^{-1/2}$ rate.

\begin{figure}[htbp]
    \centering
    \includegraphics[width=0.9\textwidth]{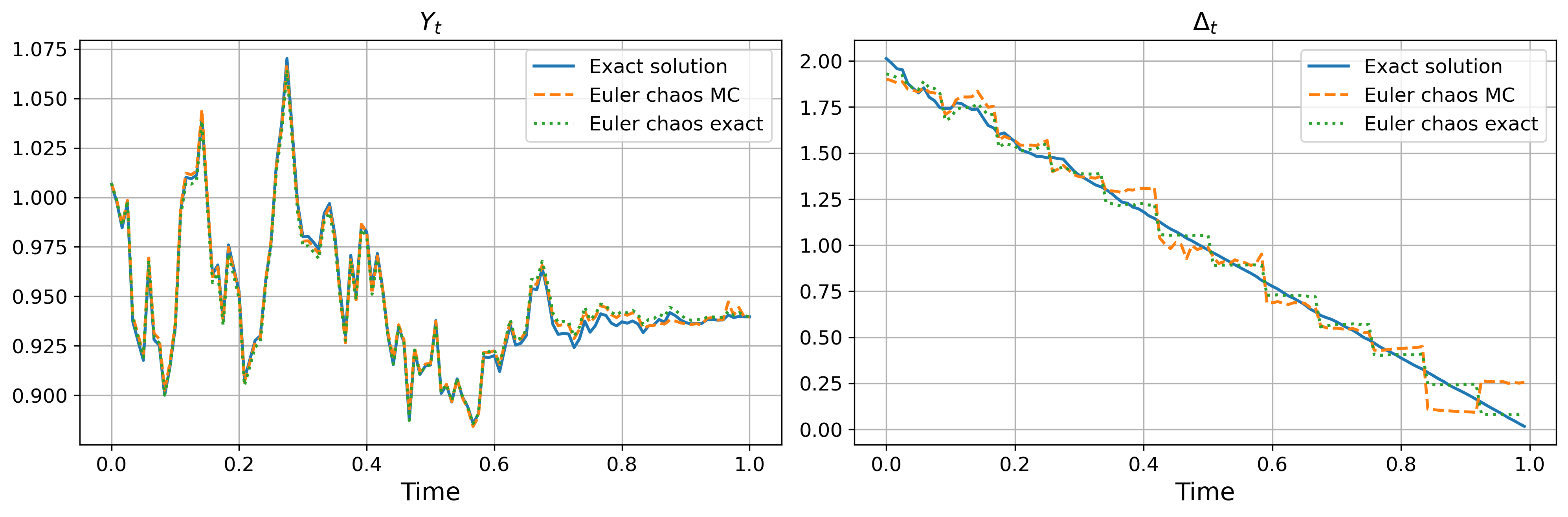}
    \caption{Sample paths of the exact price $Y$ and hedging strategy $\Delta=Z/(\sigma S)$, together with the Euler--chaos approximations obtained using Monte Carlo and exact chaos coefficients.}
    \label{fig:sample_paths_euler_mc_exact}
\end{figure}

Finally, Figure~\ref{fig:sample_paths_euler_mc_exact} illustrates the pathwise performance of the method for a single trajectory.

\subsection{Example 2 -- Non-linear BSDE with a path-dependent Volterra terminal condition}

We consider the one-dimensional BSDE with driver $f(t,y,z)=\cos(y+z)$ and terminal condition
\begin{gather*}
    \xi
    =
    \int_0^T (B_t^H)^2\,dt,
    \qquad
    B_t^H
    =
    \sqrt{2H}
    \int_0^t (t-s)^{H-\frac12}\,dB_s,
\end{gather*}
with $H=0.75$. The Volterra process $B^H$ is non-Markovian and does not admit a natural low-dimensional Markovian representation. We approximate the time integral defining $\xi$ by the trapezoidal rule on a uniform grid with $100$ intervals, which could be embedded into a $100$-dimensional state. Our method avoids this enlargement by constructing the chaos approximation directly from the one-dimensional driving Brownian motion.

As an independent benchmark, we use the Donsker-type approximation of \textcite{BriandDelyonMemin2001}, in which Brownian motion is replaced by a symmetric random walk and the resulting discrete BSDE is solved exactly on the binary tree. Extrapolating the values obtained for the finest discretizations gives
\begin{gather*}
    Y_0^{\mathrm{ref}}\approx 1.1360,
    \qquad
    Z_0^{\mathrm{ref}}\approx -0.3420.
\end{gather*}

We compare the Euler- and Picard-based chaos schemes in terms of both accuracy and computational cost. In the Picard-based scheme, the time integrals arising in the Picard iteration are discretized on a uniform grid. For consistency, we denote the number of time steps by $m$, since it plays a role analogous to $m$ in the Euler scheme. We use $N=10^6$ Monte Carlo samples and perform $100$ independent runs for each configuration. We vary one parameter at a time: $m$ with $(P,M)=(2,30)$, $M$ with $(P,m)=(2,120)$, and $P$ with $(M,m)=(15,120)$.

\begin{figure}[t]
    \centering
    \begin{subfigure}{0.47\textwidth}
        \centering
        \includegraphics[width=\textwidth]{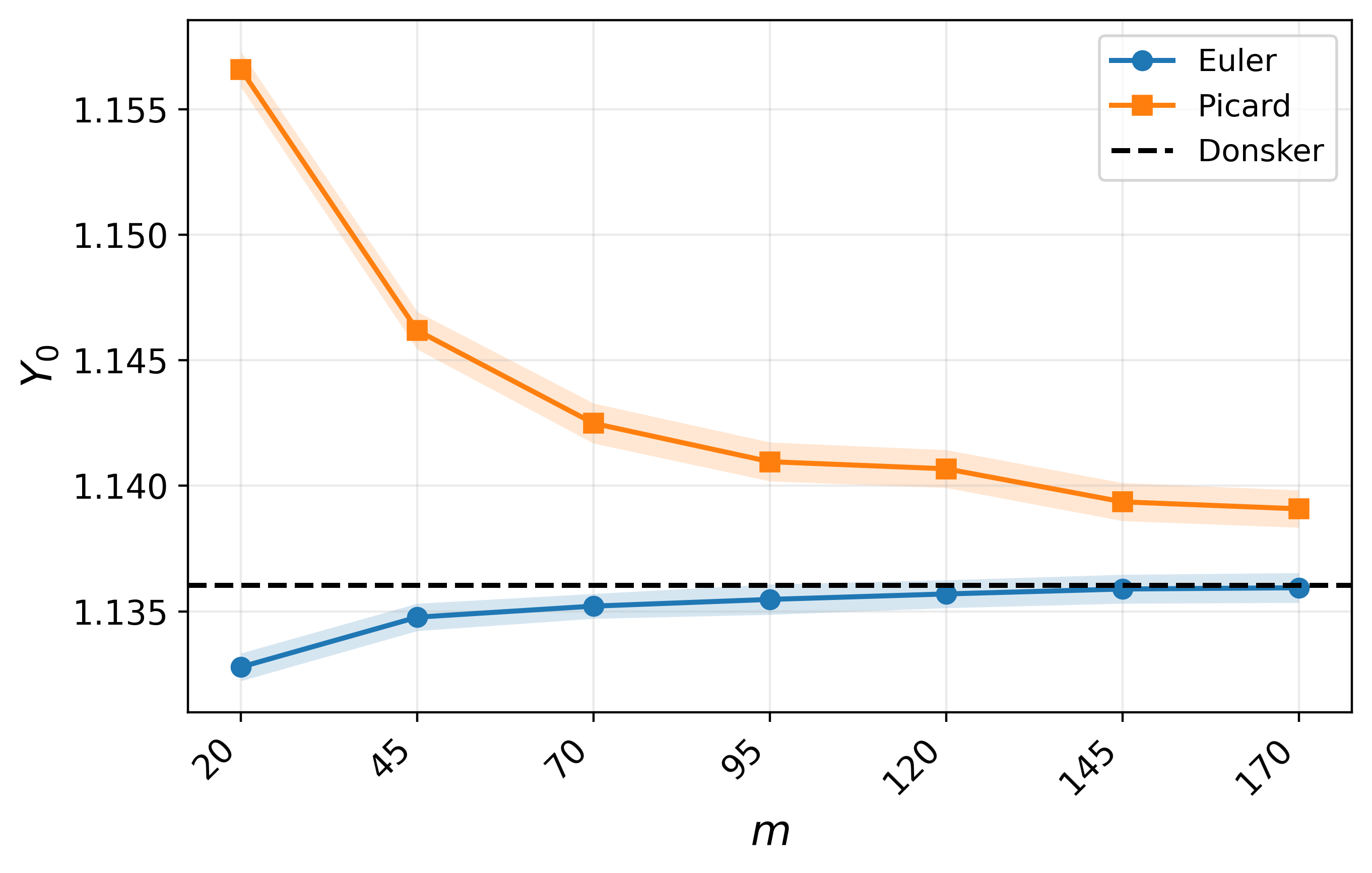}
        \caption{$Y_0$, varying $m$.}
        \label{fig:example2-m-y0}
    \end{subfigure}
    \hfill
    \begin{subfigure}{0.47\textwidth}
        \centering
        \includegraphics[width=\textwidth]{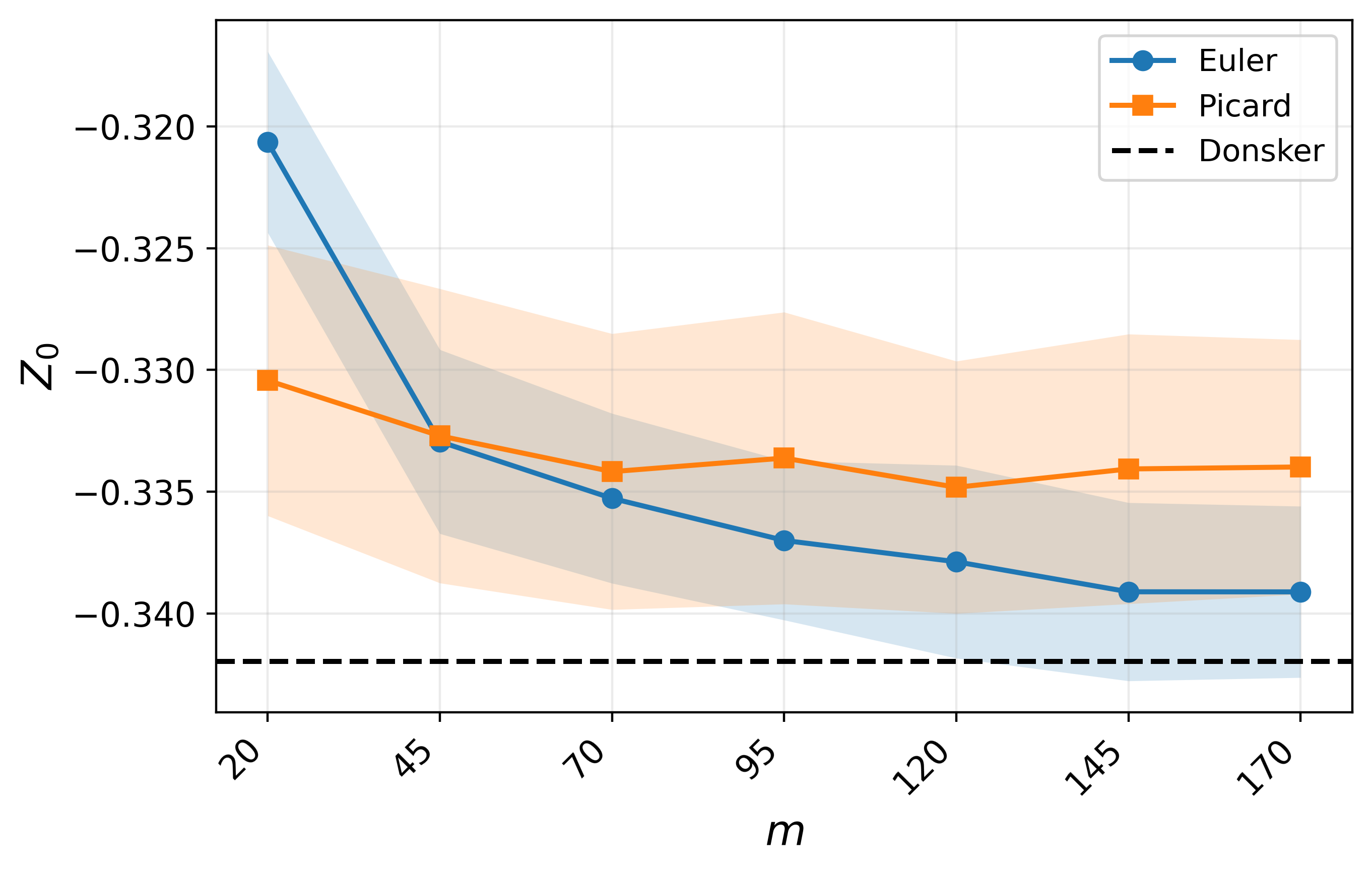}
        \caption{$Z_0$, varying $m$.}
        \label{fig:example2-m-z0}
    \end{subfigure}

    \vspace{0.1em}

    \begin{subfigure}{0.47\textwidth}
        \centering
        \includegraphics[width=\textwidth]{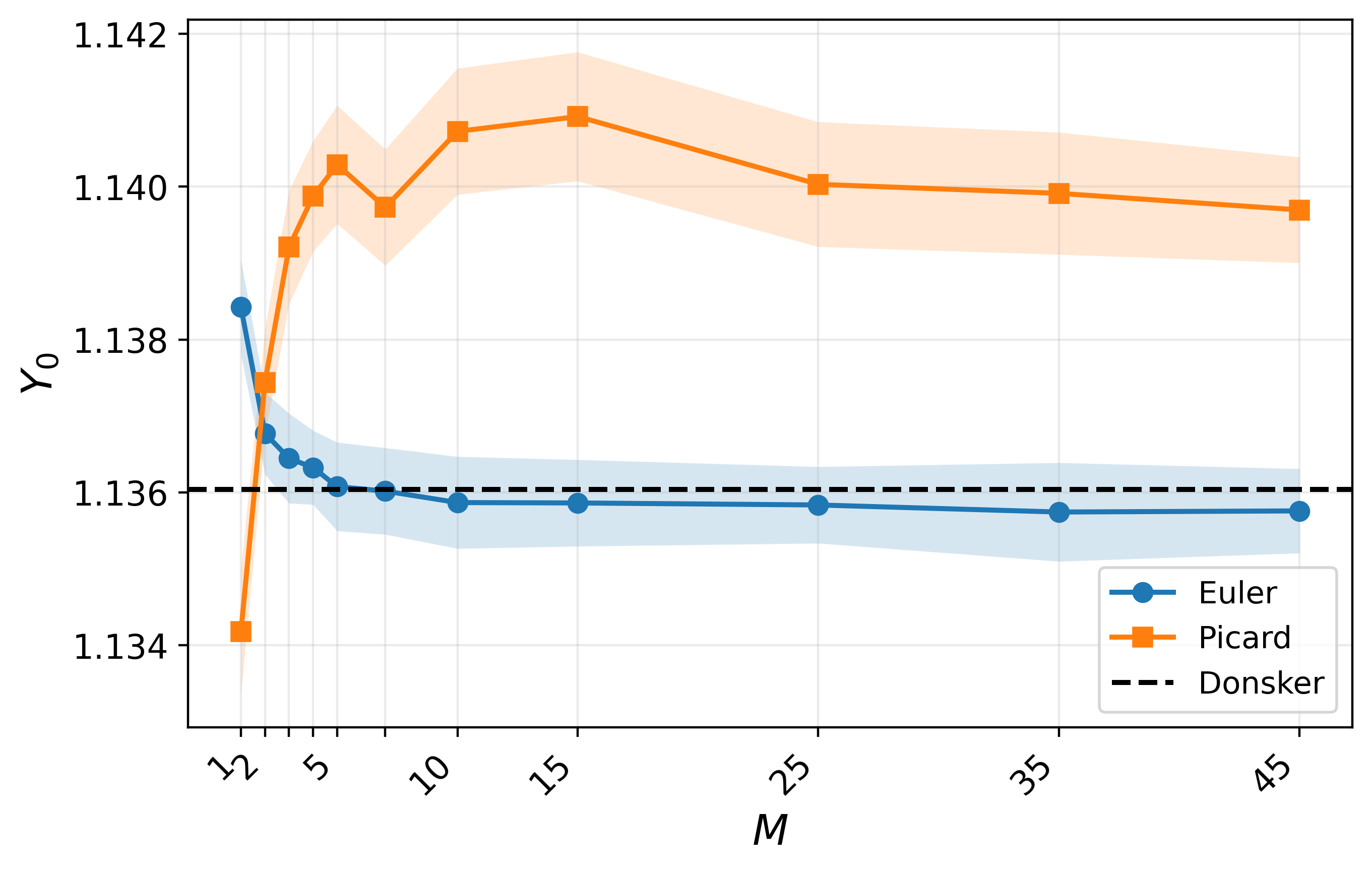}
        \caption{$Y_0$, varying $M$.}
        \label{fig:example2-M-y0}
    \end{subfigure}
    \hfill
    \begin{subfigure}{0.47\textwidth}
        \centering
        \includegraphics[width=\textwidth]{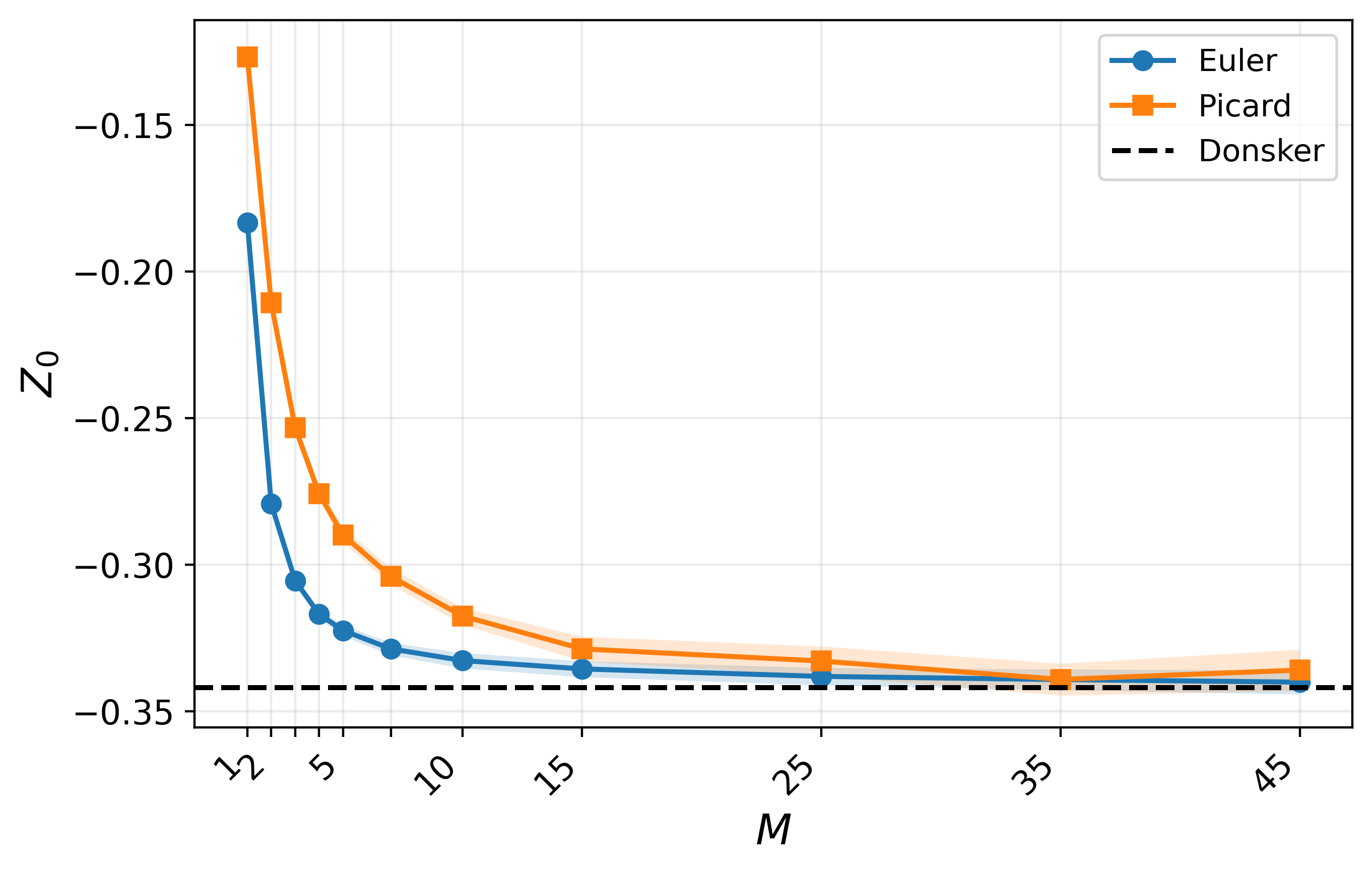}
        \caption{$Z_0$, varying $M$.}
        \label{fig:example2-M-z0}
    \end{subfigure}

    \vspace{0.1em}

    \begin{subfigure}{0.47\textwidth}
        \centering
        \includegraphics[width=\textwidth]{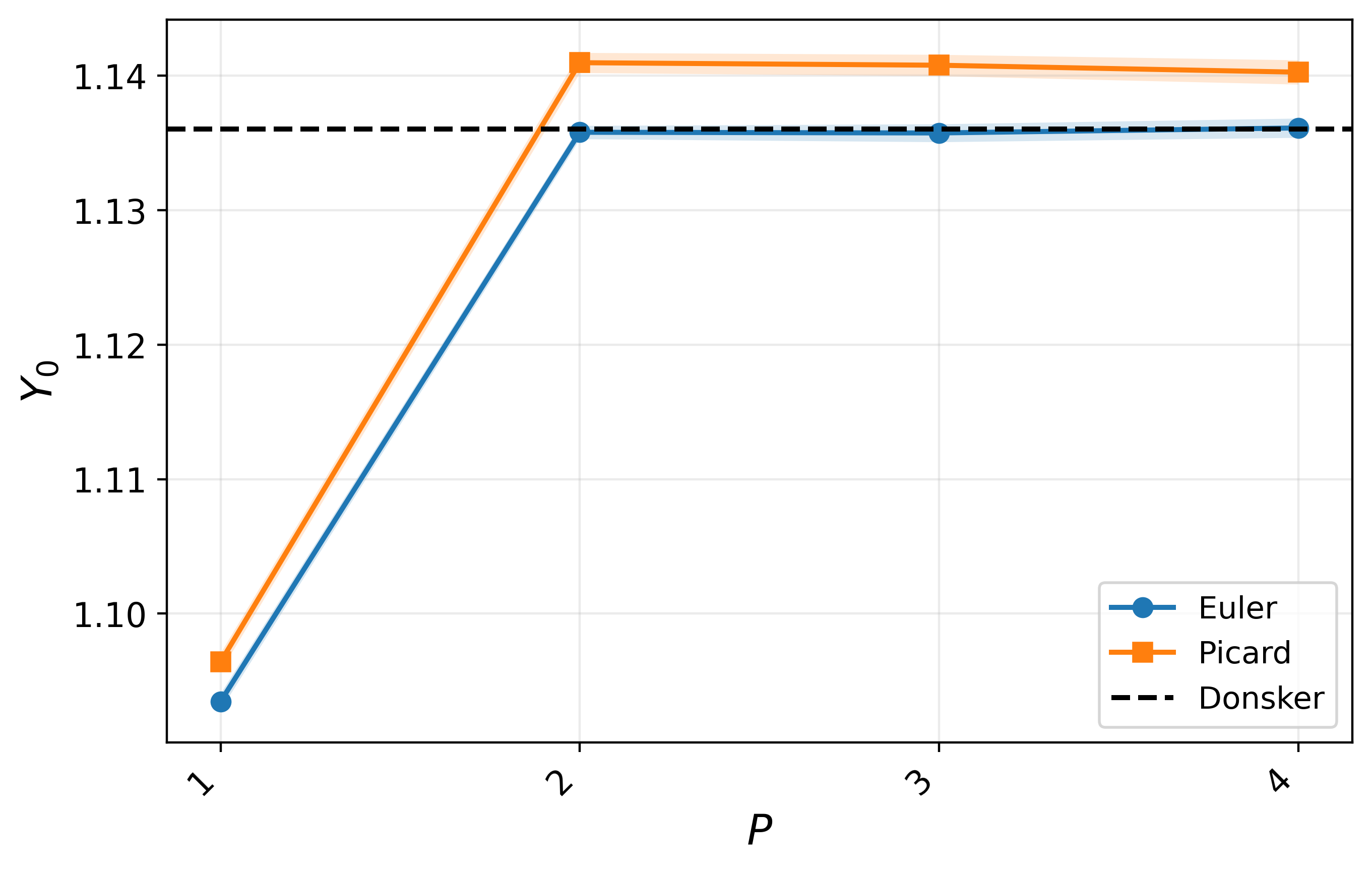}
        \caption{$Y_0$, varying $P$.}
        \label{fig:example2-P-y0}
    \end{subfigure}
    \hfill
    \begin{subfigure}{0.47\textwidth}
        \centering
        \includegraphics[width=\textwidth]{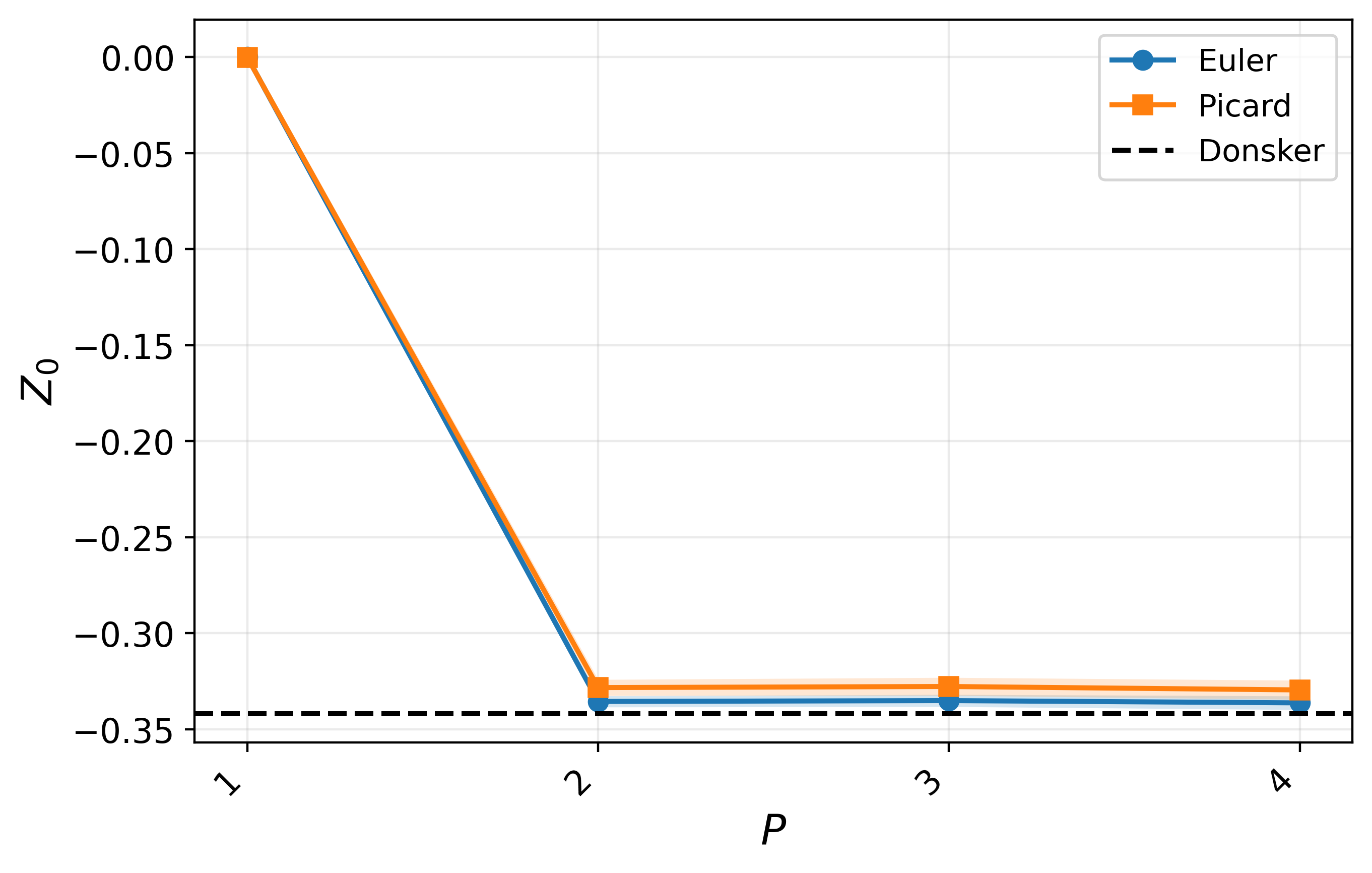}
        \caption{$Z_0$, varying $P$.}
        \label{fig:example2-P-z0}
    \end{subfigure}

    \caption{Initial-value approximations for Example~2. The solid curves show the empirical means over $100$ independent runs, the shaded regions show one standard deviation, and the dashed lines are the Donsker benchmarks.}
    \label{fig:example2-accuracy}
\end{figure}

Figure~\ref{fig:example2-accuracy} shows the approximations of $Y_0$ and $Z_0$. The Euler scheme is generally closer to the Donsker benchmark, particularly as $m$ and $M$ increase. When varying $P$, most of the improvement occurs from $P=1$ to $P=2$, after which both methods are essentially stable. This is consistent with the terminal condition being represented exactly in chaos order for $P\geq2$, although the non-linear driver may generate higher-order components.

\begin{figure}[t]
    \centering
    \begin{subfigure}{0.47\textwidth}
        \centering
        \includegraphics[width=\textwidth]{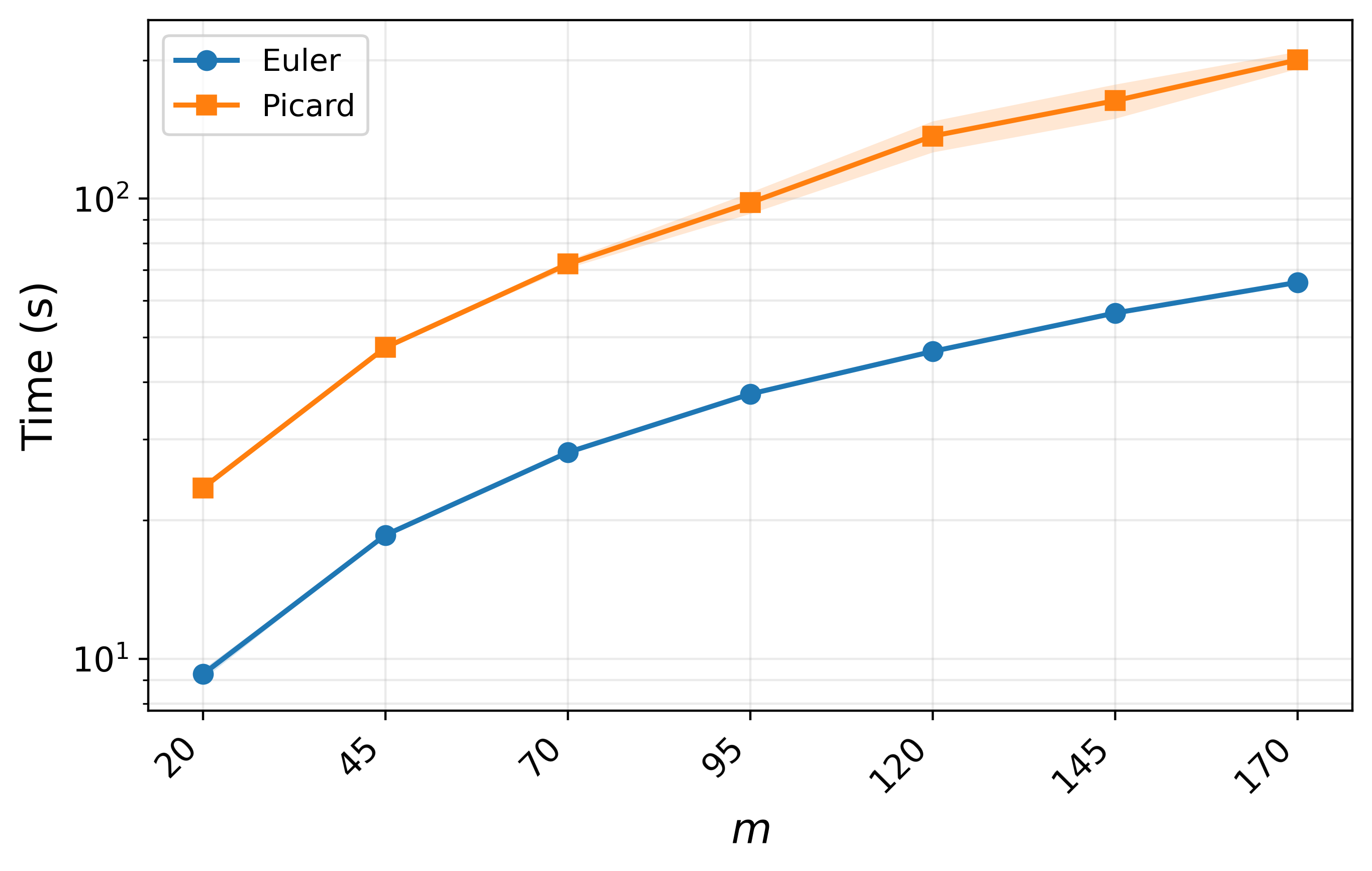}
        \caption{Time, varying $m$.}
        \label{fig:example2-m-runtime}
    \end{subfigure}
    \hfill
    \begin{subfigure}{0.47\textwidth}
        \centering
        \includegraphics[width=\textwidth]{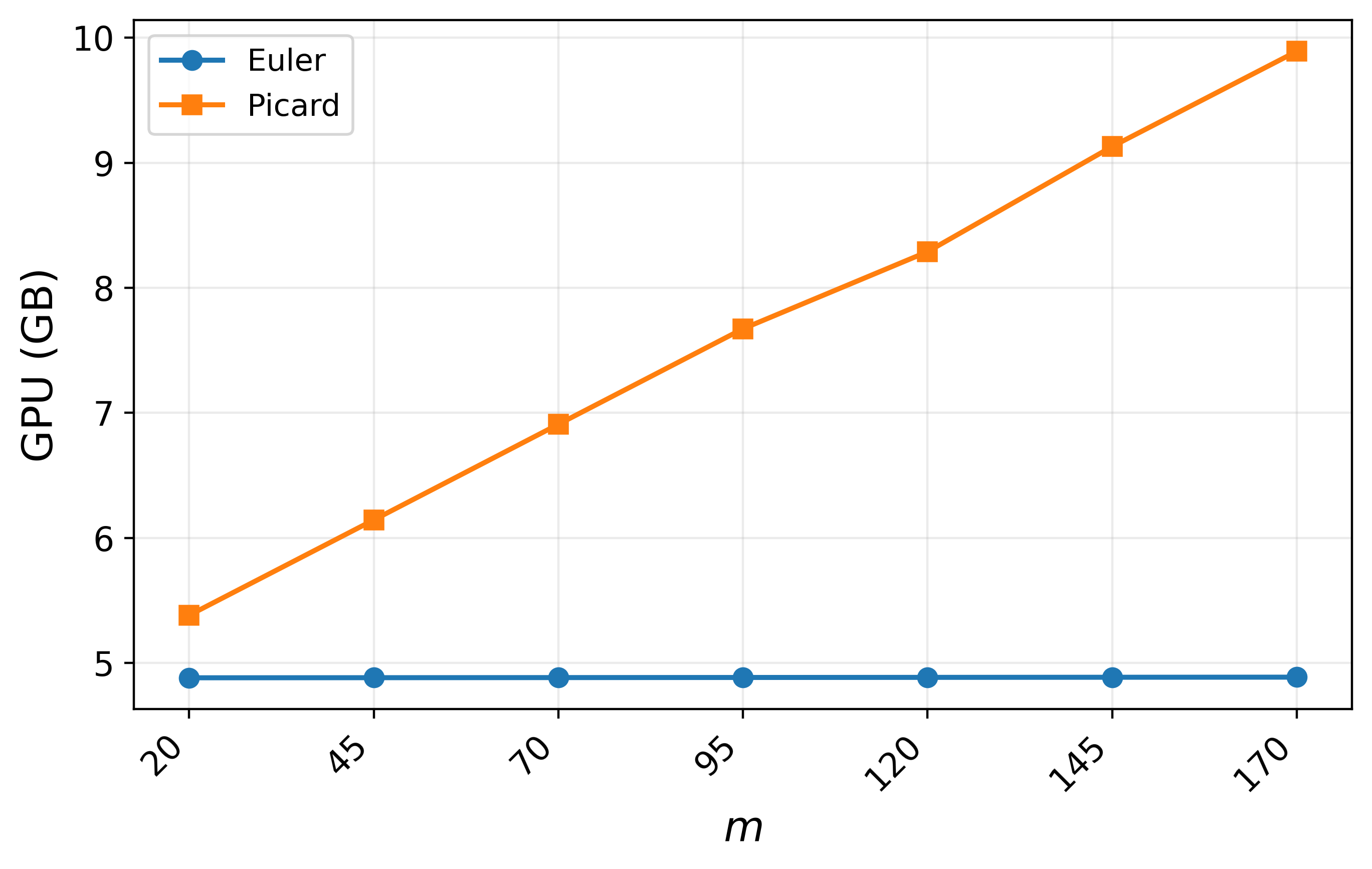}
        \caption{GPU memory, varying $m$.}
        \label{fig:example2-m-memory}
    \end{subfigure}

    \vspace{0.1em}

    \begin{subfigure}{0.47\textwidth}
        \centering
        \includegraphics[width=\textwidth]{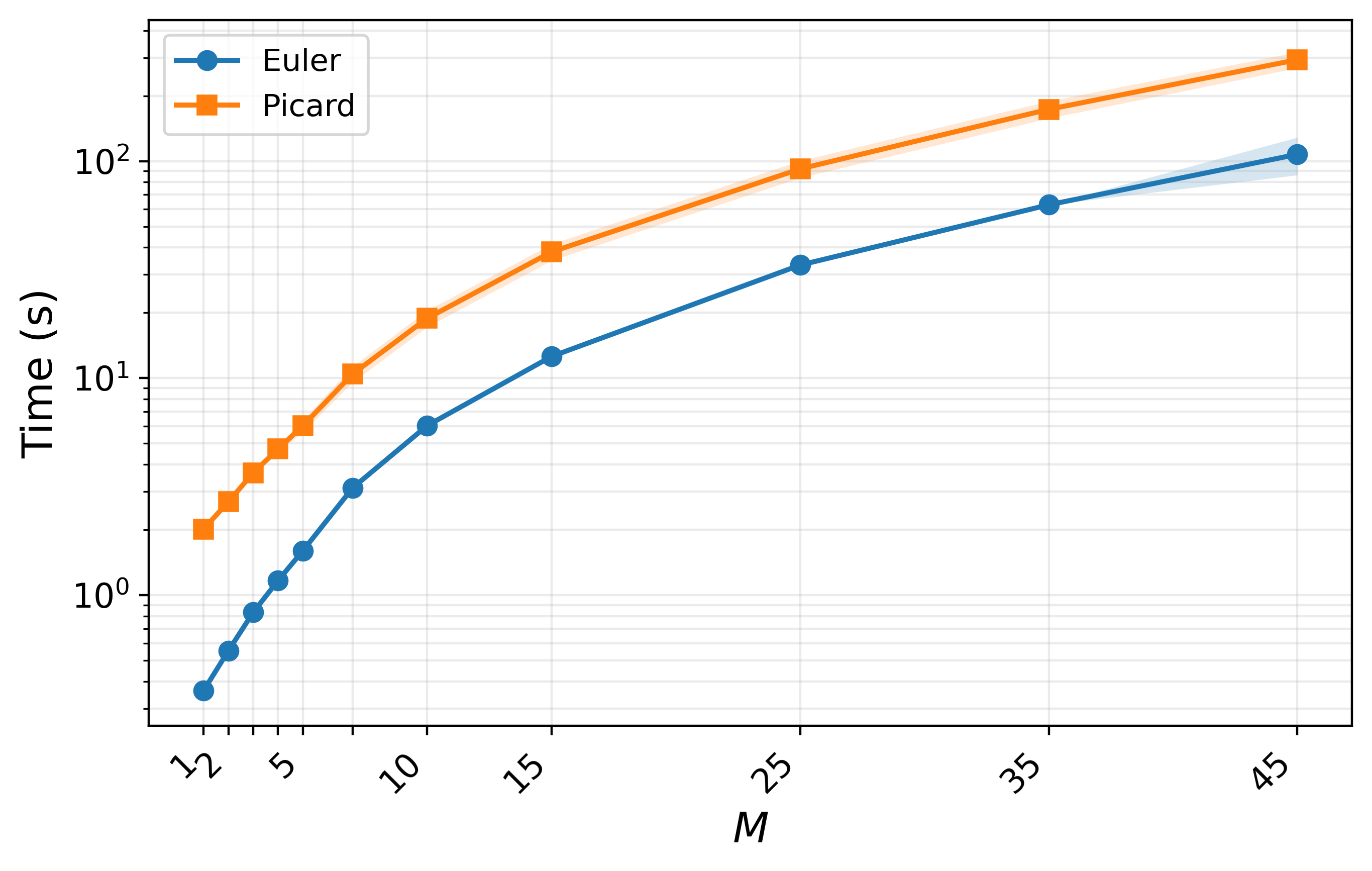}
        \caption{Time, varying $M$.}
        \label{fig:example2-M-runtime}
    \end{subfigure}
    \hfill
    \begin{subfigure}{0.47\textwidth}
        \centering
        \includegraphics[width=\textwidth]{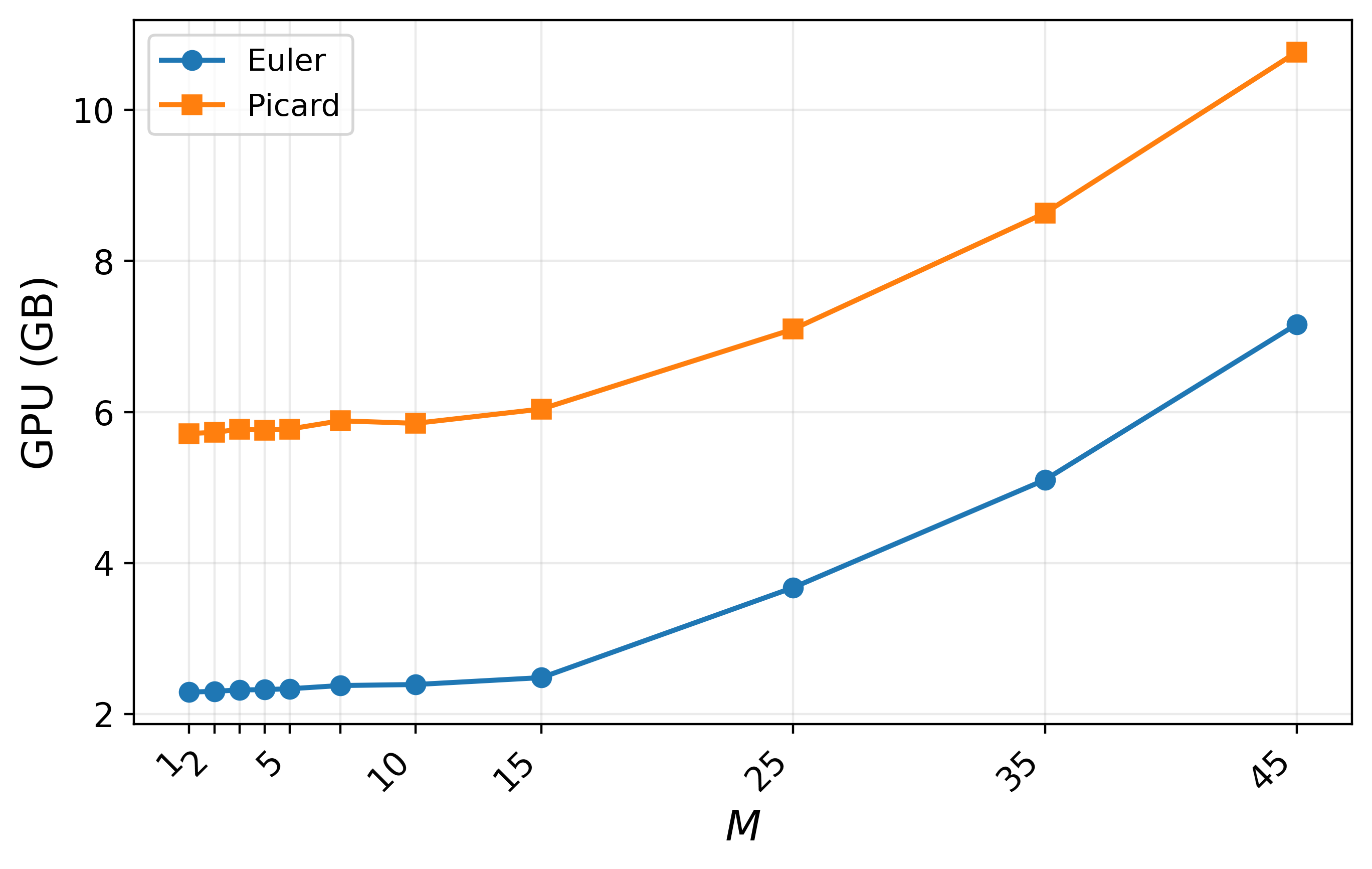}
        \caption{GPU memory, varying $M$.}
        \label{fig:example2-M-memory}
    \end{subfigure}

    \vspace{0.1em}

    \begin{subfigure}{0.47\textwidth}
        \centering
        \includegraphics[width=\textwidth]{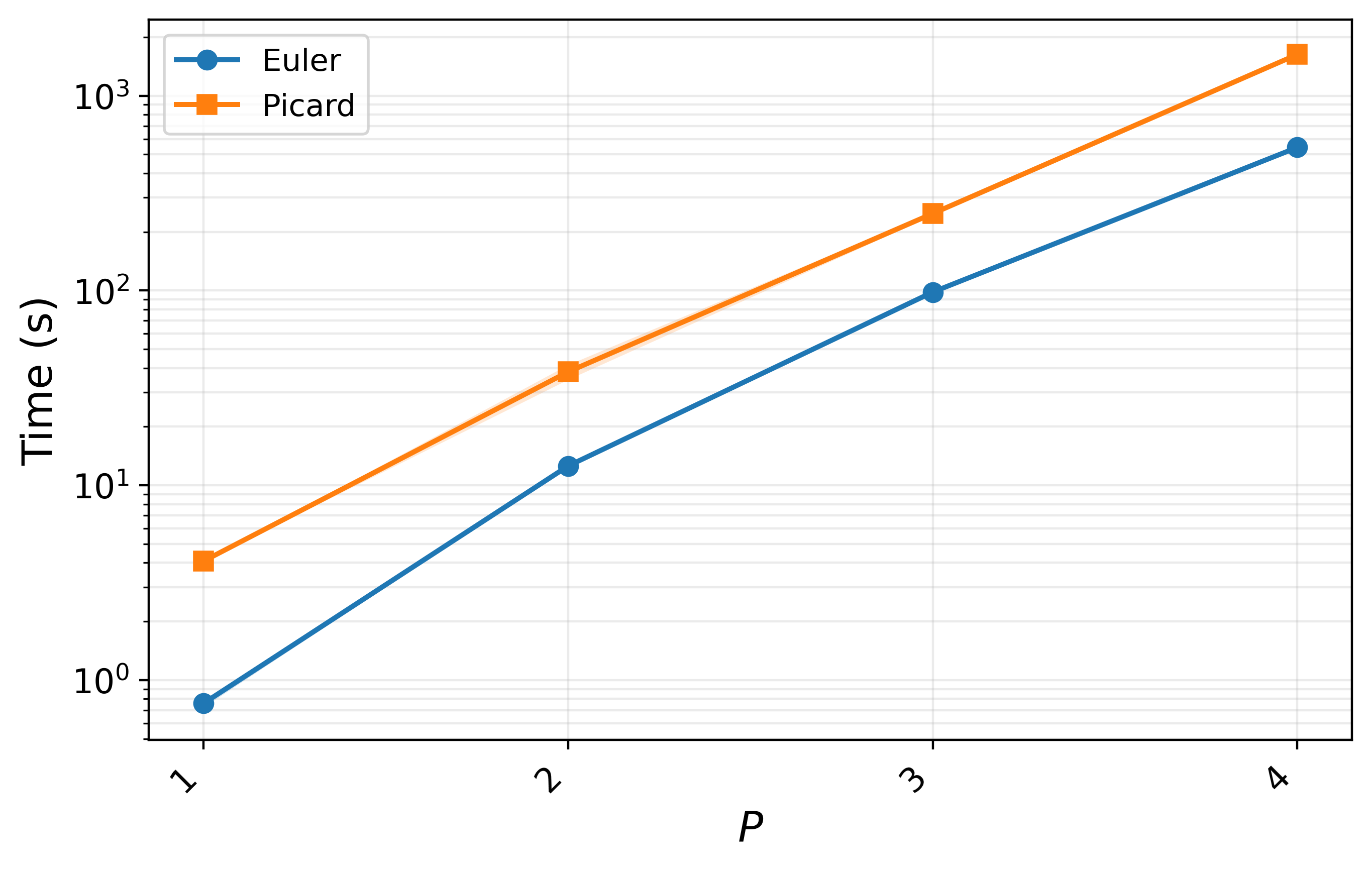}
        \caption{Time, varying $P$.}
        \label{fig:example2-P-runtime}
    \end{subfigure}
    \hfill
    \begin{subfigure}{0.47\textwidth}
        \centering
        \includegraphics[width=\textwidth]{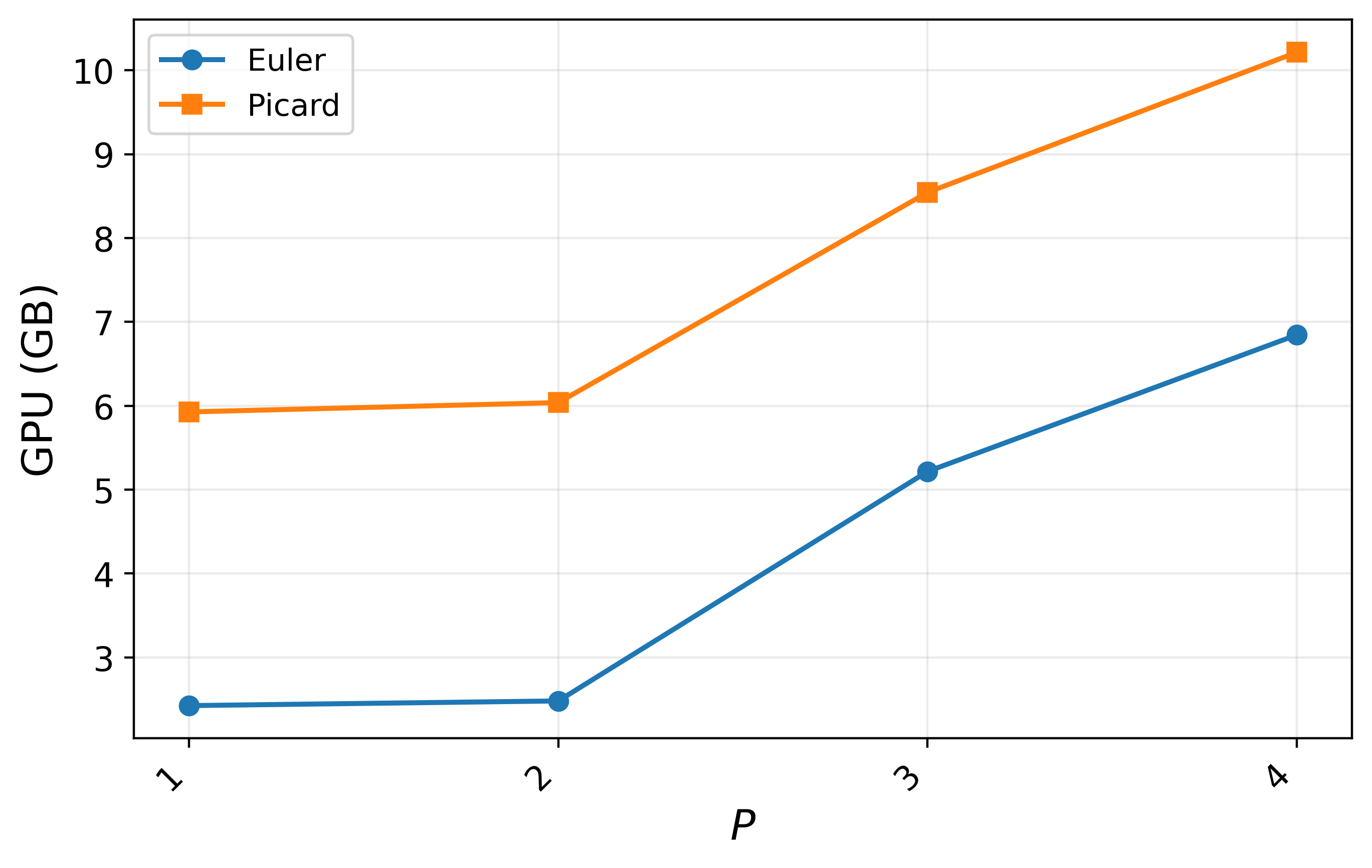}
        \caption{GPU memory, varying $P$.}
        \label{fig:example2-P-memory}
    \end{subfigure}

    \caption{Computational cost of the Euler and Picard schemes in Example~2. Left: runtime. Right: peak GPU memory.}
    \label{fig:example2-performance}
\end{figure}

As shown in Figure~\ref{fig:example2-performance}, the Euler scheme is consistently faster and requires less GPU memory than the Picard scheme across all three parameter sweeps, with the difference becoming more pronounced as the approximation parameters increase.

\subsection{Example 3 -- Linear pricing and hedging in the rough Bergomi model}

We consider the pricing and hedging of a European call option in a rough Bergomi market under the risk-neutral probability measure $\mathbb{Q}$. Let $B=(B^{1},B^{2})$ be a two-dimensional Brownian motion under $\mathbb{Q}$. We then consider
\begin{flalign*}
\frac{dS_t}{S_t}
&=
r\,dt
+
\sqrt{V_t}
\left(
\rho\,dB_t^{1}
+
\sqrt{1-\rho^2}\,dB_t^{2}
\right),
\\
V_t
&=
\nu_0(t)
\exp\left(
\eta\sqrt{2H}
\int_0^t
(t-s)^{H-\frac12}\,dB_s^{1}
-
\frac{\eta^2}{2}t^{2H}
\right),
\end{flalign*}
where $H \in (0,\frac12)$, $\eta > 0$, $\rho \in (-1,1)$, and $\nu_0 \colon [0,T_V] \to (0,\infty)$. We also denote by
\begin{gather*}
\nu_t(u)
\coloneqq
\mathbb E_{\mathbb Q}
\left[
V_u
\,\middle|\,
\mathcal F_t
\right],
\qquad
0\leq t\leq u\leq T_V,
\end{gather*}
the forward variance curve. The rough Bergomi model is non-Markovian and does not admit the finite-dimensional forward representation required by the Markovian methods discussed above; see \textcite{Bayer02062016}.

Following \textcite{FukasawaHorvathTankov2023Hedging}, we assume that, in addition to the underlying asset, a variance swap with fixed maturity $T_V>T$ is continuously tradable. We denote its price process by
\begin{gather*}
U_t
=
e^{-r(T_V-t)}
\mathbb E_{\mathbb Q}
\left[
\int_0^{T_V} V_u\,du
\,\middle|\,
\mathcal F_t
\right],
\qquad
0\leq t\leq T,
\end{gather*}
whose dynamics are
\begin{gather*}
dU_t
=
rU_t\,dt+\gamma_t^U\,dB_t^1,
\qquad
\gamma_t^U
=
e^{-r(T_V-t)}
\eta\sqrt{2H}
\int_t^{T_V}
(u-t)^{H-\frac12}\nu_t(u)\,du.
\end{gather*}
Since $\gamma_t^U>0$ for $t\leq T<T_V$, the volatility matrix of $(S,U)$ is invertible, and the resulting market is complete.

Let $\xi$ be any contingent claim written on $S$. The pricing and hedging of $\xi$ can be related to the solution to the BSDE
\begin{gather*}
Y_t
=
\xi
-
\int_t^T rY_s\,ds
-
\int_t^T Z_s^{1}\,dB_s^{1}
-
\int_t^T Z_s^{2}\,dB_s^{2},
\qquad
0\leq t\leq T.
\end{gather*}
Indeed, the price of $\xi$ at time $t$ is given by $Y_t$, and the strategy $(\Delta_t,\Lambda_t)$ in $(S,U)$ is
\begin{gather*}
\Delta_t
=
\frac{Z_t^{2}}
{S_t\sqrt{V_t}\sqrt{1-\rho^2}},
\qquad
\Lambda_t
=
(\gamma_t^U)^{-1}
\left(
Z_t^{1}
-
\frac{\rho}{\sqrt{1-\rho^2}}Z_t^{2}
\right).
\end{gather*}

In the numerical experiment, we take $\xi=(S_T-K)^+$, with $T=1$ and $K=1.1$. We also set $r=0.01$, $S_0=1$, $H=0.25$, $\eta=1.5$, $\rho=-0.8$, and take the flat initial forward variance curve $\nu_0(t)=0.04$ for $t\in[0,T_V]$.

We obtain benchmark values for the initial price and stock hedge from the Monte Carlo representations
\begin{gather*}
Y_0
=
e^{-rT}
\mathbb{E}_{\mathbb{Q}}
\left[
(S_T-K)^+
\right],
\qquad
\Delta_0
=
\frac{e^{-rT}}{S_0}
\mathbb{E}_{\mathbb{Q}}
\left[
S_T\mathbf{1}_{\{S_T>K\}}
\right],
\end{gather*}
where the second formula follows from \textcite{BroadieGlasserman1996Estimating}. The stock process is simulated using the same $50$-step discretization and hybrid scheme as in the Euler--chaos experiment. We use $10^8$ independent Monte Carlo samples.

For the Euler scheme via the Wiener chaos decomposition, we take $P=3$, $M=15$, and a uniform partition with $m=50$ time steps. The chaos coefficients are estimated using $N=10^6$ Monte Carlo samples. We perform $100$ independent runs of the algorithm. Table~\ref{tab:rough-bergomi-euler} reports the empirical mean and standard deviation of the resulting approximations and absolute relative errors.

\begin{table}[t]
\centering
\begin{tabular}{lccc}
\toprule
Quantity
& Monte Carlo reference
& Euler--chaos
& Relative error \\
\midrule
$Y_0$
& $0.026143 \pm 4.5\times 10^{-6}$
& $0.026144 \pm 5.1\times 10^{-5}$
& $(0.152 \pm 0.122)\%$ \\
$\Delta_0$
& $0.39406 \pm 4.9\times 10^{-5}$
& $0.3673 \pm 2.0\times 10^{-3}$
& $(6.80 \pm 0.51)\%$ \\
\bottomrule
\end{tabular}
\caption{Euler--chaos approximation for the rough Bergomi example over $100$ independent runs. For the Monte Carlo reference, the reported uncertainty is the standard error. The Euler--chaos approximations and their absolute relative errors are reported as empirical mean $\pm$ standard deviation across the independent runs.}
\label{tab:rough-bergomi-euler}
\end{table}

The initial price is accurately reproduced, with a small dispersion across the independent runs. The initial stock hedge has an average relative error of approximately $6.8\%$, with a comparatively small run-to-run standard deviation.

\subsection{Example 4 -- Non-linear pricing and hedging in an SVV model}

We conclude with an example based on the sandwiched Volterra volatility (SVV) model introduced in \textcite{DiNunnoMishuraYurchenko2024SVV}; see also \textcite{DiNunnoMishuraYurchenkoTytarenko2026Volterra}. Let $B=(B^1,B^2)$ be a two-dimensional Brownian motion and consider
\begin{align*}
\frac{dS_t}{S_t}
&=
r\,dt
+
\sigma_t
\left(
\rho\,dB_t^1
+
\sqrt{1-\rho^2}\,dB_t^2
\right),
\\
\sigma_t
&=
\sigma_0
+
\int_0^t b(\sigma_s)\,ds
+
\frac{0.8}{\Gamma(H+\frac12)}
\int_0^t
(t-s)^{H-\frac12}\,dB_s^1,
\end{align*}
where
\[
b(\sigma)
=
\frac{\kappa_1}{(\sigma-\underline{\sigma})^\gamma}
-
\frac{\kappa_2}{(\overline{\sigma}-\sigma)^\gamma}
+
\kappa_0-\kappa_3\sigma.
\]
We take $H=0.2$, $\sigma_0=1.2$, $\underline{\sigma}=0.05$, $\overline{\sigma}=2.5$, $\gamma=5$, $\kappa_0=0.0125$, $\kappa_1=\kappa_2=0.005$, $\kappa_3=0.025$, and $\rho=-0.7$. The singular drift keeps the volatility between the prescribed boundaries, while the Volterra noise makes the model genuinely path dependent and prevents a natural low-dimensional Markovian representation; see \textcite{DiNunnoYurchenko2026SVV}.

To complete the market, we introduce a stylized volatility-factor security
\[
\frac{dU_t}{U_t}
=
r\,dt+\eta_U\,dB_t^1,
\qquad
U_0=1,
\]
with $\eta_U=0.3$. Since $S_t$, $U_t$, $\sigma_t$, and $\eta_U$ are positive and $|\rho|<1$, the volatility matrix of $(S,U)$ is invertible.

We consider the European straddle $\xi=|S_T-K|$, with $T=1$ and $S_0=K=1$. Here, $r$ denotes the reference rate in the risk-neutral dynamics. Positive cash balances earn the lending rate $\underline{r}$, whereas negative balances are financed at the borrowing rate $\overline{r}$, where $\underline{r}\leq r\leq\overline{r}$. For $z=(z^1,z^2)$, matching the Brownian exposures gives
\[
\Delta_t S_t
=
\frac{z^2}{\sigma_t\sqrt{1-\rho^2}},
\qquad
\Lambda_t U_t
=
\frac{1}{\eta_U}
\left(
z^1
-
\frac{\rho}{\sqrt{1-\rho^2}}z^2
\right).
\]
Hence, setting
\[
\Psi_t(z)
=
\frac{z^2}{\sigma_t\sqrt{1-\rho^2}}
+
\frac{1}{\eta_U}
\left(
z^1
-
\frac{\rho}{\sqrt{1-\rho^2}}z^2
\right),
\]
the associated BSDE driver is
\[
f(t,y,z;\sigma_t)
=
-\underline{r}y
+
(\underline{r}-r)\Psi_t(z)
+
(\overline{r}-\underline{r})
\bigl(\Psi_t(z)-y\bigr)^+.
\]
For the choice $r=\underline{r}$ used below, this simplifies to
\[
f(t,y,z;\sigma_t)
=
-\underline{r}y
+
(\overline{r}-\underline{r})
\bigl(\Psi_t(z)-y\bigr)^+.
\]

We set $r=\underline{r}=0.01$ and consider $\overline{r}\in\{0.01,0.03,0.06,0.10\}$. The forward model and the BSDE Euler scheme are simulated on uniform grids with $100$ time steps. The volatility process is approximated using the drift-implicit scheme of \textcite{DiNunnoMishuraYurchenko2023Implicit}. For the Euler--chaos method, we take $P=4$, $M=15$, and estimate the chaos coefficients using $10^6$ Monte Carlo samples.

When $\overline{r}=\underline{r}$, the BSDE is linear and independent Monte Carlo benchmarks are available:
\[
Y_0
=
e^{-\underline{r}T}
\mathbb{E}
\left[
|S_T-K|
\right],
\qquad
\Delta_0
=
\frac{e^{-\underline{r}T}}{S_0}
\mathbb{E}
\left[
S_T\operatorname{sgn}(S_T-K)
\right].
\]
The stock process is simulated using the same $100$-step discretization and drift-implicit scheme as in the Euler--chaos experiment. Using $10^6$ independent Monte Carlo samples gives the comparison in Table~\ref{tab:svv-linear-benchmark}.

\begin{table}[ht]
\centering
\begin{tabular}{lcc}
\toprule
Quantity
& Monte Carlo reference
& Euler--chaos \\
\midrule
$Y_0$
& $0.9598 \pm 0.0014$
& $0.9599$ \\
$\Delta_0$
& $0.5329 \pm 0.0018$
& $0.5056$ \\
\bottomrule
\end{tabular}
\caption{Monte Carlo benchmark and Euler--chaos approximation for $r=\underline{r}=\overline{r}=0.01$. The uncertainty reported for the Monte Carlo reference is the standard error.}
\label{tab:svv-linear-benchmark}
\end{table}

The initial-price estimates agree within the Monte Carlo uncertainty, while the stock hedge differs by approximately $5.1\%$, similarly to the rough Bergomi experiment.

\begin{figure}[ht]
\centering
\begin{subfigure}{0.48\textwidth}
\includegraphics[width=\textwidth]{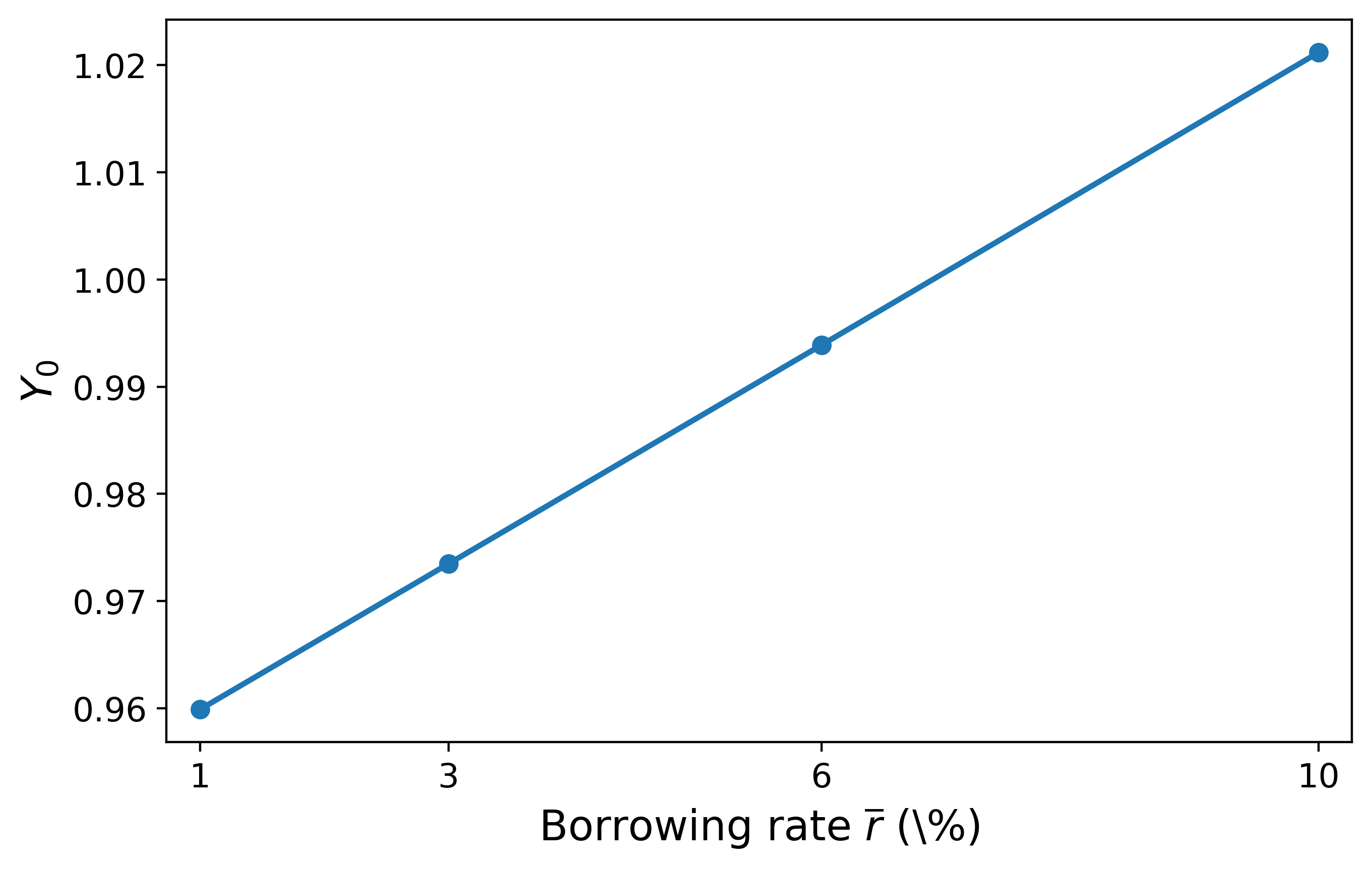}
\caption{$Y_0$.}
\end{subfigure}
\hfill
\begin{subfigure}{0.48\textwidth}
\includegraphics[width=\textwidth]{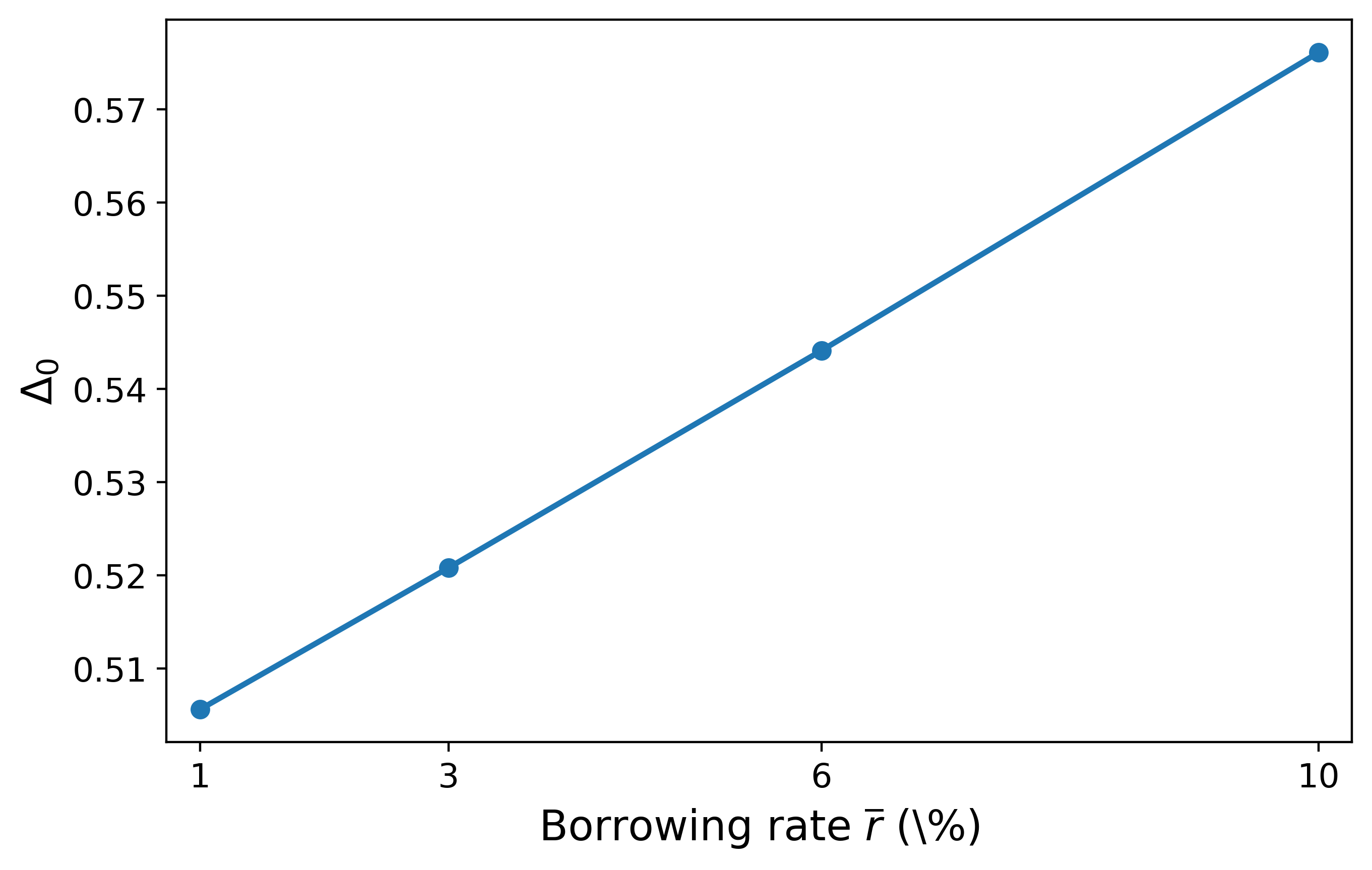}
\caption{$\Delta_0$.}
\end{subfigure}
\caption{Euler--chaos initial value and stock hedge as functions of the borrowing rate $\overline{r}$, with $r=\underline{r}=0.01$.}
\label{fig:svv-funding-sensitivity}
\end{figure}

We next keep $r=\underline{r}=0.01$ fixed and increase $\overline{r}$. We use common random numbers for all borrowing rates. Figure~\ref{fig:svv-funding-sensitivity} shows that both the initial price and the stock hedge increase smoothly with the borrowing rate.

\begin{figure}[ht]
\centering
\begin{subfigure}{0.48\textwidth}
\includegraphics[width=\textwidth]{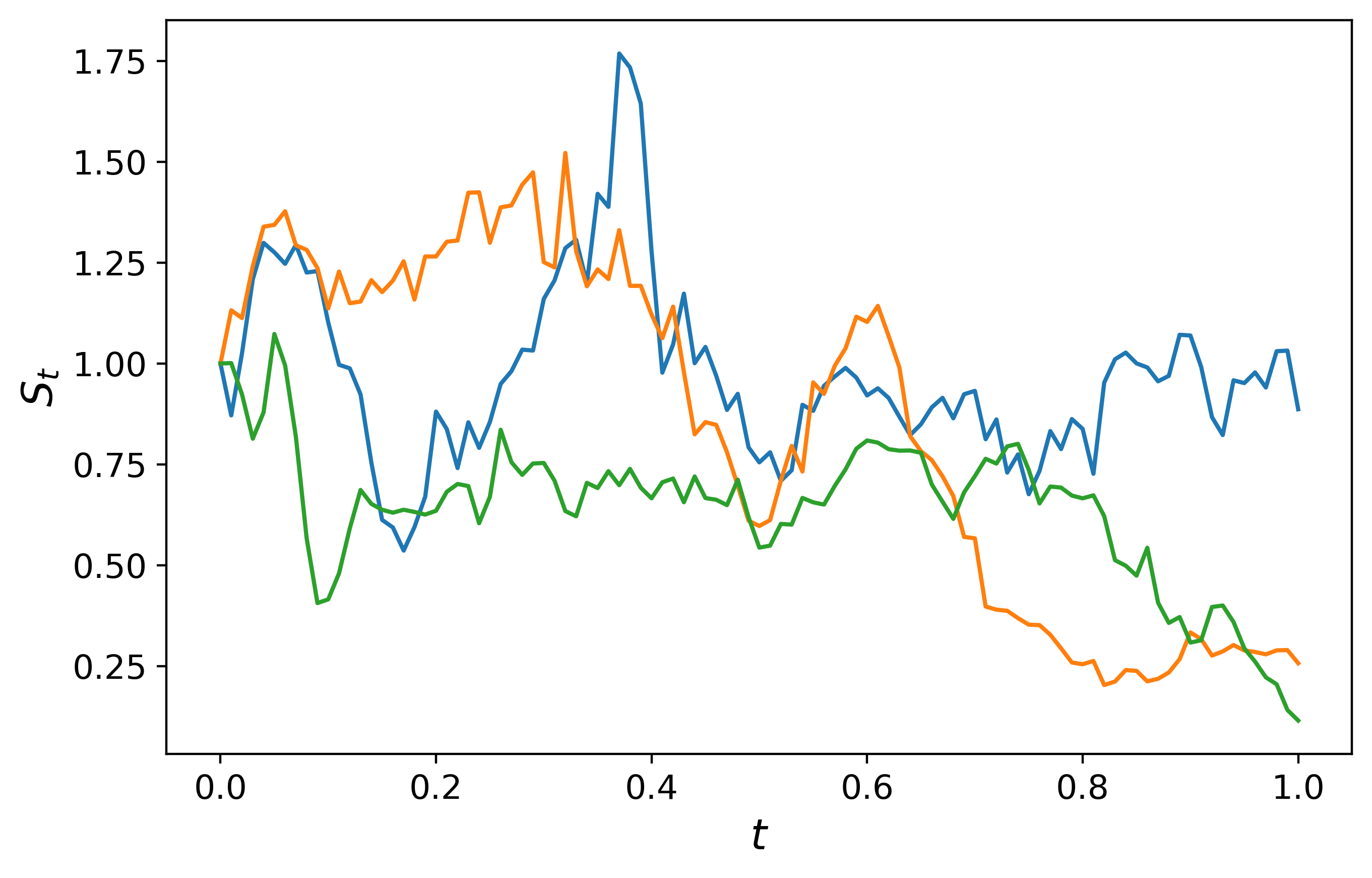}
\caption{$S_t$.}
\end{subfigure}
\hfill
\begin{subfigure}{0.48\textwidth}
\includegraphics[width=\textwidth]{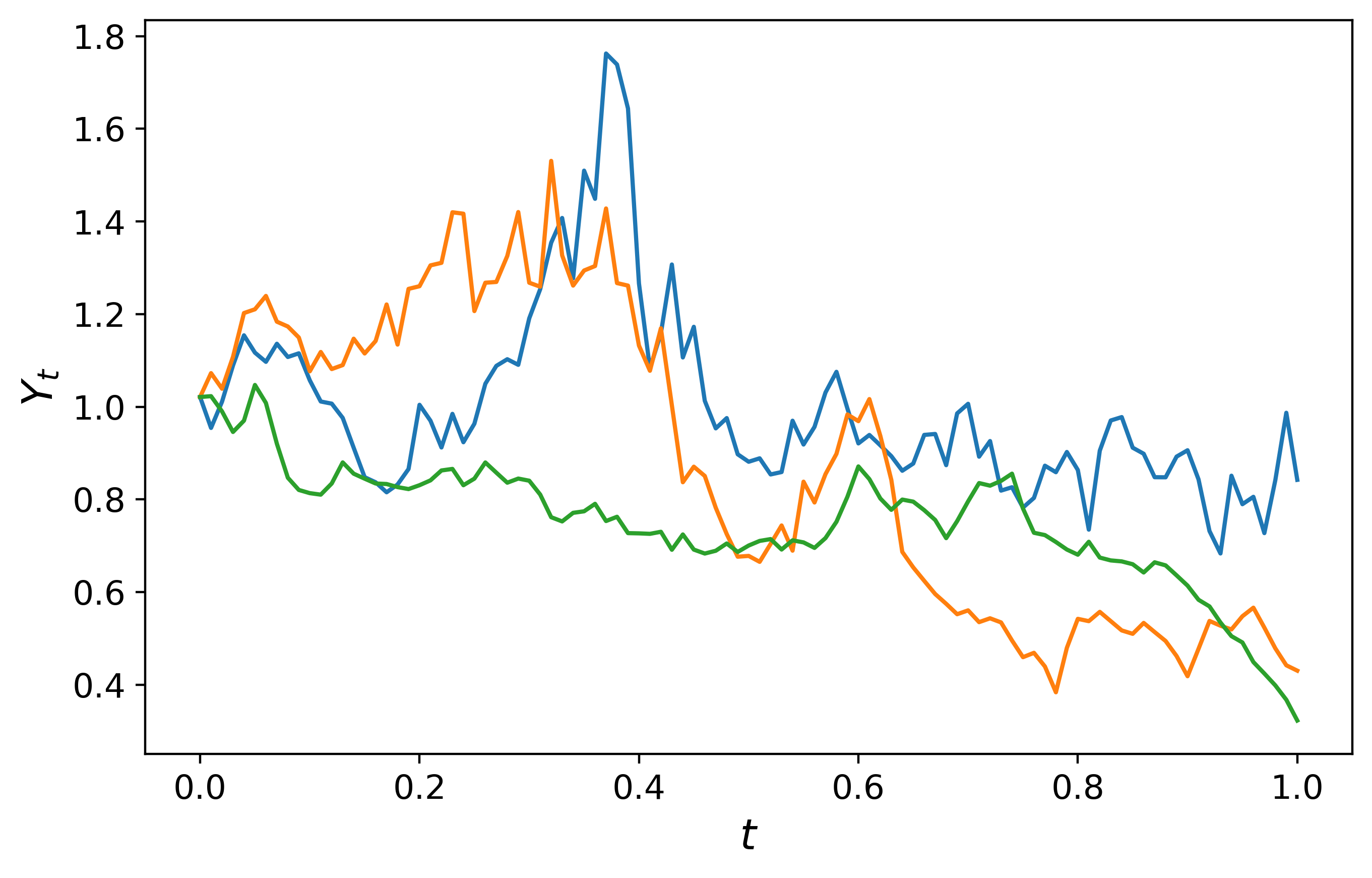}
\caption{$Y_t$.}
\end{subfigure}
\caption{Sample paths of the stock and Euler--chaos value processes for $\overline{r}=0.10$ and $r=\underline{r}=0.01$.}
\label{fig:svv-sample-paths}
\end{figure}

Finally, Figure~\ref{fig:svv-sample-paths} shows sample paths of the stock and value processes for $\overline{r}=0.10$.

\section*{Acknowledgments}

The authors thank Lukas Gonon for suggesting the inclusion of random generators, and more generally the participants of the Workshop on Mathematical Finance and Related Issues, held in Osaka in 2026 and organized by Masaaki Fukasawa, for interesting discussions and comments that helped improve the paper.

Giulia Di Nunno is also affiliated with the Department of Business and Management Science, NHH Norwegian School of Economics, Bergen, Norway. Pere Diaz-Lozano is the corresponding author.

This work was carried out within the SURE-AI Center, funded by the Research Council of Norway, grant no.~357482.

\appendix

\section{Multi-dimensional Wiener chaos formulae}\label{section: multidimensional formulas}

We now present the framework for a $d$-dimensional Brownian motion. Let $(P,M)\in\mathbb N^2$ and let $\overline{\pi}_i$, $i=1,\dots,m$, be constructed as in Section~\ref{sec: description basis}. Let $(e_l)_{1\leq l\leq d}$ be the canonical basis of $\mathbb R^d$.

For each $i\in\{1,\dots,m\}$, define
\begin{gather}\label{truncated basis appendix}
    h_{j,l}^i(t)\coloneqq e_l\,\mathbf 1_{(s_{j-1}^i,s_j^i]}(t)/\sqrt{\delta_j^i}, \quad \delta_j^i\coloneqq s_j^i-s_{j-1}^i, \quad j=1,\dots,M(i),\quad l=1,\dots,d.
\end{gather}
The family $(h_{j,l}^i)_{1\leq j\leq M(i),\,1\leq l\leq d}$ is an orthonormal system in $\mathbb H$ and can be completed to an orthonormal basis of $\mathbb H$. Moreover,
\begin{gather*}
    B(h_{j,l}^i)\coloneqq\int_0^T h_{j,l}^i(s)\cdot dB_s=\frac{B_{s_j^i}^l-B_{s_{j-1}^i}^l}{\sqrt{\delta_j^i}}, \quad j=1,\dots,M(i),\quad l=1,\dots,d.
\end{gather*}

For $a=(a^1,\dots,a^d)\in\mathcal A_{M(i)}^d$, write $a^l=(a_j^l)_{j\geq1}$, $|a|\coloneqq|a^1|+\cdots+|a^d|$, and $a!\coloneqq\prod_{l=1}^d\prod_{j=1}^{M(i)}a_j^l!$. Throughout this section, all sums indexed by $a$ are understood to run over $a\in\mathcal A_{M(i)}^d$.

We denote by $\mathcal C_P^{\overline{\pi}_i}(\cdot)$ the truncated chaos decomposition with basis given by \eqref{truncated basis appendix}. Thus, for $F\in L^2(\mathcal F_{t_i})$,
\begin{flalign*}
    \mathcal C_P^{\overline{\pi}_i}(F)
    &\coloneqq d_0+\sum_{k=1}^P\sum_{|a|=k}d_a\prod_{l=1}^d\prod_{j=1}^{M(i)}H_{a_j^l}\big(B(h_{j,l}^i)\big),\\
    d_0&\coloneqq\mathbb E(F),\qquad d_a\coloneqq a!\,\mathbb E\Big(F\times\prod_{l=1}^d\prod_{j=1}^{M(i)}H_{a_j^l}\big(B(h_{j,l}^i)\big)\Big).
\end{flalign*}

We now state the corresponding multi-dimensional formulae. For each $i\in\{1,\dots,m\}$, let $d_a^i$ denote the chaos coefficient of $\mathcal C_P^{\overline{\pi}_i}(F_i^{\pi,\theta})$ associated with $a\in\mathcal A_{M(i)}^d$. Related formulae appear in \textcite[Proposition B.1]{BriandPhilippe2014SOBB}.

\begin{proposition}
Let $i\in\{1,\dots,m\}$ and $r\in\{1,\dots,M(i)\}$ be such that $t\in[t_{i-1},t_i)\cap(s_{r-1}^i,s_r^i]$. Then
\begin{flalign*}
Y_t^{\pi,\theta}
&=
d_0^i+\sum_{k=1}^P\sum_{\substack{a\in\mathcal A_{M(i)}^{d,r}\\|a|=k}}d_a^i\prod_{j<r}\prod_{l=1}^dH_{a_j^l}\big(B(h_{j,l}^i)\big)\times\prod_{l=1}^d\Bigg(\frac{t-s_{r-1}^i}{\delta_r^i}\Bigg)^{a_r^l/2}H_{a_r^l}\Bigg(\frac{B_t^l-B_{s_{r-1}^i}^l}{\sqrt{t-s_{r-1}^i}}\Bigg),\\
(Z_t^{\pi,\theta})^\gamma
&=
(\delta_r^i)^{-1/2}\sum_{k=1}^P\sum_{\substack{a\in\mathcal A_{M(i)}^{d,r}\\|a|=k\\a_r^\gamma>0}}d_a^i\prod_{j<r}\prod_{l=1}^dH_{a_j^l}\big(B(h_{j,l}^i)\big)\\
&\quad\times\Bigg(\frac{t-s_{r-1}^i}{\delta_r^i}\Bigg)^{(a_r^1+\cdots+a_r^d-1)/2}H_{a_r^\gamma-1}\Bigg(\frac{B_t^\gamma-B_{s_{r-1}^i}^\gamma}{\sqrt{t-s_{r-1}^i}}\Bigg)\times\prod_{l\neq\gamma}H_{a_r^l}\Bigg(\frac{B_t^l-B_{s_{r-1}^i}^l}{\sqrt{t-s_{r-1}^i}}\Bigg),
\end{flalign*}
for $\gamma=1,\dots,d$, where
\[
\mathcal A_{M(i)}^{d,r}
\coloneqq
\{a\in\mathcal A_{M(i)}^d:a_j^l=0\text{ for all }j>r,\ l=1,\dots,d\}.
\]
At $t=0$, we use the continuous extensions $Y_0^{\pi,\theta}=d_0^1$ and $(Z_0^{\pi,\theta})^\gamma=(\delta_1^1)^{-1/2}d_{e_{1,\gamma}}^1$, where $e_{1,\gamma}$ denotes the multi-index whose only non-zero component is $a_1^\gamma=1$.
\end{proposition}

\begin{proposition}\label{prop:formula-for-scalarZ}
Fix $\gamma\in\{1,\dots,d\}$, let $i\in\{2,\dots,m\}$, and let $u$ be such that $t_{i-1}\in(s_{u-1}^i,s_u^i]$. Let
\begin{gather*}
    C_1(i)\coloneqq\frac{s_u^i-t_{i-1}}{\sqrt{\delta_u^i}}, \quad C_2(i)\coloneqq\frac{t_{i-1}-s_{u-1}^i}{\delta_u^i}.
\end{gather*}
Then $(\sprod{Z}_{i-1}^{\pi,\theta})^\gamma=$
\begin{flalign*}
\frac{1}{\Delta_i}\sum_{k=1}^P
&\Bigg\{
C_1(i)
\sum_{\substack{a\in\mathcal A_{M(i)}^{d,u}\\|a|=k\\a_u^\gamma>0}}
d_a^i
\prod_{j<u}\prod_{l=1}^dH_{a_j^l}\big(B(h_{j,l}^i)\big)
\prod_{l\neq\gamma}H_{a_u^l}\Bigg(
\frac{B_{t_{i-1}}^l-B_{s_{u-1}^i}^l}
{\sqrt{t_{i-1}-s_{u-1}^i}}
\Bigg)
\times\big[C_2(i)\big]^{\frac{a_u^1+\cdots+a_u^d-1}{2}}\\
&\quad\times
H_{a_u^\gamma-1}\Bigg(
\frac{B_{t_{i-1}}^\gamma-B_{s_{u-1}^i}^\gamma}
{\sqrt{t_{i-1}-s_{u-1}^i}}
\Bigg)
+
\sum_{r=u+1}^{M(i)}
\sqrt{\delta_r^i}
\sum_{\substack{a\in\mathcal A_{M(i)}^{d,u,r,\gamma}\\|a|=k}}
d_a^i
\prod_{j<u}\prod_{l=1}^dH_{a_j^l}\big(B(h_{j,l}^i)\big)\\
&\quad\times
\big[C_2(i)\big]^{\frac{a_u^1+\cdots+a_u^d}{2}}
\prod_{l=1}^dH_{a_u^l}\Bigg(
\frac{B_{t_{i-1}}^l-B_{s_{u-1}^i}^l}
{\sqrt{t_{i-1}-s_{u-1}^i}}
\Bigg)
\Bigg\},
\end{flalign*}
where $\mathcal A_{M(i)}^{d,u}\coloneqq\{a\in\mathcal A_{M(i)}^d:a_j^l=0\text{ for all }j>u,\ l=1,\dots,d\}$ and
$\mathcal A_{M(i)}^{d,u,r,\gamma}\coloneqq\{a\in\mathcal A_{M(i)}^d:a_r^\gamma=1,\ a_j^l=0\text{ whenever }j>u\text{ and }(j,l)\neq(r,\gamma)\}$.
\end{proposition}

\begin{proof}
We start by decomposing the integrals and using conditional Fubini,
\begin{flalign*}
\mathbb E_{t_{i-1}}\Big(\int_{t_{i-1}}^{t_i}Z_t^{\pi,\theta}\,dt\Big)
&=
\int_{t_{i-1}}^{t_i}\mathbb E_{t_{i-1}}\big(Z_t^{\pi,\theta}\big)\,dt\\
&=
\int_{t_{i-1}}^{s_u^i}\mathbb E_{t_{i-1}}\big(Z_t^{\pi,\theta}\big)\,dt
+
\sum_{r=u+1}^{M(i)}
\int_{s_{r-1}^i}^{s_r^i}
\mathbb E_{t_{i-1}}\big(Z_t^{\pi,\theta}\big)\,dt.
\end{flalign*}
We first focus on the terms inside the summation. For
$t\in[t_{i-1},t_i)\cap(s_{r-1}^i,s_r^i]$,
\begin{align*}
\mathbb E_{t_{i-1}}\big((Z_t^{\pi,\theta})^\gamma\big)
&=
(\delta_r^i)^{-1/2}
\sum_{k=1}^P
\sum_{\substack{a\in\mathcal A_{M(i)}^{d,r}\\|a|=k\\a_r^\gamma>0}}
d_a^i
\prod_{j<u}\prod_{l=1}^dH_{a_j^l}\big(B(h_{j,l}^i)\big)
\prod_{l=1}^d
\mathbb E_{t_{i-1}}\big(H_{a_u^l}(B(h_{u,l}^i))\big)\\
&\quad\times
\prod_{u<j<r}\prod_{l=1}^d
\mathbb E\big(H_{a_j^l}(B(h_{j,l}^i))\big)
\Bigg(\frac{t-s_{r-1}^i}{\delta_r^i}\Bigg)^{\frac{a_r^\gamma-1}{2}}
\mathbb E\Bigg[
H_{a_r^\gamma-1}\Bigg(
\frac{B_t^\gamma-B_{s_{r-1}^i}^\gamma}
{\sqrt{t-s_{r-1}^i}}
\Bigg)
\Bigg]\\
&\quad\times
\prod_{l\neq\gamma}
\Bigg(\frac{t-s_{r-1}^i}{\delta_r^i}\Bigg)^{a_r^l/2}
\mathbb E\Bigg[
H_{a_r^l}\Bigg(
\frac{B_t^l-B_{s_{r-1}^i}^l}
{\sqrt{t-s_{r-1}^i}}
\Bigg)
\Bigg].
\end{align*}
Hence, in order for the previous expression to be different from zero, we need:
\begin{itemize}
    \item $a_j^l=0$ for $u<j<r$ and $l=1,\dots,d$.
    \item $a_r^\gamma=1$.
    \item $a_r^l=0$ for $l\neq\gamma$.
\end{itemize}
Since
\begin{gather*}
\mathbb E_{t_{i-1}}\big(H_{a_u^l}(B(h_{u,l}^i))\big)
=
\Bigg(\frac{t_{i-1}-s_{u-1}^i}{\delta_u^i}\Bigg)^{a_u^l/2}
H_{a_u^l}\Bigg(
\frac{B_{t_{i-1}}^l-B_{s_{u-1}^i}^l}
{\sqrt{t_{i-1}-s_{u-1}^i}}
\Bigg),
\end{gather*}
we get
\begin{flalign*}
\mathbb E_{t_{i-1}}\big((Z_t^{\pi,\theta})^\gamma\big)
=
(\delta_r^i)^{-1/2}
\sum_{k=1}^P
&\sum_{\substack{a\in\mathcal A_{M(i)}^{d,u,r,\gamma}\\|a|=k}}
d_a^i
\prod_{j<u}\prod_{l=1}^dH_{a_j^l}\big(B(h_{j,l}^i)\big)\\
&\times
\prod_{l=1}^d
\Bigg(\frac{t_{i-1}-s_{u-1}^i}{\delta_u^i}\Bigg)^{a_u^l/2}
H_{a_u^l}\Bigg(
\frac{B_{t_{i-1}}^l-B_{s_{u-1}^i}^l}
{\sqrt{t_{i-1}-s_{u-1}^i}}
\Bigg).
\end{flalign*}
Notice that this does not depend on $t$. Therefore, for $u<r$,
\begin{flalign*}
\int_{s_{r-1}^i}^{s_r^i}
\mathbb E_{t_{i-1}}\big((Z_t^{\pi,\theta})^\gamma\big)\,dt
=
\sqrt{\delta_r^i}
\sum_{k=1}^P
&\sum_{\substack{a\in\mathcal A_{M(i)}^{d,u,r,\gamma}\\|a|=k}}
d_a^i
\prod_{j<u}\prod_{l=1}^dH_{a_j^l}\big(B(h_{j,l}^i)\big)\\
&\times
\prod_{l=1}^d
\Bigg(\frac{t_{i-1}-s_{u-1}^i}{\delta_u^i}\Bigg)^{a_u^l/2}
H_{a_u^l}\Bigg(
\frac{B_{t_{i-1}}^l-B_{s_{u-1}^i}^l}
{\sqrt{t_{i-1}-s_{u-1}^i}}
\Bigg).
\end{flalign*}

Now, for $t\in[t_{i-1},t_i)\cap(s_{u-1}^i,s_u^i]$, we have that
$\mathbb E_{t_{i-1}}\big((Z_t^{\pi,\theta})^\gamma\big)=$
\begin{flalign*}
(\delta_u^i)^{-1/2}
\sum_{k=1}^P
&\sum_{\substack{a\in\mathcal A_{M(i)}^{d,u}\\|a|=k\\a_u^\gamma>0}}
d_a^i
\prod_{j<u}\prod_{l=1}^dH_{a_j^l}\big(B(h_{j,l}^i)\big)
\Bigg(\frac{t-s_{u-1}^i}{\delta_u^i}\Bigg)^{\frac{a_u^\gamma-1}{2}}\\
&\times
\mathbb E_{t_{i-1}}\Bigg[
H_{a_u^\gamma-1}\Bigg(
\frac{B_t^\gamma-B_{s_{u-1}^i}^\gamma}
{\sqrt{t-s_{u-1}^i}}
\Bigg)
\Bigg]
\prod_{l\neq\gamma}
\Bigg(\frac{t-s_{u-1}^i}{\delta_u^i}\Bigg)^{a_u^l/2}\\
&\times
\mathbb E_{t_{i-1}}\Bigg[
H_{a_u^l}\Bigg(
\frac{B_t^l-B_{s_{u-1}^i}^l}
{\sqrt{t-s_{u-1}^i}}
\Bigg)
\Bigg].
\end{flalign*}
Computing the conditional expectations as before, the previous expression equals
\begin{flalign*}
(\delta_u^i)^{-1/2}
\sum_{k=1}^P
&\sum_{\substack{a\in\mathcal A_{M(i)}^{d,u}\\|a|=k\\a_u^\gamma>0}}
d_a^i
\prod_{j<u}\prod_{l=1}^dH_{a_j^l}\big(B(h_{j,l}^i)\big)
\Bigg(\frac{t-s_{u-1}^i}{\delta_u^i}\Bigg)^{\frac{a_u^\gamma-1}{2}}
\Bigg(\frac{t_{i-1}-s_{u-1}^i}{t-s_{u-1}^i}\Bigg)^{\frac{a_u^\gamma-1}{2}}\\
&\times
H_{a_u^\gamma-1}\Bigg(
\frac{B_{t_{i-1}}^\gamma-B_{s_{u-1}^i}^\gamma}
{\sqrt{t_{i-1}-s_{u-1}^i}}
\Bigg)
\prod_{l\neq\gamma}
\Bigg(\frac{t-s_{u-1}^i}{\delta_u^i}\Bigg)^{a_u^l/2}
\Bigg(\frac{t_{i-1}-s_{u-1}^i}{t-s_{u-1}^i}\Bigg)^{a_u^l/2}\\
&\times
H_{a_u^l}\Bigg(
\frac{B_{t_{i-1}}^l-B_{s_{u-1}^i}^l}
{\sqrt{t_{i-1}-s_{u-1}^i}}
\Bigg).
\end{flalign*}
Cancelling terms and integrating, we get that
$\int_{t_{i-1}}^{s_u^i}\mathbb E_{t_{i-1}}\big((Z_t^{\pi,\theta})^\gamma\big)\,dt=$
\begin{flalign*}
& (s_u^i-t_{i-1})(\delta_u^i)^{-1/2}
\sum_{k=1}^P
\sum_{\substack{a\in\mathcal A_{M(i)}^{d,u}\\|a|=k\\a_u^\gamma>0}}
d_a^i
\prod_{j<u}\prod_{l=1}^dH_{a_j^l}\big(B(h_{j,l}^i)\big)
\Bigg(\frac{t_{i-1}-s_{u-1}^i}{\delta_u^i}\Bigg)^{\frac{a_u^\gamma-1}{2}}\\
& \qquad \times
H_{a_u^\gamma-1}\Bigg(
\frac{B_{t_{i-1}}^\gamma-B_{s_{u-1}^i}^\gamma}
{\sqrt{t_{i-1}-s_{u-1}^i}}
\Bigg)
\prod_{l\neq\gamma}
\Bigg(\frac{t_{i-1}-s_{u-1}^i}{\delta_u^i}\Bigg)^{a_u^l/2}
H_{a_u^l}\Bigg(
\frac{B_{t_{i-1}}^l-B_{s_{u-1}^i}^l}
{\sqrt{t_{i-1}-s_{u-1}^i}}
\Bigg).
\end{flalign*}
Combining this identity with the contribution from the intervals $(s_{r-1}^i,s_r^i]$, $r=u+1,\dots,M(i)$, and dividing by $\Delta_i$ gives the desired formula.
\end{proof}

\section{Proofs}\label{Proofs section convergence analysis}

\subsection{Proofs of Subsection \ref{subsection: Error given by the truncation of the chaos}}

\subsubsection*{Proof of Lemma \ref{lemma error implementation}}

For $i\in\{1,\dots,m\}$, set
\[
\Delta^\theta Y_{t_i}^{\pi}
:=
Y_{t_i}^{\pi}-Y_{t_i}^{\pi,\theta},
\qquad
\Delta^\theta Z^{\pi}
:=
Z^{\pi}-Z^{\pi,\theta},
\qquad
\Delta^\theta\sprod{Z}_i^\pi
:=
\sprod{Z}_i^\pi-\sprod{Z}_i^{\pi,\theta},
\]
and
\[
R_i
:=
\mathbb{E}\left(
    \abs{f_i^{\pi,\theta}-f_i^{M,\pi,\theta}}^2
    +
    \abs{
        f_i^{M,\pi,\theta}
        -
        \mathcal{C}_{P}(f_i^{M,\pi,\theta})
    }^2
\right).
\]
Since
$f_i^{M,\pi,\theta}
=
\mathbb{E}(f_i^{\pi,\theta}\mid\mathcal{G}_i^M)$
and
$\mathcal{C}_{P}(f_i^{M,\pi,\theta})$
is $\mathcal{G}_i^M$-measurable, orthogonality of conditional expectation yields
\[
\mathbb{E}\left(
    \abs{
        f_i^{\pi,\theta}
        -
        \mathcal{C}_{P}(f_i^{M,\pi,\theta})
    }^2
\right)
=
R_i.
\]

Define
\[
I_m:=\mathbb{E}\big(\abs{\Delta^\theta Y_{t_m}^{\pi}}^2\big),
\qquad
I_i
:=
\mathbb{E}\left(
    \abs{\Delta^\theta Y_{t_i}^{\pi}}^2
    +
    \int_{t_i}^{t_{i+1}}
    \abs{\Delta^\theta Z_u^\pi}^2\,du
\right),
\quad i=0,\dots,m-1.
\]
By the definition,
\[
\Delta^\theta Y_{t_{i-1}}^\pi
=
F_i^\pi-F_i^{\pi,\theta}
+
\Delta_i\left[
    f_i^{\pi,\theta}
    -
    \mathcal{C}_{P}(f_i^{M,\pi,\theta})
\right]
-
\int_{t_{i-1}}^{t_i}
\Delta^\theta Z_u^\pi\,dB_u.
\]
Using the orthogonality of the It\^o integral and Young's inequality, for
$|\pi|$ sufficiently small,
\begin{equation}\label{eq:Ii-1-start}
I_{i-1}
\leq
(1+\Delta_i)
\mathbb{E}\big(
    \abs{F_i^\pi-F_i^{\pi,\theta}}^2
\big)
+
2\Delta_iR_i.
\end{equation}

The Lipschitz continuity of $f$, together with Young's inequality, gives
\begin{equation}\label{eq:Fi-diff-simplified}
\mathbb{E}\big(
    \abs{F_i^\pi-F_i^{\pi,\theta}}^2
\big)
\leq
(1+C\Delta_i)
\mathbb{E}\big(
    \abs{\Delta^\theta Y_{t_i}^\pi}^2
\big)
+
\frac{\Delta_i}{2L}
\mathbb{E}\big(
    \abs{\Delta^\theta\sprod{Z}_i^\pi}^2
\big),
\end{equation}
for a constant $C$ depending only on $L$ and $[f]_L$. Moreover, for
$i=1,\dots,m-1$,
\[
\mathbb{E}\big(
    \abs{\Delta^\theta\sprod{Z}_i^\pi}^2
\big)
\leq
\frac{1}{\Delta_{i+1}}
\mathbb{E}\left(
    \int_{t_i}^{t_{i+1}}
    \abs{\Delta^\theta Z_u^\pi}^2\,du
\right).
\]
Hence, using $\Delta_i/\Delta_{i+1}\leq L$,
\begin{equation}\label{eq:main-step}
I_{i-1}
+
\frac12
\mathbb{E}\left(
    \int_{t_i}^{t_{i+1}}
    \abs{\Delta^\theta Z_u^\pi}^2\,du
\right)
\leq
(1+C\Delta_i)I_i
+
2\Delta_iR_i,
\end{equation}
for $i=1,\dots,m-1$. For $i=m$, since
$\sprod{Z}_m^\pi=\sprod{Z}_m^{\pi,\theta}=0$,
\begin{equation}\label{eq:main-step-terminal}
I_{m-1}
\leq
(1+C\Delta_m)I_m
+
2\Delta_mR_m.
\end{equation}
Since
\[
I_m
=
\mathbb{E}\big(
    \abs{\xi-\mathcal{C}_{P}^{\overline{\pi}_m}(\xi)}^2
\big),
\]
the discrete Gronwall lemma yields
\begin{equation}\label{eq:maxI}
\max_{0\leq i\leq m}I_i
\leq
C\left\{
    \mathbb{E}\big(
        \abs{\xi-\mathcal{C}_{P}^{\overline{\pi}_m}(\xi)}^2
    \big)
    +
    \sum_{i=1}^m\Delta_iR_i
\right\}.
\end{equation}

For $t\in[t_i,t_{i+1})$,
\[
\Delta^\theta Y_t^\pi
=
\Delta^\theta Y_{t_i}^\pi
+
\int_{t_i}^t\Delta^\theta Z_u^\pi\,dB_u.
\]
Thus, Doob's inequality and It\^o's isometry imply
\[
\mathbb{E}\left(
    \sup_{t\in[t_i,t_{i+1}]}
    \abs{\Delta^\theta Y_t^\pi}^2
\right)
\leq
5\max_{0\leq j\leq m}I_j.
\]
Finally, summing \eqref{eq:main-step} and
\eqref{eq:main-step-terminal}, and using \eqref{eq:maxI} and
$\sum_i\Delta_i=T$, gives
\[
\mathbb{E}\left(
    \int_0^T\abs{\Delta^\theta Z_u^\pi}^2\,du
\right)
\leq
C\left\{
    \mathbb{E}\big(
        \abs{\xi-\mathcal{C}_{P}^{\overline{\pi}_m}(\xi)}^2
    \big)
    +
    \sum_{i=1}^m\Delta_iR_i
\right\}.
\]
Combining the last two estimates proves the result.

\subsubsection*{Proof of Lemma \ref{lemma generator basis truncation}}

Since $Y_{t_i}^{\pi,\theta}$ and $\sprod{Z}_i^{\pi,\theta}$ are
$\mathcal G_i^M$-measurable,
\[
    f_i^{M,\pi,\theta}
    =
    \mathbb E\big(
        f_i^{\pi,\theta}\mid\mathcal G_i^M
    \big).
\]
Let $\Pi^{\overline{\pi}_i}$ denote the orthogonal projection of $\mathbb H$ onto $\operatorname{span}\{h_1^i,\dots,h_{M(i)}^i\}$. Hence, by the same argument as in Lemma~\ref{lemma truncation chaos M},
\[
    \mathbb E\big(
        |f_i^{\pi,\theta}-f_i^{M,\pi,\theta}|^2
    \big)
    \leq
    \mathbb E\big(
        \norm{(I-\Pi^{\overline{\pi}_i})Df_i^{\pi,\theta}}_{\mathbb H}^2
    \big).
\]
By the Malliavin chain rule,
\begin{flalign*}
    D_s f_i^{\pi,\theta}
=
 (D_s f)(t_i,Y_{t_i}^{\pi,\theta},\sprod{Z}_i^{\pi,\theta})
& +
\partial_y f(t_i,Y_{t_i}^{\pi,\theta},\sprod{Z}_i^{\pi,\theta})
D_sY_{t_i}^{\pi,\theta} \\ &
+
\partial_z f(t_i,Y_{t_i}^{\pi,\theta},\sprod{Z}_i^{\pi,\theta})
D_s\sprod{Z}_i^{\pi,\theta}.
\end{flalign*}
Both $Y_{t_i}^{\pi,\theta}$ and $\sprod{Z}_i^{\pi,\theta}$ are
$\mathcal G_i^M$-measurable. Their Malliavin derivatives lie in\linebreak
$\operatorname{span}\{h_1^i,\dots,h_{M(i)}^i\}$. Therefore,
\[
    (I-\Pi^{\overline{\pi}_i})Df_i^{\pi,\theta}
    =
    (I-\Pi^{\overline{\pi}_i})
    \big(
        (Df)(t_i,Y_{t_i}^{\pi,\theta},\sprod{Z}_i^{\pi,\theta})
    \big).
\]
Using Assumption~\ref{main assumption convergence rate}(iii) and arguing as in
Lemma~\ref{lemma approximation holder}, we obtain
\[
    \mathbb E\big(
        \abs{f_i^{\pi,\theta}-f_i^{M,\pi,\theta}}^2
    \big)
    \leq
    C|\overline{\pi}_i|^{2\beta_f}
    \leq
    C M^{-2\beta_f}.
\]
Multiplying by $\Delta_i$ and summing over $i$ concludes the proof.

\subsubsection*{Proof of Lemma \ref{lem:malliavin-energy}}

We first prove \eqref{eq:energy-estimates-Malliavin-Y-Z}. For $i=m,\dots,1$, the recursion of the scheme gives
\begin{gather*}
    Y_{t_{i-1}}^{\pi,\theta}
    +
    \int_{t_{i-1}}^{t_i}Z_r^{\pi,\theta}\,dB_r
    =
    Y_{t_i}^{\pi,\theta}
    +
    \Delta_i\mathcal C_P(f_i^{M,\pi,\theta}).
\end{gather*}
Differentiating in the Malliavin sense for $s<t_{i-1}$ gives
\begin{gather*}
    D_sY_{t_{i-1}}^{\pi,\theta}
    +
    \int_{t_{i-1}}^{t_i}D_sZ_r^{\pi,\theta}\,dB_r
    =
    D_sY_{t_i}^{\pi,\theta}
    +
    \Delta_iD_s\mathcal C_P(f_i^{M,\pi,\theta}),
\end{gather*}
where we used (3.6) on page 57 of \textcite{Nualart_Nualart_2018}. Hence, by
Itô's isometry, Young's inequality and Lemma
\ref{lemma inequality chaos truncation},
\begin{flalign*}
    \mathbb E\Big(
        \norm{DY_{t_{i-1}}^{\pi,\theta}}_{\mathbb H}^2
        +
        \int_0^{t_{i-1}}\int_{t_{i-1}}^{t_i}
        |D_sZ_r^{\pi,\theta}|^2\,dr\,ds
    \Big)
    &\leq
    \Big(1+\frac{\Delta_i}{\epsilon}\Big)
    \mathbb E\big(
        \norm{DY_{t_i}^{\pi,\theta}}_{\mathbb H}^2
    \big)\\
    &\quad+
    \Delta_i(\Delta_i+\epsilon)
    \mathbb E\big(
        \norm{Df_i^{M,\pi,\theta}}_{\mathbb H}^2
    \big).
\end{flalign*}

Moreover, the Malliavin derivative commutes with the conditional expectation defining
$\sprod Z_{i-1}^{\pi,\theta}$, so conditional Jensen, Cauchy--Schwarz and the mesh
regularity give, for $i=2,\dots,m$,
\begin{gather*}
    \frac{\Delta_{i-1}}{L}
    \mathbb E\big(
        \norm{D\sprod Z_{i-1}^{\pi,\theta}}_{\mathbb H}^2
    \big)
    \leq
    \mathbb E\Big(
        \int_0^{t_{i-1}}\int_{t_{i-1}}^{t_i}
        |D_sZ_r^{\pi,\theta}|^2\,dr\,ds
    \Big).
\end{gather*}
Set
\begin{gather*}
    I_i
    :=
    \mathbb E\Big(
        \norm{DY_{t_i}^{\pi,\theta}}_{\mathbb H}^2
        +
        \frac{\Delta_i}{L}
        \norm{D\sprod Z_i^{\pi,\theta}}_{\mathbb H}^2
    \Big),
    \qquad i=1,\dots,m,
\end{gather*}
with $\sprod Z_m^{\pi,\theta}=0$.

We next bound the last term above. Since
\[
    f_i^{M,\pi,\theta}
    =
    \mathbb E\big(
        f_i^{\pi,\theta}\mid\mathcal G_i^M
    \big),
\]
the $\mathbb D^{1,2}$-contraction property of conditional expectation and the
Malliavin chain rule yield
\begin{flalign}
    \mathbb E\big(
        \norm{Df_i^{M,\pi,\theta}}_{\mathbb H}^2
    \big)
    &\leq
    \mathbb E\big(
        \norm{Df_i^{\pi,\theta}}_{\mathbb H}^2
    \big) \nonumber\\
    &\leq
    C\left\{
        1+
        \mathbb E\big(
            \norm{DY_{t_i}^{\pi,\theta}}_{\mathbb H}^2
        \big)
        +
        \mathbb E\big(
            \norm{D\sprod Z_i^{\pi,\theta}}_{\mathbb H}^2
        \big)
    \right\},
    \label{eq:Df-chain-rule-bound}
\end{flalign}
where we used Assumption~\ref{main assumption convergence rate}(iii) to control
the intrinsic Malliavin derivative of $f$, and the boundedness of its derivatives
with respect to $(y,z)$.

Choosing $\epsilon>0$ sufficiently small and then $|\pi|\leq\epsilon$, we obtain
for $i=m,\dots,2$,
\begin{gather}\label{eq:I-recursion-tight}
    I_{i-1}
    +
    \frac{\Delta_i}{2L}
    \mathbb E\big(
        \norm{D\sprod Z_i^{\pi,\theta}}_{\mathbb H}^2
    \big)
    \leq
    (1+C\Delta_i)I_i+C\Delta_i.
\end{gather}
Since
\begin{gather*}
    I_m
    =
    \mathbb E\big(
        \norm{D\mathcal C_P^{\overline{\pi}_m}(\xi)}_{\mathbb H}^2
    \big)
    \leq
    \mathbb E\big(
        \norm{D\xi}_{\mathbb H}^2
    \big),
\end{gather*}
the discrete Gronwall inequality gives
\begin{gather}\label{eq:gronwall-maxI}
    \max_{1\leq i\leq m}I_i
    \leq
    C\left(
        1+
        \mathbb E\big(
            \norm{D\xi}_{\mathbb H}^2
        \big)
    \right).
\end{gather}
Summing \eqref{eq:I-recursion-tight} and using
$\Delta_1\mathbb E(\norm{D\sprod Z_1^{\pi,\theta}}_{\mathbb H}^2)
\leq LI_1$ therefore yields
\begin{gather*}
    \sum_{i=1}^m\Delta_i
    \mathbb E\big(
        \norm{D\sprod Z_i^{\pi,\theta}}_{\mathbb H}^2
    \big)
    \leq
    C\left(
        1+
        \mathbb E\big(
            \norm{D\xi}_{\mathbb H}^2
        \big)
    \right),
\end{gather*}
which, together with \eqref{eq:gronwall-maxI}, proves
\eqref{eq:energy-estimates-Malliavin-Y-Z}.

Finally, multiplying \eqref{eq:Df-chain-rule-bound} by $\Delta_i$, summing over
$i$, and using \eqref{eq:energy-estimates-Malliavin-Y-Z}, we obtain
\begin{gather*}
    \sum_{i=1}^m\Delta_i
    \mathbb E\big(
        \norm{Df_i^{M,\pi,\theta}}_{\mathbb H}^2
    \big)
    \leq
    C\left(
        1+
        \mathbb E\big(
            \norm{D\xi}_{\mathbb H}^2
        \big)
    \right),
\end{gather*}
which proves \eqref{eq:energy-estimates-malliavin-f}.

\subsection{Additional results and proofs of Subsection \ref{subsection: error montecarlo}}

\begin{lemma}\label{lemma MC}
We have
\begin{gather}\label{second equality lemma}
    \overline{\mathbb E}\Big(
        \abs{
        \widehat{\mathcal C}_{P}^{\overline{\pi}_m}(\xi)
        -
        \mathcal C_{P}^{\overline{\pi}_m}(\xi)
        }^2
    \Big)
    =
    \frac{1}{N}\mathbf V_m(\xi),
\end{gather}
and, for $i=1,\dots,m$,
\begin{gather}\label{first equality lemma}
    \overline{\mathbb E}\Big(
        \abs{
        \widehat{\mathcal C}_{P}^{\overline{\pi}_i}(\widehat F_i^{\pi,\theta})
        -
        \mathcal C_{P}^{\overline{\pi}_i}(F_i^{\pi,\theta})
        }^2
    \Big)
    \leq
    \overline{\mathbb E}\big(
        \abs{\widehat F_i^{\pi,\theta}-F_i^{\pi,\theta}}^2
    \big)
    +
    \frac{\Delta_i^2}{N}
    \mathbf V_i(\widehat f_i^{\pi,\theta}).
\end{gather}
\end{lemma}

\begin{proof}
We begin by proving \eqref{second equality lemma}. Using that
$\{\sqrt{a!}\,H_a^{(m,0)}\}_{a\in\mathcal A_{M(m)}}$ is orthonormal and that $\hat d_a^\xi$ is independent of $B^{m,0}=B$,
\begin{flalign*}
    \overline{\mathbb E}\Big(
        \abs{
        \widehat{\mathcal C}_{P}^{\overline{\pi}_m}(\xi)
        -
        \mathcal C_{P}^{\overline{\pi}_m}(\xi)
        }^2
    \Big)
    =
    \sum_{k=0}^P
    \sum_{\substack{a\in\mathcal A_{M(m)}\\|a|=k}}
    \frac{1}{a!}
    \overline{\mathbb E}\big(
        \abs{\hat d_a^\xi-d_a^\xi}^2
    \big),
\end{flalign*}
where
$d_a^\xi\coloneqq a!\,\overline{\mathbb E}\big(\xi\times H_a^{(m,0)}\big)$.
By the definition of $\hat d_a^\xi$ and the i.i.d. property of the samples,
\begin{gather*}
    \frac{1}{a!}
    \overline{\mathbb E}\big(
        \abs{\hat d_a^\xi-d_a^\xi}^2
    \big)
    =
    \frac{a!}{N}
    \overline{\mathbb V}\big(
        \xi\times H_a^{(m,0)}
    \big).
\end{gather*}
Summing over $a$ proves \eqref{second equality lemma}.

We now prove \eqref{first equality lemma}. Again, by orthonormality and the independence of $\hat d_a^i$ from $B^{i,0}=B$,
\begin{flalign*}
    \overline{\mathbb E}\Big(
        \abs{
        \widehat{\mathcal C}_{P}^{\overline{\pi}_i}(\widehat F_i^{\pi,\theta})
        -
        \mathcal C_{P}^{\overline{\pi}_i}(F_i^{\pi,\theta})
        }^2
    \Big)
    =
    \sum_{k=0}^P
    \sum_{\substack{a\in\mathcal A_{M(i)}\\|a|=k}}
    \frac{1}{a!}
    \overline{\mathbb E}\big(
        \abs{\hat d_a^i-d_a^i}^2
    \big).
\end{flalign*}
Using orthogonality of conditional expectation,
\begin{gather*}
    \overline{\mathbb E}\big(
        \abs{\hat d_a^i-d_a^i}^2
    \big)
    =
    \overline{\mathbb E}\big(
        \abs{\hat d_a^i-\overline{\mathbb E}(\hat d_a^i\mid\mathcal F^{i+1\leq})}^2
    \big)
    +
    \overline{\mathbb E}\big(
        \abs{\overline{\mathbb E}(\hat d_a^i\mid\mathcal F^{i+1\leq})-d_a^i}^2
    \big).
\end{gather*}
Conditionally on $\mathcal F^{i+1\leq}$, the random variables
$\widehat f_i^{\pi,\theta,n}H_a^{(i,n)}$, $n=1,\dots,N$, are i.i.d., since the pairs $(B^{i,n},f^{i,n})$ are independent copies of $(B,f)$. Hence,
\begin{flalign*}
    \frac{1}{a!}
    \overline{\mathbb E}\big(
        \abs{\hat d_a^i-\overline{\mathbb E}(\hat d_a^i\mid\mathcal F^{i+1\leq})}^2
    \big)
    &=
    \Delta_i^2\frac{a!}{N}
    \overline{\mathbb E}\Big[
        \overline{\mathbb V}\big(
            \widehat f_i^{\pi,\theta}H_a^{(i,0)}
            \mid\mathcal F^{i+1\leq}
        \big)
    \Big]\\
    &\leq
    \Delta_i^2\frac{a!}{N}
    \overline{\mathbb V}\big(
        \widehat f_i^{\pi,\theta}H_a^{(i,0)}
    \big).
\end{flalign*}
Moreover,
\begin{gather*}
    \overline{\mathbb E}(\hat d_a^i\mid\mathcal F^{i+1\leq})
    =
    a!\,
    \overline{\mathbb E}\Big(
        \widehat F_i^{\pi,\theta}H_a^{(i,0)}
        \,\big|\,
        \mathcal F^{i+1\leq}
    \Big).
\end{gather*}
Therefore, by orthonormality and the non-expansiveness of the chaos projection,
\begin{flalign*}
    \sum_{k=0}^P
    \sum_{\substack{a\in\mathcal A_{M(i)}\\|a|=k}}
    \frac{1}{a!}
    \overline{\mathbb E}\big(
        \abs{d_a^i-\overline{\mathbb E}(\hat d_a^i\mid\mathcal F^{i+1\leq})}^2
    \big)
    \leq
    \overline{\mathbb E}\big(
        \abs{F_i^{\pi,\theta}-\widehat F_i^{\pi,\theta}}^2
    \big).
\end{flalign*}
Combining the previous estimates and summing over $a$ proves
\eqref{first equality lemma}.
\end{proof}

\subsubsection*{Proof of Lemma \ref{lemma error implementation MC}}

All $L^p$-norms appearing in this proof are taken with respect to $\overline{\mathbb P}$.

For $i=1,\dots,m$, set
\[
    \mathcal K_i
    \coloneqq
    \sum_{k=0}^P
    \sum_{\substack{a\in\mathcal A_{M(i)}\\|a|=k}}
    a!\abs{H_a^{(i,0)}}^2.
\]
Since $\sqrt{a!}\,H_a^{(i,0)}$ belongs to the chaos of order $|a|$ and has unit $L^2$-norm, Gaussian hypercontractivity gives
\[
    \norm{\mathcal K_i}_{L^2}
    \leq
    \sum_{k=0}^P
    \sum_{\substack{a\in\mathcal A_{M(i)}\\|a|=k}}
    \norm{\sqrt{a!}\,H_a^{(i,0)}}_{L^4}^2
    \leq
    3^P\binom{P+M}{M}.
\]
Consequently, if $G$ is, conditionally on $\mathcal F^{i+1\leq}$, a finite chaos decomposition of order at most $P$ with respect to $B^{i,0}$, then conditional Cauchy--Schwarz and hypercontractivity yield
\begin{gather}\label{weighted finite chaos estimate}
    \overline{\mathbb E}\big(\abs{G}^2\mathcal K_i\big)
    \leq
    \Lambda_{P,M}\overline{\mathbb E}(\abs{G}^2).
\end{gather}

Assumption~\ref{main assumption convergence rate}(ii),(iv) implies
$\sup_{t\in[0,T]}\norm{f(t,0,0)}_{L^q}<\infty$. Since $q>4$, Hölder's inequality gives
\[
    \overline{\mathbb E}
    \big(
        \abs{f(t_i,0,0)}^2\mathcal K_i
    \big)
    \leq
    \norm{f(t_i,0,0)}_{L^q}^2
    \norm{\mathcal K_i}_{L^{q/(q-2)}}
    \leq
    C\Lambda_{P,M}.
\]
The random variables $\widehat Y_{t_i}^{\pi,\theta}$ and
$\sprod{\widehat Z}_i^{\pi,\theta}$ satisfy the condition of
\eqref{weighted finite chaos estimate}. Therefore, by the Lipschitz continuity of $f$,
\begin{gather}\label{bound V unbounded}
    \mathbf V_i(\widehat f_i^{\pi,\theta})
    \leq
    C\Lambda_{P,M}
    \left(
        1
        +
        \overline{\mathbb E}\abs{\widehat Y_{t_i}^{\pi,\theta}}^2
        +
        \overline{\mathbb E}\abs{\sprod{\widehat Z}_i^{\pi,\theta}}^2
    \right).
\end{gather}
The same argument and Hölder's inequality show that
\[
    \mathbf V_m(\xi)
    \leq
    \overline{\mathbb E}\big(\abs{\xi}^2\mathcal K_m\big)
    \leq
    3^P\binom{P+M}{M}\norm{\xi}_{L^q}^2.
\]

The standard stability estimate for the exact-coefficient scheme, using the contraction property of the chaos projections and the Lipschitz continuity of $f$, gives
\begin{gather}\label{stability exact chaos scheme MC}
    \max_{1\leq i\leq m}
    \mathbb E\abs{Y_{t_i}^{\pi,\theta}}^2
    +
    \mathbb E\left(
        \int_0^T\abs{Z_u^{\pi,\theta}}^2\,du
    \right)
    \leq C.
\end{gather}
In particular, conditional Jensen and the mesh regularity imply
\begin{gather}\label{weighted stability exact chaos scheme MC}
    \sum_{i=1}^m
    \Delta_i^2
    \left(
        1
        +
        \mathbb E\abs{Y_{t_i}^{\pi,\theta}}^2
        +
        \mathbb E\abs{\sprod Z_i^{\pi,\theta}}^2
    \right)
    \leq
    C\abs{\pi}.
\end{gather}

Set
\[
    \widehat\Delta Y_{t_i}^{\pi,\theta}
    \coloneqq
    Y_{t_i}^{\pi,\theta}
    -
    \widehat Y_{t_i}^{\pi,\theta},
    \qquad
    \widehat\Delta Z_u^{\pi,\theta}
    \coloneqq
    Z_u^{\pi,\theta}
    -
    \widehat Z_u^{\pi,\theta},
\]
and
\[
    J_i
    \coloneqq
    \overline{\mathbb E}\left(
        \int_{t_i}^{t_{i+1}}
        \abs{\widehat\Delta Z_u^{\pi,\theta}}^2\,du
    \right),
    \qquad
    I_i
    \coloneqq
    \overline{\mathbb E}
    \abs{\widehat\Delta Y_{t_i}^{\pi,\theta}}^2
    +
    J_i,
\]
for $i=0,\dots,m-1$, with
$I_m\coloneqq\overline{\mathbb E}
\abs{\widehat\Delta Y_{t_m}^{\pi,\theta}}^2$.

By It\^o's isometry, Lemma~\ref{lemma MC}, the Lipschitz continuity of $f$, and the mesh regularity, for $|\pi|$ sufficiently small,
\begin{gather}\label{MC recursion unbounded}
    I_{i-1}
    +
    \frac12J_i
    \leq
    (1+C\Delta_i)I_i
    +
    \frac{\Delta_i^2}{N}
    \mathbf V_i(\widehat f_i^{\pi,\theta}),
    \qquad i=1,\dots,m-1,
\end{gather}
while the analogous inequality without $J_m$ holds for $i=m$.

Moreover,
\[
    \overline{\mathbb E}
    \abs{\widehat Y_{t_i}^{\pi,\theta}}^2
    \leq
    2\mathbb E\abs{Y_{t_i}^{\pi,\theta}}^2
    +
    2I_i,
\]
and, for $i=1,\dots,m-1$,
\[
    \overline{\mathbb E}
    \abs{\sprod{\widehat Z}_i^{\pi,\theta}}^2
    \leq
    2\mathbb E\abs{\sprod Z_i^{\pi,\theta}}^2
    +
    \frac{2I_i}{\Delta_{i+1}}.
\]
Combining these estimates with \eqref{bound V unbounded}, using
$\Delta_i/\Delta_{i+1}\leq L$ and $\abs{\pi}\leq1$, gives
\[
    I_{i-1}
    +
    \frac12J_i
    \leq
    \left(
        1+C\Delta_i+
        \frac{C\Lambda_{P,M}}{N}\Delta_i
    \right)I_i
    +
    \frac{C\Lambda_{P,M}}{N}
    \Delta_i^2
    \left(
        1+
        \mathbb E\abs{Y_{t_i}^{\pi,\theta}}^2+
        \mathbb E\abs{\sprod Z_i^{\pi,\theta}}^2
    \right).
\]
Since Lemma~\ref{lemma MC} gives
$I_m=\mathbf V_m(\xi)/N$, the discrete Gronwall inequality and
\eqref{weighted stability exact chaos scheme MC} yield
\[
    \max_{0\leq i\leq m}I_i
    \leq
    \frac{C}{N}
    \exp\left(
        \frac{C\Lambda_{P,M}}{N}
    \right)
    \left[
        \mathbf V_m(\xi)
        +
        \Lambda_{P,M}\abs{\pi}
    \right].
\]
Doob's inequality gives the corresponding estimate for the supremum of the $Y$-error on each interval. Finally, summing \eqref{MC recursion unbounded} and using the previous bound gives the same estimate for the integral of the $Z$-error. This proves \eqref{error given by MC}.

\printbibliography[title={References}]

\end{document}